\documentclass[preprint,11pt]{article}
\usepackage{amssymb}
\usepackage{graphicx}
\usepackage{subfig}
\usepackage{amsmath}
\usepackage[mathscr]{euscript}
\usepackage{lscape}
\usepackage{mathtools}
\usepackage{xcolor}
\usepackage{footnote}
\usepackage{float}
\usepackage{array,multirow}
\usepackage{adjustbox}

\DeclareMathAlphabet\mathbfcal{OMS}{cmsy}{b}{n}
\newcommand{\RN}[1]{\textit{\uppercase\expandafter{\romannumeral#1}}}
\begin{document}

\title{\vskip-2cm A general higher-order shell theory for compressible isotropic hyperelastic materials using orthonormal moving frame}
\author{A Arbind\footnote{archana.iitg@gmail.com \newline Preprint submitted to a journal}, J N Reddy, and A R Srinivasa\\
Advanced Computational Mechanics Laboratory\\
J. Mike Walker '66 Department of Mechanical Engineering\\
Texas A{\&}M University, College Station.\\
}
\date{\today}
\maketitle

\begin{abstract}
\noindent  The primary objective of this study is three-fold: (i) to present a general higher-order shell theory to analyze large deformations of thin or thick shell structures made of general compressible hyperelastic materials; (ii) to formulate an efficient shell theory using the orthonormal moving frame, and (iii) to develop and apply the nonlinear weak-form Galerkin finite element model for the proposed shell theory. The displacement field of the line normal to the shell reference surface is approximated by the Taylor series/Legendre polynomials in the thickness coordinate of the shell. The use of an orthonormal moving frame makes it possible to represent kinematic quantities (e.g., the determinant of the deformation gradient) in a far more efficient manner compared to the non-orthogonal covariant bases. Kinematic quantities for the shell deformation are obtained in a novel way in the surface coordinate described in the appendix of this study with the help of exterior calculus. Further, the governing equation of the shell deformation has been derived in the general surface coordinates. To obtain the nonlinear solution in the quasi-static cases, we develop the weak-form finite element model in which the reference surface of the shell is modeled exactly. The general invariant based compressible hyperelastic material model is considered. The formulation presented herein can be specialized for various other nonlinear compressible hyperelastic constitutive models, for example, in bio-mechanics and other soft-material problems (e.g., compressible neo-Hookean material, compressible Mooney-Rivlin material, Saint-Venant Kirchhoff model, and others). Various numerical examples are presented to verify and validate the formulation presented in this study. The scope of potential extensions are outlined in the final section of this study.
\end{abstract}
\vskip2mm \noindent {\it Keywords.}~ higher-order shell theory; thin  and thick shell structures; compressible isotropic hyperelastic material; curved tubular shells; orthonormal moving frames; Cartan's moving frame.

\section{Introduction}

The governing equations for curved shells and their numerical implementation have long been a challenging task in the mechanics community since the first general derivation of the theory was posited by Love \cite{love1888xvi}. Earlier attempts at shell theories (See \cite{naghdi1973theory}, \cite{ventsel2002thin} and references therein.) are based on Gauss's theory of surfaces \cite{gauss1902general} in which the kinematics of deformation is measured via metric and curvature variables (in the first and second fundamental differential forms) of the mid surface of the shell structure in the course of deformation. Further, following the works of Green and Zerna \cite{green1971theoretical}, the equations of motion for shells have been developed using convected coordinates (\cite{naghdi1973theory}). Most of these earlier works in the shell theory were carried out  for linear elastic materials.

 In recent years,  shells made of nonlinear hyper elastic materials have received a great deal of attention  in view of potential applications to soft and bio-materials.
Basar and Ding \cite{basar1996finite} developed a large strain shell model for thin shell structure on the basis of a quadratic displacement approximation in
thickness coordinate by neglecting transverse shear strains for incompressible hyperelastic materials. This leads to a three-parameter theory governed by mid-surface displacements.
Başar and Itskov \cite{bacsar1998finite} presented a thin shell theory for the Ogden material model for rubber‐like shells. They transformed the  strain energy density, which is function of principal stretches, in terms of the invariants of right Cauchy–Green tensor and thus bypassing the need for eigenvalue calculation in their formulation. They have also presented an algorithm to deal with  eigenvalue coalescence for the stretches.
Campello {\it et. al} \cite{campello2011exact} presented a nonlinear shell dynamics for thin shell assuming the Rodrigues rotation vector for the rotation vector field of linear material and neo-Hookean material.
Kiendl {\it et. al} \cite{kiendl2015isogeometric} presented Kirchhoff–Love shell formulations for thin shells for general compressible and incompressible hyperelastic materials using isogeometric analysis. They have used linear thickness strain along with the  Kirchhoff–Love hypothesis and statically condensed the thickness deformation to express the shell equations in terms of metric and curvature of the mid surface of the shell; they have applied the thin shell theory to dynamic simulation of a bio-prosthetic heart valve.
Luo {\it et. al} \cite{luo2016nonlinear} carried out  a nonlinear static and dynamic analysis for hyperelastic thin shells via the absolute nodal coordinate formulation considering the  Kirchhoff–Love hypothesis.
Betsch {\it et. al} developed a shell element for large deformation based on an extensible director approach for compressible and incompressible hyperelastic material. This shell theory accounts for constant thickness stretch through the  thickness via director stretch along with large rotation.
Song and Dai \cite{song2016consistent} have developed consistent models for thin shells via high-order expansion coefficients from the 3-D equations for compressible hyperelastic  cylindrical and spherical shell structures. Further, Li {\it et. al} \cite{li2019consistent} extended this study to incompressible materials.

While much work has been done for hyperelastic thin shells, higher-order shell theories for hyperelastic materials have received much less attention. Such  higher-order shell theories can model transverse normal and shear strain components via a higher-order displacement field approximation. In the case of soft material, such as rubber or biological materials, higher-order shear and transverse deformation shell theories are of considerable  importance as soft shells undergo considerable thickness deformation under loading. In addition, many biological shells are multilayered and so require nonlinear interpolations through the thickness.

Reddy and coworkers have developed several higher-order theories and the corresponding finite element models for large deformation in shell structure for linear material, functionally graded material and laminated composites (see \cite{amabili2020}). Arciniega and Reddy (see \cite{arciniega2007tensor}, \cite{arciniega2006tensor}) formulated a tensor-based 7-parameter shell theory and its finite element model using 3D linear constitutive relation between the Second Piola-Kirchhoff (PK) Stress and the Green-St. Venant strain (St. Venant-Kirchhoff nonlinear material model) with fully nonlinear geometry. Amabili and Reddy \cite{amabili2010new} developed a consistent higher-order shell theory for von Kármán nonlinearity. Amabili \cite{amabili2015non} developed a geometrically non-linear shell theory which allowed third order thickness and shear deformation using 8-parameter displacement field and concluded that such theory for linear material can predict the thickness deformation correctly. He has calculated the solution in terms of Fourier bases. Rivera and Reddy developed 7-parameter and 12-parameter shell theories (see \cite{rivera2016stress}, \cite{rivera2016new}, \cite{rivera2017nonlinear}) for functionally graded material and laminated composite shell structures using St. Venant-Kirchhoff constitutive model. Rivera, Reddy, and Amabili \cite{rivera2020continuum} developed a new 8-parameter shell theory for St. Venant-Kirchhoff constitutive model, which allows the use of a third‐order thickness stretch kinematics, which avoids Poisson's locking.
Amabili, Breslavsky, and Reddy \cite{amabili2019nonlinear} developed a 9-parameter shell theory for circular cylindrical shell considering incompressible neo-Hookean material. This shell theory is higher-order in both shear and thickness deformations where the four parameters describing the thickness deformation are obtained directly from the incompressibility condition.

We note that unlike the case of thin shells, the approaches for higher-order shell modeling have been based on specific geometries (cylindrical, spherical, etc.) and/or special constitutive models.
For a general hyperelastic material model, no such higher order shell theory is reported in the literature which can account for large deformation as well as shear and thickness deformation. This is partly because of the complexity of the formulation for a general hyperelastic material and a general shell geometry is a daunting task with extremely complex formulations due to the curvilinear coordinates and non-orthonormal base vectors.

In this study, we derive such a general higher order shell theory which is based on a general polynomial expansion of the transverse and in-plane displacement components with a general order of approximation. Rather than obtaining the equation through asymptotic expansions or integration of the 3-D equations, we derive the equations and boundary conditions for the FEM formulation directly using the principal of virtual work.

Undoubtedly, such generality in the nonlinear material model and deformation approximation of the proposed shell theory will invite enormous complexity in terms of kinematics and numerical model formulation in the general curvilinear coordinate system of the curved surface. The natural covariant bases for general curvilinear coordinates, in general, are non-orthogonal. Due to the presence of non-orthonormal basis vectors, the calculation of even the simplest kinematic quantities such as deformation gradients require extensive computations involving the non-identity metric tensor (see Table\footnote{Here the repeated index implies summation over the range of the index, namely, 1 to 3 for three-dimensional space. Also, ${\bf e}_i$ and ${\bf e}^i$ are the covariant and contravariant basis vectors of the natural covariant frame, respectively; $g_{ij}$ are the components of the covariant metric tensor; $\hat{\bf e}_i$ are the orthonormal bases of orthonormal (Cartan's) moving frame; and $e^{ijk}$ or $e_{ijk}$ is the permutation symbol.}.~\ref{Table:1}). The complexity in the numerical model of such higher-order shell theories can be fathomed by the fact that the virtual internal energy contains approximately 20000 terms for the 7-parameter shell theory (see Chapter 9, Reddy \cite{reddy2015introduction}) even for St. Venant-Kirchhoff constitutive model, which is considered to be the simplest nonlinear hyperelastic material model. The complexity of numerical formulation gets compounded when one considers the {\it general invariant based hyperelastic} nonlinear material model, which involves terms, such as, determinants and other invariants of right Cauchy-Green deformation tensor.
\begin{table}
\centering
\caption{Comparison for the kinematic quantities or tensor operations for covariant frame and orthonormal frame.}\label{Table:1}
\begin{tabular}{ ||c||c|| }
 \hline
 Natural covariant frame & Orthonormal moving frame \\[2ex]
  \hline\hline
${\bf A}=A^{ij}{\bf e}_i{\bf e}_j=A_{ij}{\bf e}^i{\bf e}^j
  =A_{i}^{.j}{\bf e}^i{\bf e}_j$
& ${\bf A}=A_{ij}\hat{\bf e}_i\hat{\bf e}_j$  \\[2ex]
$
  {\bf A}\cdot{\bf B}=A^{ij}{\bf e}_i{\bf e}_j\cdot B^{kl}{\bf e}_k{\bf e}_l
  =A^{ij}B^{kl}g_{jk}{\bf e}_i{\bf e}_l
$ & ${\bf A}\cdot{\bf B}=A_{ij}B_{jk}\hat{\bf e}_i\hat{\bf e}_j$  \\[2ex]
$
\text{det}({\bf A}) = e^{ijk}A^1_iA^2_jA^3_k
=e^{ijk}A^{1l}A^{2m}A^{3n}g_{li}g_{mj}g_{nk}
$ & $\text{det}({\bf A}) = e_{ijk}A_{1i}A_{2j}A_{3k}$\\[2ex]
 \hline
\end{tabular}
\end{table}

To circumvent the difficulty of the non-orthogonal bases in the curvilinear coordinate system, in this study, we adopt the orthonormal (Cartan's) moving frame and derive all the kinematic invariants such as deformation gradient, determinant of displacement gradient, subsequently, the governing equation and the numerical model. In the shell theory literature, Knowles and Reissner \cite{knowles1956derivation} have derived the shell theory in orthonormal bases by restricting the surface coordinates lines (or coordinates) as orthogonal lines (or coordinates). For example, $(\theta,s)$ coordinate of surfaces of revolution (see Appendix A) are examples of orthogonal coordinates. However, in the best of the author's knowledge, there is no study reported in the literature, which derives the governing equations of shells and its numerical model with general non-orthogonal coordinates via the orthonormal moving frame. Appendix A of this study also presents a novel way\footnote{The method is slightly different from than process involving the structure matrix proposed by Darboux and Cartan for doing calculus with an orthonormal moving frame. However, one can find the similarity, on careful observation, between the two approaches as both reach the same results.}, in conjunction with exterior calculus, to derive the kinematics of the deformation considering the orthonormal basis. The derivation presents the kinematics for a general curved tubular surface with a varying radius which encompasses a wide variety of curved shell surfaces, spherical shell, and also for plates.

\begin{figure}[htbp]
\centering{\includegraphics[trim=4cm 6cm 5cm 6cm,clip,scale=0.7]{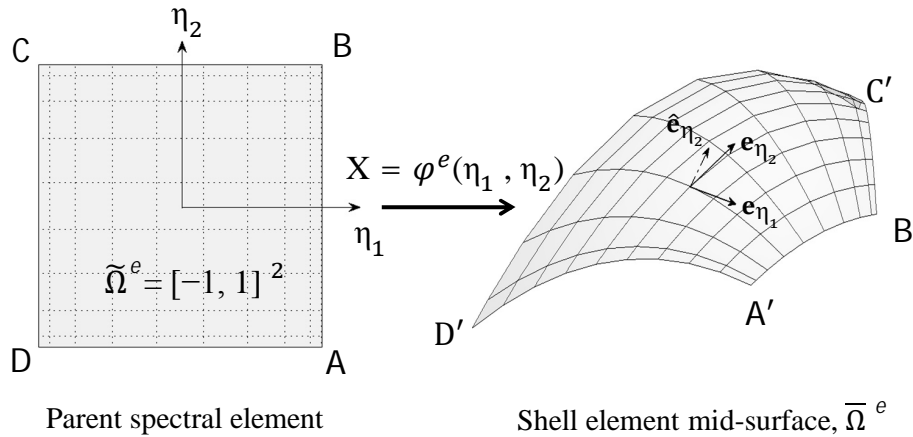}}
\caption{Isoparametric mapping of the shell element mid surface from the parent element. Here, we see that the surface coordinate lines in the physical shell element surface are not orthogonal and the bases (${\bf e}_{\eta_1},{\bf e}_{\eta_2}$) are not orthonormal. So, we generate ($\hat{\bf e}_{\eta_1},\hat{\bf e}_{\eta_2}$) as orthonormal moving frame on each shell surface element where $\hat{\bf e}_{\eta_1}$ is the unite vector along ${\bf e}_{\eta_1}$}.
\label{fig:1}
\end{figure}

The non-orthogonal covariant coordinate bases arise very often in various cases, such as arbitrarily curved tubular shells or the computational shell theory (see \cite{reddy2015introduction}, \cite{arciniega2007tensor}). In the later, we model the mid-surface of the shell element using the parent element via an isoparametric map, as shown in the Fig.\ref{fig:1}., where the coordinates $(\eta_1, \eta_2)$ of the parent spectral element becomes the surface coordinate of the shell element in the physical space. Also, we note that these coordinate lines are not always orthogonal to each other, in general. Hence, in the general case, the surface basis vectors generated by these coordinates are not orthogonal. As stated before, the computation in non-orthogonal bases is a tedious task due to the involvement of metric tensor, hence the use of the orthogonal moving frame on each shell element, separately, would reduce the complexity of the tensor operation in the formulations in the finite element model.

This theory can be specialized to specific higher-order shell theory, using specific values of order of approximation of displacement field, to different higher order shell theory as demanded by the problem at hand.

Furthermore, in contrast to  the existing computational shell theory, in this study, the geometry of the shell structure has been modeled exactly (see Appendix A for the geometries covered in this study). The components of displacement field along the orthogonal surface bases are approximated using polynomial expansion, in contrast to the existing computational shell theory where the components of displacement field along the global Euclidean bases are approximated. This further reduces the finite element model's complexity as no transformation is needed from Euclidean to the surface coordinates for each element at each nonlinear iteration. In this study, such transformations are only required to plot the final deformed shape of the shell surface. However, the proposed approach is limited to surfaces with continuous orthogonal bases. Moreover, there exist many applications where the elemental surface interpolation is imperative (see Fig.\ref{fig:1}), and we can extend the current study for approximated elemental shell surface geometry with the approximation of displacement component in global Euclidean coordinate system.

The outline of this study is as follows. We first present a brief review of the hyperelastic material model and various stress measures for the material model along with the governing equations in Section 2. Then in Section 3 we introduce the kinematics, governing equation for General higher-order shell theory in general curvilinear surface coordinates considering the orthonormal moving frame. In the appendix of this study, the kinematics specified to any arbitrary curved tubular shell with the variable radius (which can also be specialized to the surface of revolution) along with spherical shells are developed with a novel approach. In Section 4, we present the weak form finite element model of the introduced higher-order shell theory. In Section 5, we specify the parameters needed for various nonlinear hyperelastic material model and the finite element analysis, followed by several numerical examples. In the Section 6, we summarize the present work and present conclusions of this study.

\section{Hyperelastic or Green elastic material}
The hyperelastic material, also called Green elastic material, is a class of elastic material where the stress tensor at any point can be derived from the strain energy stored in the deformed body. The strain energy density functional of such material (see \cite{spencer2004continuum}) is given as follows:
\begin{eqnarray}\label{eq:1}
\psi = \hat{\psi}(I_C,II_C,III_C)
\end{eqnarray}
where $I_C,II_C$ and $III_C$ are the first, second, and third principal invariants of the right Cauchy--Green deformation tensor, ${\bf C}={\bf F}^{\rm T}{\bf F}$ defined as:
\begin{eqnarray}\label{eq:2}
I_C = \text{tr}({\bf C}),\quad II_C={1\over 2}\left((\text{tr}({\bf C}))^2-\text{tr}({\bf C}^2)\right),\quad III_C=\text{det}({\bf C})=J^2
\end{eqnarray}
where $J$ is the determinant of the deformation tensor ${\bf F}$. The strain energy density function can also be expressed in terms of other mutually independent invariants of ${\bf C}$ such as $(I_1,I_2,I_3)$, which are defined as follows:
\begin{eqnarray}\label{eq:3}
I_1= I_C =\text{tr}({\bf C}) ,\quad I_2 = \text{tr}({\bf C}^2),\quad I_3=\text{tr}({\bf C}^3)
\end{eqnarray}
Let us express the strain energy density functional for hyperelastic material as follows:
\begin{eqnarray}\label{eq:4}
\psi = \psi(I_1,I_2,J)
\end{eqnarray}
The derivative of the invariants of ${\bf C}$ with respect to the deformation tensor ${\bf F}$ are
\begin{eqnarray}\label{eq:5}
\frac{dI_1}{d{\bf F}} = 2{\bf F},\quad \frac{dI_2}{d{\bf F}} = 4{\bf F}{\bf C},\quad \frac{dJ}{d{\bf F}} = J{\bf F}^{\rm -T}
\end{eqnarray}
The governing equation for the compressible hyperelastic material body for the static case is given by
\begin{eqnarray}\label{eq:6}
-\text{Div}({\bf P})=\rho_0 {\bf b}
\end{eqnarray}
where ${\bf b}$ is the body force vector measured per unit mass of the body and ${\bf P}$ is the first Piola--Kirchhoff stress tensor
\begin{equation}\label{eq:7}
{\bf P} =  J\,\left(\frac{\partial \psi}{\partial J}\right) {\bf F}^{-\rm T}+
2\frac{\partial \psi} {\partial I_1}{\bf F}+4 \frac{\partial \psi}{ \partial I_2}{\bf F}{\bf C}
\end{equation}
The Cauchy stress tensor, ${\boldsymbol \sigma}$ is
\begin{align}\label{eq:8}
{\boldsymbol \sigma}= J^{-1} {\bf P} {\bf F}^{\rm T}
= \beta_3{\bf I} +\frac{1}{J}\beta_1{\bf B}+\frac{1}{J}\beta_2{\bf B}^2
\end{align}
where
\begin{eqnarray}\label{eq:9}
&&\beta_1=2 \frac{\partial \psi}{\partial I_1},\quad \beta_2=4 \frac{\partial \psi} {\partial I_2},\quad \beta_3= \frac{\partial \psi}{\partial J}
\end{eqnarray}
and ${\bf B}={\bf F}{\bf F}^{\rm T}$ is the left  Cauchy--Green deformation tensor. To ensure zero stress in the natural configuration, we require
\begin{equation}\label{eq:10}
\beta_1+\beta_2+\beta_3=0
\end{equation}
Also, the boundary condition is given as:
\begin{eqnarray}\label{eq:11}
[{\bf P}]{\bf N}={\bf q}
\end{eqnarray}
where ${\bf N}$ is the unit vector normal to the boundary surface in the reference configuration and ${\bf q}$ is the surface traction (transformed (or pulled) back to the reference configuration) acting on the boundary surface of the body.
\section{General higher order shell theory}
In this section, we derive the governing equation for the arbitrary shell structure in the general surface coordinates considering Cartan's moving frame.
\subsection{Curvilinear coordinate system and orthonormal moving frame}
\begin{figure}[htbp]
\centering{\includegraphics[trim=0cm 3cm 0cm 2cm,clip,scale=0.50]{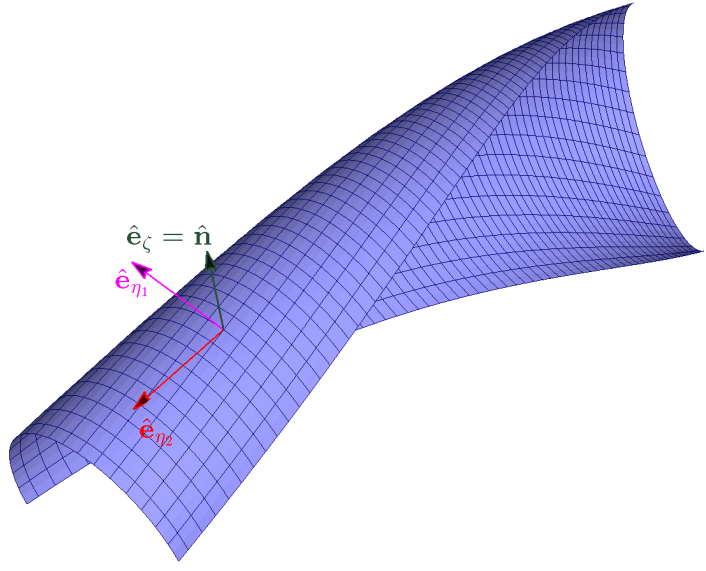}}
\caption{Arbitrary curved shell surface with Cartan's moving frame. Here, we note that the surface bases are not necessarily along the coordinate lines, contrary to the natural covariant frame.}
\label{fig:2}
\end{figure}

Let $(\eta_1,\eta_2,\zeta)$ constitute a curvilinear coordinate system for the shell structure in three-dimensional space, where $(\eta_1,\eta_2)$ are the surface coordinates on the reference-surface\footnote{Generally, the mid surface of the shell structure could be considered as the reference surface of the shell structure. However, this is not necessary.} of the shell structure, whereas $\zeta$ is the thickness coordinate measured along the normal direction to the shell reference surface. We will take $\zeta$-coordinate as zero at the reference surface of the shell structure. For any general surface coordinates $(\eta_1,\eta_2)$, the corresponding covariant basis may or may not be orthogonal. So, in this study, we will consider a {\it non-coordinate} orthonormal moving frame (see Fig.~\ref{fig:2}) for the ease of the derivation of the governing equation and finite element model, as it would become evident in the following sections.
Now, let ${\bf R}$ be the position vector of any arbitrary point ${\cal P}$; then at that point we define the following orthonormal basis\footnote{If we assume the surface coordinate $\eta_1$ as the arc length coordinate of the $\eta_1$ coordinate line, then the resulting moving frame would become the Darboux frame (see \cite{darboux1896leccons}), which is a surface analog of the Serre-Frenet frame of the space curve.}, which may or may not align with the covariant bases for the surface and normal coordinates assumed:
\begin{eqnarray}\label{eq:12}
\hat{\bf e}_{\zeta}=\hat{\bf n} ={{\bf R}_{,\zeta}\over||{\bf R}_{,\zeta}||},\quad
\hat{\bf e}_{\eta_1} ={{\bf R}_{,\eta_1}\over||{\bf R}_{,\eta_1}||},\quad \hbox{and}\quad \hat{\bf e}_{\eta_2}=\hat{\bf e}_{\zeta}\times\hat{\bf e}_{\eta_1}
\end{eqnarray}
where ${\bf R}_{,\zeta}$, for example, represents the derivative of ${\bf R}$ with respect to $\zeta$. The set $(\hat{\bf e}_{\eta_1},\hat{\bf e}_{\eta_2},\hat{\bf e}_{\zeta})$ forms a right-handed orthonormal basis, which does not necessarily align with the covariant bases; when these orthonormal bases do not align with coordinate bases (or covariant bases) then they form a {\it non-coordinate} orthonormal moving frame, for example, in the case of general closed curved pipe surface (see appendix A).
\subsection{Displacement field}
Now, in the orthonormal coordinate system, we approximate the displacement field of apoint on the line normal to the reference surface of the shell in its full generality as
\begin{equation}\label{eq:13}
{\bf u}=u_{\eta_1}\,{\hat {\bf e}_{\eta_1}}+u_{\eta_2} \,{\hat {\bf e}_{\eta_2}}+u_\zeta\,{\hat {\bf e}_\zeta}
\end{equation}
where
\begin{equation}\label{eq:14}
u_{\eta_1}=[{\bf A}_{\eta_1}]_{_{(1\times n)}} \{{\boldsymbol \Phi}_{\eta_1}\}_{_{(n\times 1)}},\quad
u_{\eta_2} =[{\bf A}_{\eta_2}]_{_{(1\times m)}} \{{\boldsymbol \Phi}_{\eta_2}\}_{_{(m\times 1)}},\quad u_\zeta=[{\bf A}_\zeta]_{_{(1\times p)}} \{{\boldsymbol \Phi}_\zeta\}_{_{(p\times 1)}}
\end{equation}
and ${\bf A}_{\eta_1},\, {\bf A}_{\eta_2}$, and ${\bf A}_{\zeta}$ are the row vectors of the basis function in the $\zeta$-coordinate for the approximation of $u_{\eta_1},\ u_{\eta_2}$, and $u_{\zeta}$, respectively, whereas ${\boldsymbol \Phi}_{\eta_1},\, {\boldsymbol \Phi}_{\eta_2}$, and ${\boldsymbol \Phi}_{\zeta}$ are the column vectors of the corresponding coefficients of the basis functions defined as follows:
\begin{eqnarray}\label{eq:15}
\begin{aligned}
{\bf A}_{\eta_1}&= \left[\begin{matrix}
                    1 & f_1(\zeta) & f_2(\zeta) & \dotsc & f_{n}(\zeta)\end{matrix}\right],&
{\boldsymbol \Phi}_{\eta_1}&=\left[\begin{matrix}
                                \phi_{\eta_1}^{\scriptscriptstyle {(0)}} & \phi_{\eta_1}^{\scriptscriptstyle {(1)}} & \phi_{\eta_1}^{\scriptscriptstyle {(2)}} & \dotsc  & \phi_{\eta_1}^{\scriptscriptstyle {(n)}}
                              \end{matrix}\right]^{\rm T}\cr
{\bf A}_{\eta_2}&= \left[\begin{matrix}
                    1 & f_1(\zeta) & f_2(\zeta) & \dotsc & f_{m}(\zeta)\end{matrix}\right],&
{\boldsymbol \Phi}_{\eta_2}&=\left[\begin{matrix}
                                \phi_{\eta_2} ^{\scriptscriptstyle {(0)}} &
                                \phi_{\eta_2} ^{\scriptscriptstyle {(1)}} & \phi_{\eta_2} ^{\scriptscriptstyle {(2)}} & \dotsc  & \phi_{\eta_2}^{\scriptscriptstyle {(m)}}\cr
                              \end{matrix}\right]^{\rm T}\cr
{\bf A}_{\zeta }&= \left[\begin{matrix}
                    1 & f_1(\zeta) & f_2(\zeta) & \dotsc & f_{p}(\zeta)\end{matrix}\right], &
{\boldsymbol \Phi}_\zeta &=\left[\begin{matrix}
                                \phi_\zeta^{\scriptscriptstyle {(0)}} & \phi_\zeta^{\scriptscriptstyle {(1)}} & \phi_\zeta^{\scriptscriptstyle {(2)}} & \dotsc  & \phi_\zeta^{\scriptscriptstyle {(p)}}
                              \end{matrix}\right]^{\rm T}
\end{aligned}
\end{eqnarray}
where $n,\,m$, and $p$ are the order of approximation of the displacement component $u_{\eta_1},\, u_{\eta_2}$, and $u_{\zeta}$, respectively.
For higher order shell theory, the basis function can be taken as a polynomial in $\zeta$ as  $f_i(\zeta)=\zeta^i$; and in this case the components of the displacement field at any point can be interpreted as the Taylor series expansion about the corresponding point of the reference surface of the shell:
\begin{equation}\label{eq:16}
u_{\eta_1}=\sum_{i=0}^n \zeta^i \phi^{\scriptscriptstyle{(i)}}_{\eta_1}(\eta_1,\eta_2 ),\quad u_{\eta_2} =\sum_{i=0}^m \zeta^i \phi^{\scriptscriptstyle{(i)}}_{\eta_2} (\eta_1,\eta_2 ),\quad u_\zeta=\sum_{i=0}^{p} \zeta^i \phi^{\scriptscriptstyle{(i)}}_\zeta(\eta_1,\eta_2 )
\end{equation}
where $\phi^{\scriptscriptstyle{(0)}}_{\eta_1}  = u(\eta_1,\eta_2)$, $\phi^{\scriptscriptstyle{(0)}}_{\eta_2}= v(\eta_1,\eta_2)$, and $\phi^{\scriptscriptstyle{(0)}}_\zeta  = w(\eta_1,\eta_2)$ are the displacements of the reference surface of the shell structure at point $(\eta_1,\eta_2)$ along $\hat{\bf e}_{\eta_1}$, $\hat{\bf e}_{\eta_2}$, and $\hat{\bf e}_\zeta$ directions, respectively. Then the various variables in the above displacement field can be expressed as follows:
\begin{equation}\label{eq:17}
\phi^{\scriptscriptstyle{(i)}}_{\eta_1} ={1\over{(i)!}} \left({\partial^i u_{\eta_1}\over \partial \zeta^i}\right)_{\zeta=0},\qquad \phi^{\scriptscriptstyle{(i)}}_{\eta_2}  ={1\over{(i)!}} \left({\partial^i u_{\eta_2} \over \partial \zeta^i}\right)_{\zeta=0},\qquad \phi^{\scriptscriptstyle{(i)}}_\zeta ={1\over{(i)!}} \left({\partial^i u_\zeta\over \partial \zeta^i}\right)_{\zeta=0}
\end{equation}
For the approximation orders $n=m=1$ and $p=2$, this higher order theory specializes to the 7-parameter shell theory (see \cite{rivera2016stress}) and for $n=m=p=3$, the presented theory reduces to 12-parameter shell theory (see \cite{rivera2016new}).

Alternatively, the basis functions $f_i(\zeta)$ can also be considered as the Legendre polynomials\footnote{In the case of the higher-order theory where the approximation order is higher than 3 or 4, the Legendre polynomials would behave better numerically due to its orthogonality property. But, we should be careful in applying the boundary conditions as the scaled Legendre polynomial is not always zero at the reference surface, which could make hinged edge type boundary condition challenging to apply.} in the following form:
\begin{eqnarray}\label{eq:18}
f_i(\zeta) = P_i\left(\bar{\zeta}\right),\quad \bar{\zeta}={\zeta- \mathcal{R}_1\over\mathcal{R}_2 }
\end{eqnarray}
where $P_i(\bar{\zeta})$  is the Legendre polynomial in $\bar{\zeta}$ and
\begin{eqnarray}\label{eq:19}
\mathcal{R}_1 = {\zeta_{t}(\eta_1,\eta_2 )+\zeta_{b}(\eta_1,\eta_2 )\over 2},\ \mathcal{R}_2 = {\zeta_{t}(\eta_1,\eta_2 )-\zeta_{b}(\eta_1,\eta_2 )\over 2}.
\end{eqnarray}
Here, $\zeta_{t}$ and $\zeta_{b}$ are the $\zeta$-coordinate\footnote{Note here that $\zeta_{t}$ and $\zeta_{b}$ could be  functions of $s$ and $\theta$. Such an assumption could help analyze shell structures having some local bulge where $\zeta_{t}$ and $\zeta_{b}$ would be a known functions of $s$ and $\theta$.} of the top and bottom surface of the shell structure, respectively. Further, the displacement vector at a point can be expressed as a column vector as follows:
\begin{equation}\label{eq:20}
{\bf u}= {\bf A}{\boldsymbol \Phi},\quad \hbox{where  } {\bf A}= \left[\begin{matrix}
        {\bf A}_{\eta_1} & {\bf 0} & {\bf 0} \\
        {\bf 0} & {\bf A}_{\eta_2}  & {\bf 0} \\
        {\bf 0} & {\bf 0} & {\bf A}_{\zeta }
      \end{matrix}\right],\quad {\boldsymbol \Phi}=\left\{\begin{matrix}
                                                       {\boldsymbol \Phi}_{\eta_1} \\
                                                       {\boldsymbol \Phi}_{\eta_2} \\
                                                       {\boldsymbol \Phi}_{\zeta}
                                                     \end{matrix}\right\},\quad
{\bf u} = \left\{\begin{matrix}u_{\eta_1} \\
                               u_{\eta_2} \\
                               u_{\zeta}
                               \end{matrix}\right\}
\end{equation}
\subsubsection{Deformation gradient}
The deformation gradient in the assumed orthonormal coordinate system can be derived with ease by using the tools of exterior calculus (see the appendix of Arbind, Srinivasa, and Reddy \cite{arbind2018curved} and Appendix A for the detailed methodology) for a given coordinate system $(\eta_1,\,\eta_2,\,\zeta)$. The deformation gradient tensor can be written as
\begin{align}\label{eq:21}
{\bf F}&=F_{ij} \hat{\bf e}_i\otimes\hat{\bf e}_j,\quad \hbox{where}\quad i,j=\eta_1,\,\eta_2,\,\zeta.
\end{align}
Then the components of deformation gradient can be expressed in column vector form:
\begin{align}\label{eq:22}
 \{\tilde{\bf F}\}&=\left[\begin{matrix}
                       F_{\eta_1\eta_1} &F_{\eta_1\eta_2} & F_{\eta_1\zeta} &F_{\eta_2 \eta_1} &F_{\eta_2\eta_2} & F_{\eta_2 \zeta} & F_{\zeta \eta_1} & F_{\zeta\eta_2}& F_{\zeta\zeta}
                     \end{matrix}\right]^T
\end{align}
We will drop curly braces from the column vectors $\{\tilde{\bf F}\}$ for the further references in this study for brevity. Further, $\tilde{\bf F}$ and its first variation can be given in terms of displacement variables, ${\boldsymbol\Phi}$ as follows:
\begin{align}\label{eq:23}
 \tilde{\bf F}&=\tilde{\bf I}+ {\bf G}_1 {\boldsymbol \Phi}+{\bf G}_2 {\boldsymbol \Phi}_{,\eta_1}+{\bf G}_3 {\boldsymbol \Phi}_{,\eta_2},\quad \delta\tilde{\bf F}={\bf G}_1 \delta{\boldsymbol \Phi}+{\bf G}_2 \delta{\boldsymbol \Phi}_{,\eta_1}+{\bf G}_3 \delta{\boldsymbol \Phi}_{,\eta_2 }
\end{align}
where the expressions for $\tilde{\bf I},\,{\bf G}_1$, ${\bf G}_2$, and ${\bf G}_3$ for different curvilinear coordinate system are given in the Appendix A. Also, the components of displacement gradient, $({\bf L}={\boldsymbol\nabla}{\bf u}={\bf F}-{\bf I})$ can be written in the form of a column vector $\tilde{\bf L}$ in a similar fashion as deformation gradient (see Eq.~\eqref{eq:22}) as follows:
\begin{eqnarray}\label{eq:24}
\tilde{\bf L} = \tilde{\bf F}-\tilde{\bf I} = {\bf G}_1 {\boldsymbol \Phi}+{\bf G}_2 {\boldsymbol \Phi}_{,\eta_1}+{\bf G}_3 {\boldsymbol \Phi}_{,\eta_2 },\quad \delta\tilde{\bf L}={\bf G}_1 \delta{\boldsymbol \Phi}+{\bf G}_2 \delta{\boldsymbol \Phi}_{,\eta_1}+{\bf G}_3 \delta{\boldsymbol \Phi}_{,\eta_2 }
\end{eqnarray}
The determinant of the deformation gradient, $J$ can be expanded in terms of invariants of the displacement gradient as follows:
\begin{eqnarray}\label{eq:25}
J=\text{det}({\bf F})=\text{det}({\bf I}+{\boldsymbol\nabla}{\bf u}) = 1+{ \RN{1}+\RN{2}+\RN{3}}
\end{eqnarray}
where $\RN{1},\, \RN{2}$, and $\RN{3}$ are the first, second, and third principal invariants of ${\boldsymbol\nabla}{\bf u}$, respectively.
These invariants can be given in terms of column vector $\tilde{\bf L}$ (which ultimately can be written in terms of vector of displacement variables ${\boldsymbol \Phi}$ using Eq.~\eqref{eq:24}) as:
\begin{align}\label{eq:26}
\RN{1}= \tilde{\bf I}^{\rm T}\tilde{\bf L},\quad
\RN{2}= {1\over 2}{\bf g}_1^{\rm T}\tilde{\bf L},\quad
\RN{3}={1\over 3}\tilde{\bf L}_{cof}^{\rm T}\tilde{\bf L}
\end{align}
where
\begin{equation}\label{eq:27}
{\bf g}_1={\bf G}_{0}\tilde{\bf L},\quad \tilde{\bf L}_{cof}={1\over 2}{\bf G}_{cof}\tilde{\bf L}
\end{equation}
and the expression of ${\bf G}_0$ and ${\bf G}_{cof}$ are given in appendix B. Now, using Eqs.~\eqref{eq:25} and \eqref{eq:26}, $J$ can be rewritten as follows:
\begin{eqnarray}\label{eq:28}
J = 1+\left(\tilde{\bf I}^{\rm T} + {1\over 2}{\bf g}_1^{\rm T} + {1\over3}\tilde{\bf L}_{cof}^{\rm T}\right)\tilde{\bf L}
\end{eqnarray}
Also, the derivative of $J$ with respect to $\tilde{\bf L}$ is:
\begin{eqnarray}\label{eq:29}
{\partial J\over\partial\tilde{\bf L}}= \tilde{\bf I} + \left({\bf G}_0 + {1\over 2}{\bf G}_{cof}\right)\tilde{\bf L} = \tilde{\bf I} + \tilde{\bf G}_0\tilde{\bf L},\quad \text{where}\quad \tilde{\bf G}_0 = {\bf G}_0 + {1\over 2}{\bf G}_{cof}
\end{eqnarray}
Further, the derivative of the invariants $I_1$ and $I_2$ (see Eq.~\eqref{eq:3}) with respect to $\tilde{\bf L}$ are
\begin{eqnarray}\label{eq:30}
{\partial I_1\over \partial\tilde{\bf L}} = 2(\tilde{\bf I}+\tilde{\bf L}),\quad {\partial I_2\over \partial\tilde{\bf L}} = 4\left(\tilde{\bf I} + {\bf B}_1\tilde{\bf L}+{1\over 2}{\bf B}_2\tilde{\bf L}+{1\over 3}{\bf B}_3\tilde{\bf L}\right)
\end{eqnarray}
where the detail derivation and definition of matrices ${\bf B}_1,\, {\bf B}_2$, and ${\bf B}_3$ are given in Appendix C. Also, the volume element in the curvilinear coordinate system is given by
\begin{align}\label{eq:31}
dV &=  g\,d\eta_1\,d\eta_2\,d\zeta
\end{align}
Here $g$ is the square root of the determinant of the covariant metric tensor of the curvilinear coordinates assumed.
\subsection{Strain energy density}
The strain energy density functional of the isotropic hyperelastic material is given in Eq.~\eqref{eq:4}. Then the first variation in strain energy density functional is
\begin{align}\label{eq:32}
\delta\psi &= \delta {\bf L}:\left({\partial \psi\over \partial I_1}{\partial I_1\over \partial {\bf L}}+{\partial \psi\over \partial I_2}{\partial I_2\over \partial {\bf L}}+\left({\partial \psi\over \partial J}\right) {\partial J\over \partial{\bf L}}\right)\notag\\[4pt]
&=\delta \tilde{\bf L}\cdot\left({\beta_1\over 2}{\partial I_1\over \partial \tilde{\bf L}}+{\beta_2\over 4}{\partial I_2\over \partial \tilde{\bf L}}+\beta_3 {\partial J\over \partial\tilde{\bf L}}\right)\notag\\[4pt]
&=\delta\tilde{\bf L}\cdot\left(\left(\beta_1+\beta_2+\beta_3\right)\tilde{\bf I}  + \beta_1 \tilde{\bf L}+ \beta_2\left({\bf B}_1+{1\over 2}{\bf B}_2+{1\over 3}{\bf B}_3\right)\tilde{\bf L}+\beta_3\tilde{\bf G}_0 \tilde{\bf L}\right)
\end{align}
where $\beta_1,\, \beta_2$, and $\beta_3$ are defined in Eq.~\eqref{eq:9}.
\subsection{Governing equation for higher-order shell theory}
To derive the governing equation for the general higher-order shell theory, let us consider that ${\bf b}$ be the body force per unit mass, and ${\bf q}$ be the traction force (transform back to reference configuration) applied on the boundary surface of the structure. Further, to obtain the governing equation, we start from the following virtual work statement (see \cite{reddy2018energy}) for given strain energy density:
\begin{align}\label{eq:33}
{\bf 0}&= \int_{\mathcal B}\left(\delta \psi-\rho_0{\bf b}\cdot\delta {\bf u}\right)\, dV - \oint_{\partial B}{\bf q}\cdot\delta{\bf u}\, dS\notag\\[4pt]
&=\int_{\mathcal B}\Big[\delta \tilde{\bf L}\cdot\left(\left(\beta_1+\beta_2+\beta_3\right)\tilde{\bf I}  + \beta_1 \tilde{\bf L}+ \beta_2\left({\bf B}_1+{1\over 2}{\bf B}_2+{1\over 3}{\bf B}_3\right)\tilde{\bf L}\right.\notag\\[4pt]&\hskip20pt\left.+\beta_3\left({\bf G}_0 + {1\over 2}{\bf G}_{cof}\right)\tilde{\bf L}\right)
-\rho_0{\bf b}\cdot\delta {\bf u}\Big]\, dV - \oint_{\partial B}{\bf q}\cdot\delta{\bf u}\, dS\notag\\[4pt]
&= \int_{A} \big[\delta{\boldsymbol\Phi}\cdot({\bf h}_1{\boldsymbol\Phi}+{\bf h}_2{\boldsymbol\Phi}_{,\eta_1}+{\bf h}_3{\boldsymbol\Phi}_{,\eta_2 })+ \delta{\boldsymbol\Phi}_{,\eta_1}\cdot({\bf h}_4{\boldsymbol\Phi}+{\bf h}_5{\boldsymbol\Phi}_{,\eta_1}+{\bf h}_6{\boldsymbol\Phi}_{,\eta_2 })\cr&\hskip30pt+ \delta{\boldsymbol\Phi}_{,\eta_2 }\cdot({\bf h}_7{\boldsymbol\Phi}+{\bf h}_8{\boldsymbol\Phi}_{,\eta_1}+{\bf h}_9{\boldsymbol\Phi}_{,\eta_2} )+\delta{\boldsymbol\Phi}\cdot({\bf f}_1-{\bf f}_0)+\delta{\boldsymbol\Phi}_{,\eta_1}\cdot{\bf f}_2+\delta{\boldsymbol\Phi}_{,\eta_2}\cdot{\bf f}_3\big] \, d\eta_1\,d\eta_2 \notag\\[4pt]
&\qquad-\oint_{\partial A}[\delta{\boldsymbol\Phi}\cdot\hat{\bf f}_l]\,dl
\end{align}
where
\begin{align}\label{eq:34}
{\bf h}_1&=\int_{\zeta_b}^{\zeta_t} {\bf G}_1^{\rm T}{\bf B}_0{\bf G}_1\, g\,d\zeta,\quad {\bf h}_2=\int_{\zeta_b}^{\zeta_t} {\bf G}_1^{\rm T}{\bf B}_0{\bf G}_2\, g\,d\zeta,\quad {\bf h}_3=\int_{\zeta_b}^{\zeta_t} {\bf G}_1^{\rm T}{\bf B}_0{\bf G}_3\, g\,d\zeta\notag\\[4pt]
{\bf h}_4&=\int_{\zeta_b}^{\zeta_t} {\bf G}_2^{\rm T}{\bf B}_0{\bf G}_1\, g\,d\zeta,\quad {\bf h}_5=\int_{\zeta_b}^{\zeta_t} {\bf G}_2^{\rm T}{\bf B}_0{\bf G}_2\, g\,d\zeta,\quad {\bf h}_6=\int_{\zeta_b}^{\zeta_t} {\bf G}_2^{\rm T}{\bf B}_0{\bf G}_3\, g\,d\zeta\notag\\[4pt]
{\bf h}_7&=\int_{\zeta_b}^{\zeta_t} {\bf G}_3^{\rm T}{\bf B}_0{\bf G}_1\, g\,d\zeta,\quad {\bf h}_8=\int_{\zeta_b}^{\zeta_t} {\bf G}_3^{\rm T}{\bf B}_0{\bf G}_2\, g\,d\zeta,\quad {\bf h}_9=\int_{\zeta_b}^{\zeta_t} {\bf G}_3^{\rm T}{\bf B}_0{\bf G}_3\, g\,d\zeta
\end{align}
and
\begin{align}\label{eq:35}
{\bf f}_1&=\int_{\zeta_b}^{\zeta_t} {\bf G}_1^{\rm T}(\beta_1+\beta_2+\beta_3 )\tilde{\bf I}\,g\, d\zeta ,& {\bf f}_2&=\int_{\zeta_b}^{\zeta_t} {\bf G}_2^{\rm T}(\beta_1+\beta_2+\beta_3 ) \tilde{\bf I} \,g\, d\zeta\cr
{\bf f}_3&=\int_{\zeta_b}^{\zeta_t} {\bf G}_3^{\rm T}(\beta_1+\beta_2+\beta_3) \tilde{\bf I} \,g\, d\zeta&&\cr
{\bf f}_0&=\int_{\zeta_b}^{\zeta_t} {\bf A}^{\rm T}\rho_0 {\bf b} \, g\, d\zeta + \sqrt{G}_i {\bf A}^{\rm T}{\bf q}_i&
\hat{\bf f}_l&=\int_{\zeta_b}^{\zeta_t}\hskip-7pt{\bf A}^{\rm T}{\bf q}_l\, d\zeta
\end{align}
\begin{equation}\label{eq:36}
{\bf B}_0=\left[\beta_1{\bf I}+\beta_2\left({\bf B}_1+{1\over 2}{\bf B}_2+{1\over 3}{\bf B}_3\right)+\beta_3\left({\bf G}_0 + {1\over 2}{\bf G}_{cof}\right)\right]
\end{equation}
Moreover, ${\bf q}_i$ and ${\bf q}_l$ are the surface tractions applied on $i$th lateral surface and edge side surfaces of the shell structure, respectively, which are transformed back to the reference surface of the structure; $G_i$ is the determinant of covariant metric tensor of the surface coordinate $(\eta_1,\eta_2)$ for the $i$th lateral boundary surface (see Appendix C of \cite{arbind2018Neo-Hookean} for general curved tubular shell surface). Next, the governing equation is obtained from Eq.~\eqref{eq:33} as follows:
\begin{eqnarray}\label{eq:37}
\delta {\boldsymbol\Phi}:\quad &&{\bf h}_1{\boldsymbol\Phi}+{\bf h}_2{\boldsymbol\Phi}_{,\eta_1}+{\bf h}_3{\boldsymbol\Phi}_{,\eta_2 }-({\bf h}_4{\boldsymbol\Phi}+{\bf h}_5{\boldsymbol\Phi}_{,\eta_1}+{\bf h}_6{\boldsymbol\Phi}_{,\eta_2 })_{,\eta_1}\notag\\[2pt]&&-({\bf h}_7{\boldsymbol\Phi}+{\bf h}_8{\boldsymbol\Phi}_{,\eta_1}+{\bf h}_9{\boldsymbol\Phi}_{,\eta_2})_{,\eta_2}={\bf f}_0-{\bf f}_1+{\bf f}_{2,\eta_1}+{\bf f}_{3,\eta_2}
\end{eqnarray}
and at the boundary line of the shell surface, the boundary conditions are:
\begin{eqnarray}\label{eq:38}
  {\boldsymbol \Phi} =0 &\text{or}&\hat{\bf f}_l +({\bf h}_4{\boldsymbol\Phi}+{\bf h}_5{\boldsymbol\Phi}_{,\eta_1}+{\bf h}_6{\boldsymbol\Phi}_{,\eta_2}+{\bf f}_2)n_{\eta_1}+({\bf h}_7{\boldsymbol\Phi}+{\bf h}_8{\boldsymbol\Phi}_{,\eta_1}+{\bf h}_9{\boldsymbol\Phi}_{,\eta_2}+{\bf f}_3 )n_{\eta_2} =0\cr&
\end{eqnarray}
where $n_{\eta_1}$ and $n_{\eta_2}$ are the components of normal direction to the edges of shell-reference surface in $(\eta_1,\, \eta_2)$ space.
\section{Weak form finite element model}
The obtained governing equation of the higher-order shell theory is the nonlinear partial differential equation in two dimensions, which can be solved numerically. In this section, we develop the weak form finite element model for isotropic compressible hyperelastic material for the higher-order shell theory.
Towards this end, we consider the following Lagrangian\footnote{The load is applied very slowly such that the kinetic energy could be neglected.} for a general finite element $\Omega^{e}$, $\gamma$:
\begin{eqnarray}\label{eq:39}
\mathcal{L}_p= \int_{\Omega^{e}}\int_{\zeta_b}^{\zeta_t} \Big[\psi(I_1,I_2,J)\Big]g\, d\zeta\,ds\,d\theta -\mathcal{V}_{\Omega^{e}}
\end{eqnarray}
where $\mathcal{V}_{\Omega^{e}}$ is the work done by external forces on the element ${\Omega^{e}}$ and $\psi$ is the strain energy density function. Next, we obtain the weak form from the above Lagrangian as follows:
\begin{align}\label{eq:40}
0&=  \int_{{\mathcal B}_e}\left({\partial \psi\over \partial \tilde{\bf L}}\right)\cdot\delta \tilde{\bf L}-\rho_0{\bf b}\cdot\delta {\bf u}\, dV - \oint_{\partial{\mathcal B}_e}{\bf q}\cdot\delta{\bf u}\, dS\notag\\[4pt]
&=\int_{\Omega_e} \big[\delta{\boldsymbol\Phi}\cdot({\bf H}_1{\boldsymbol\Phi}+{\bf H}_2{\boldsymbol\Phi}_{,\eta_1}+{\bf H}_3{\boldsymbol\Phi}_{,\eta_2})+ \delta{\boldsymbol\Phi}_{,\eta_1}\cdot({\bf H}_4{\boldsymbol\Phi}+{\bf H}_5{\boldsymbol\Phi}_{,\eta_1}+{\bf H}_6{\boldsymbol\Phi}_{,\eta_2})\cr&\hskip20pt+ \delta{\boldsymbol\Phi}_{,\eta_2 }\cdot({\bf H}_7{\boldsymbol\Phi}+{\bf H}_8{\boldsymbol\Phi}_{,\eta_1}+{\bf H}_9{\boldsymbol\Phi}_{,\eta_2} )+\delta{\boldsymbol\Phi}\cdot({\bf f}_1-{\bf f}_0)+\delta{\boldsymbol\Phi}_{,\eta_1}\cdot{\bf f}_2+\delta{\boldsymbol\Phi}_{,\eta_2}\cdot{\bf f}_3\big] \, d\eta_1\,d\eta_2\cr
\end{align}
In the above equation, ${\bf I}$ is a $(9\times9)$ identity matrix. Now, we approximate the degrees of freedom vector as
\begin{equation}\label{eq:41}
{\boldsymbol \Phi}(\eta_1,\eta_2) = {\boldsymbol \Psi}(\eta_1,\eta_2){\bf U}
\end{equation}
where ${\boldsymbol \Psi}(\eta_1,\eta_2)$ is the matrix of interpolation functions, which are functions of the coordinates $(\eta_1,\eta_2)$; ${\bf U}$ is a vector of the nodal values of the variables corresponding to displacement components, given as follows:
\begin{align}
{\boldsymbol \Psi}&=\left[\begin{array}{cccccccccc} \psi^{\scriptscriptstyle{(1)}}_1 &\hdots&\psi^{\scriptscriptstyle{(1)}}_{\tilde n_1} & 0  & \hdots  & 0 & \hdots  & 0 & \hdots &0\\
0  & \hdots  & 0 &\psi^{\scriptscriptstyle{(2)}}_1 &\hdots&\psi^{\scriptscriptstyle{(2)}}_{\tilde n_2} & \hdots  & 0 & \hdots &0\\
\vdots  & \ddots  & \vdots &\vdots &\ddots&\vdots & \ddots  & \vdots & \ddots &\vdots\\
  0  & \hdots  & 0 & 0&\hdots & 0 &\hdots&\psi^{\scriptscriptstyle{(\hat{n})}}_1 &\hdots&\psi^{\scriptscriptstyle{(\hat{n})}}_{\tilde n_p}
  \end{array}\right ]\label{eq:42}\\
\textbf{U}&=\left[\begin{array}{cccccccccc} u_{1_1} & \hdots & u_{1_{\tilde n_1}}  &  u_{2_1} & \hdots & u_{2_{\tilde n_2}}  & \hdots & u_{n_1} & \hdots & u_{\hat n_{\tilde n_n}} \end{array}\right ]^{\rm T}\label{eq:43}
\end{align}
where $\tilde n_1,\tilde n_2, \ldots \tilde n_{\hat{n}}$ are the number of nodal values of  $u_1,u_2,\ldots,u_{\hat{n}}$, respectively, and $\hat{n}(=n+m+p+3)$ is the total number of Dofs at any node. Also,
\begin{equation}\label{eq:44}
\begin{aligned}
&u_1&= \phi_{\eta_1}^{\scriptscriptstyle{(0)}},\quad & u_2 & = \phi_{\eta_1}^{\scriptscriptstyle{(1)}},\quad &\cdots &u_{n+1}&=\phi_{\eta_1}^{\scriptscriptstyle{(n)}}\\
&u_{n+2}&= \phi_{\eta_2}^{\scriptscriptstyle{(0)}},\quad    & u_{n+3} & = \phi_{\eta_2}^{\scriptscriptstyle{(1)}},\quad &\cdots&u_{n+m+2}& = \phi_{\eta_2}^{\scriptscriptstyle{( m)}}\\
&u_{n+m+3}&= \phi_{\zeta}^{\scriptscriptstyle{(0)}},\quad    & u_{n+m+4} & = \phi_{\zeta}^{\scriptscriptstyle{(1)}},\quad &\cdots&u_{n+m+p+3}& = \phi_{\zeta}^{\scriptscriptstyle{(p)}}
\end{aligned}
\end{equation}
We substitute the approximation of the displacement variables and $\delta{\boldsymbol \Phi}^a={\boldsymbol \Psi} \hat{\bf I}$ (where $ \hat{\bf I}$ is the column vector with all element unity and as many elements as the columns of ${\boldsymbol \Psi}$) into the weak form in Eq.~\eqref{eq:38} to arrive at the following finite element equations:
\begin{equation}\label{eq:45}
{\bf K} {\bf U} - {\bf f} = {\bf 0}
\end{equation}
where ${\bf K}$ and ${\bf f}$ are the stiffness matrix and force vector (both the stiffness matrix and force vector are nonlinear as they depends on the current solution vector ${\bf U}$), respectively, and they are defined as follows:
\begin{align}\label{eq:46}
{\bf K} &= \int_{\Omega_e}\Big[ {\boldsymbol \Psi}^{\rm T}\left({\bf H}_1{\boldsymbol \Psi}+ {\bf H}_2{\boldsymbol \Psi}_{,\eta_1}+ {\bf H}_3{\boldsymbol \Psi}_{,\eta_2} \right)+{\boldsymbol \Psi}^{\rm T}_{,\eta_1}\left({\bf H}_4{\boldsymbol \Psi}+ {\bf H}_5{\boldsymbol \Psi}_{,\eta_1} + {\bf H}_6{\boldsymbol \Psi}_{,\eta_2}\right)\notag\\&\hskip40pt+{\boldsymbol \Psi}^{\rm T}_{,\eta_2}\left({\bf H}_7{\boldsymbol \Psi}+ {\bf H}_8{\boldsymbol \Psi}_{,\eta_1} + {\bf H}_9{\boldsymbol \Psi}_{,\eta_2}\right)\Big]\, d\eta_1\,d\eta_2\notag\\[-8pt]
&\notag\\[-8pt]
{\bf f} &= \int_{\Omega_e}{\boldsymbol \Psi}^{\rm T}({\bf f}_0-{\bf f}_{1})-{\boldsymbol \Psi}^{\rm T}_{,\eta_1}\,{\bf f}_{2}-{\boldsymbol \Psi}^{\rm T}_{,\eta_2}\,{\bf f}_3\,d\eta_1\,d\eta_2
\end{align}
with
\begin{align}\label{eq:47}
{\bf H}_1&=\int_{\zeta_b}^{\zeta_t}  {\bf G}_1^{\rm T}\hat{\bf B}_0{\bf G}_1\,g d\zeta ,& {\bf H}_2&=\int_{\zeta_b}^{\zeta_t}  {\bf G}_1^{\rm T}\hat{\bf B}_0{\bf G}_2\,g\, d\zeta,& {\bf H}_3&=\int_{\zeta_b}^{\zeta_t}  {\bf G}_1^{\rm T}\hat{\bf B}_0{\bf G}_3\,g\, d\zeta\notag\\[4pt]
{\bf H}_4&=\int_{\zeta_b}^{\zeta_t}  {\bf G}_2^{\rm T}\hat{\bf B}_0{\bf G}_1\,g\, d\zeta ,& {\bf H}_5&=\int_{\zeta_b}^{\zeta_t}  {\bf G}_2^{\rm T}\hat{\bf B}_0{\bf G}_2\,g\, d\zeta,& {\bf H}_6&=\int_{\zeta_b}^{\zeta_t}  {\bf G}_2^{\rm T}\hat{\bf B}_0{\bf G}_3\,g\, d\zeta\notag\\[4pt]
{\bf H}_7&=\int_{\zeta_b}^{\zeta_t}  {\bf G}_3^{\rm T}\hat{\bf B}_0{\bf G}_1\,g\, d\zeta ,& {\bf H}_8&=\int_{\zeta_b}^{\zeta_t}  {\bf G}_3^{\rm T}\hat{\bf B}_0{\bf G}_2\,g\, d\zeta,& {\bf H}_9&=\int_{\zeta_b}^{\zeta_t}  {\bf G}_3^{\rm T}\hat{\bf B}_0{\bf G}_3\,g\, d\zeta
\end{align}
\begin{align}\label{eq:48}
{\bf f}_1&=\int_{\zeta_b}^{\zeta_t}  {\bf G}_1^{\rm T}\left(\beta_1 +\beta_2+\beta_3\right)\tilde{\bf I}\,g\,d\zeta\notag\\[4pt]
{\bf f}_2&=\int_{\zeta_b}^{\zeta_t}  {\bf G}_2^{\rm T}\left(\beta_1 +\beta_2+\beta_3\right) \tilde{\bf I} \,g\, d\zeta\notag\\[4pt]
{\bf f}_3&=\int_{\zeta_b}^{\zeta_t}  {\bf G}_3^{\rm T}\left(\beta_1 +\beta_2+\beta_3\right) \tilde{\bf I} \,g\, d\zeta
\end{align}
and
\begin{eqnarray}\label{eq:49}
\hat{\bf B}_0&=&\beta_1{\bf I}+\beta_2\left({\bf B}_1+{1\over 2}{\bf B}_2+{1\over 3}{\bf B}_3\right)+\beta_3\left({\bf G}_0 + {1\over 2}{\bf G}_{cof}\right)
\end{eqnarray}
The vectors ${\bf f}_0$ and $\hat{\bf f}_0$ are the same as given in Eq.~\eqref{eq:35}. The stiffness matrices are nonlinear and not symmetric. Again, we can apply Newton's method (see \cite{reddy2015introduction}) to solve the nonlinear FE equation. The tangent matrix ${\bf T}$ in this case is given as
\begin{align}\label{eq:50}
{\bf T}&= \int_{\Omega_e}\Big[ {\boldsymbol \Psi}^{\rm T}\left(\tilde{\bf H}_1{\boldsymbol \Psi}+ \tilde{\bf H}_2{\boldsymbol \Psi}_{,\eta_1}+ \tilde{\bf H}_3{\boldsymbol \Psi}_{,\eta_2} \right)+{\boldsymbol \Psi}^{\rm T}_{,\eta_1}\left(\tilde{\bf H}_4{\boldsymbol \Psi}+ \tilde{\bf H}_5{\boldsymbol \Psi}_{,\eta_1} + \tilde{\bf H}_6{\boldsymbol \Psi}_{,\eta_2}\right)\notag\\&\hskip40pt+{\boldsymbol \Psi}^{\rm T}_{,\eta_2}\left(\tilde{\bf H}_7{\boldsymbol \Psi}+ \tilde{\bf H}_8{\boldsymbol \Psi}_{,\eta_1} + \tilde{\bf H}_9{\boldsymbol \Psi}_{,\eta_2}\right)\Big]\, d\eta_1\,d\eta_2- {\bf T}_f
\end{align}
and $\tilde {\bf H}_i$ are defined as
\begin{align}\label{eq:51}
\tilde{\bf H}_1&=\int_{\zeta_b}^{\zeta_t}  {\bf G}_1^{\rm T}\tilde{\bf B}_0{\bf G}_1\,g\, d\zeta ,& \tilde{\bf H}_2&=\int_{\zeta_b}^{\zeta_t}  {\bf G}_1^{\rm T}\tilde{\bf B}_0{\bf G}_2\,g\, d\zeta,& \tilde{\bf H}_3&=\int_{\zeta_b}^{\zeta_t}  {\bf G}_1^{\rm T}\tilde{\bf B}_0{\bf G}_3\,g\, d\zeta\notag\\[4pt]
\tilde{\bf H}_4&=\int_{\zeta_b}^{\zeta_t}  {\bf G}_2^{\rm T}\tilde{\bf B}_0{\bf G}_1\,g\, d\zeta ,& \tilde{\bf H}_5&=\int_{\zeta_b}^{\zeta_t}  {\bf G}_2^{\rm T}\tilde{\bf B}_0{\bf G}_2\,g\, d\zeta,& \tilde{\bf H}_6&=\int_{\zeta_b}^{\zeta_t}  {\bf G}_2^{\rm T}\tilde{\bf B}_0{\bf G}_3\,g\, d\zeta\notag\\[4pt]
\tilde{\bf H}_7&=\int_{\zeta_b}^{\zeta_t}  {\bf G}_3^{\rm T}\tilde{\bf B}_0{\bf G}_1\,g\, d\zeta ,& \tilde{\bf H}_8&=\int_{\zeta_b}^{\zeta_t}  {\bf G}_3^{\rm T}\tilde{\bf B}_0{\bf G}_2\,g\, d\zeta,& \tilde{\bf H}_9&=\int_{\zeta_b}^{\zeta_t}  {\bf G}_3^{\rm T}\tilde{\bf B}_0{\bf G}_3\,g\, d\zeta
\end{align}
where
\begin{align}\label{eq:52}
\tilde{\bf B}_0&= \beta_1{\bf I}+\beta_2\left({\bf B}_1+{\bf B}_2+{\bf B}_3\right)+\beta_3\left({\bf G}_0 + {\bf G}_{cof}\right)\notag\\[2pt]
&+\left[\beta_{1,1}\tilde{\bf F}+\beta_{2,1}\tilde{\bf L}_{_{(FC)}}+\beta_{3,1}\left(\tilde{\bf I} + \tilde{\bf G}_0\tilde{\bf L}\right)\right]\tilde{\bf F}^{\rm T}\notag\\[2pt]
&+\left[\beta_{1,2}\tilde{\bf F}+\beta_{2,2}\tilde{\bf L}_{_{(FC)}}+\beta_{3,2}\left(\tilde{\bf I} + \tilde{\bf G}_0\tilde{\bf L}\right)\right]\tilde{\bf L}_{_{(FC)}}^{\rm T}\notag\\[2pt]
&+\left[\beta_{1,3}\tilde{\bf F}+\beta_{2,3}\tilde{\bf L}_{_{(FC)}}+\beta_{3,3}\left(\tilde{\bf I} + \tilde{\bf G}_0\tilde{\bf L}\right)\right]\left(\tilde{\bf I} + \tilde{\bf G}_0\tilde{\bf L}\right)^{\rm T}
\end{align}\label{eq:53}
where $\tilde{\bf L}_{_{(FC)}}$ (see Appendix B) is defined as follows:
\begin{align}
\tilde{\bf L}_{_{(FC)}} &= \left(\tilde{\bf I} + {\bf B}_1\tilde{\bf L}+{1\over 2}{\bf B}_2\tilde{\bf L}+{1\over 3}{\bf B}_3\tilde{\bf L}\right)
\end{align}
and
\begin{align}\label{eq:54}
\beta_{n,1}&=2{\partial \beta_n\over \partial\, I_1},\quad \beta_{n,2}=4{\partial \beta_n\over \partial\,I_2},\quad
\beta_{n,3}={\partial \beta_n\over \partial J}, \quad \text{for}\quad n=1,2,3
\end{align}
and in Eq.~\eqref{eq:50} the term ${\bf T}_f$ comes from the derivative of the force vector with respect to ${\bf U}$ (see Eq. (58) in \cite{arbind2018Neo-Hookean}). Also, we note that the tangent matrix is symmetric if ${\bf T}_f$ is zero.
\section{Specialization to various models of isotropic hyperelastic material}
\subsection{Saint Venant–Kirchhoff nonlinear material model}
For the Saint Venant-–Kirchhoff nonlinear material model, the second Piola--Kirchhoff stress tensor and the Green--Lagrange strain tensor are related linearly in the same way as in stress--strain relation for linear elasticity. This model is the simplest model of nonlinear hyperelastic material. In this case, the strain energy density per unit reference volume is given by
\begin{eqnarray}\label{eq:55}
\psi={1\over2}{\bf S}:{\bf E} = {\mu\over 4}(I_2-2I_1+3)+{\lambda\over 8} (I_1-3)^2
\end{eqnarray}
where ${\bf S}$ and ${\bf E}$ are the second Piola--Kirchhoff stress tensor, and the Green--Lagrange strain tensor, respectively; $\mu$ and $\lambda$ are the  La\'me parameters. The $\beta$'s and its derivatives used in the finite element model, described in section 4.1, can be specialized for this material model as follows:
\begin{eqnarray}\label{eq:56}
&&\beta_1= {\lambda\over2}(I_1-3)-\mu,\quad\beta_2=\mu,\quad \beta_{1,1}=\lambda \cr
&&\beta_3=\beta_{3,i}=\beta_{2,i}=\beta_{1,2}=\beta_{1,3}=0,\quad \hbox{where } i=1,2,3.
\end{eqnarray}
\subsection{Compressible neo-Hookean model}
For the compressible neo-Hookean solid, the strain energy density per unit reference volume is given by
\begin{eqnarray}\label{eq:57}
\psi={\lambda\over 2}\log^2(J)-\mu \log(J)+{\mu\over 2}(I_1-3)
\end{eqnarray}
where $\lambda$ and $\mu$ are the La\'me parameters. The $\beta$'s and their derivatives can be specialized for this model as follows:
\begin{eqnarray}\label{eq:58}
&&\beta_1= \mu,\quad \beta_3={1\over J}(\lambda \log(J)-\mu),\quad
\beta_{3,3}={\lambda\over J^2}-{2\over J^2}(\lambda \log(J)-\mu)\cr\cr
&&\beta_2=\beta_{1,1}=\beta_{1,3}=\beta_{3,1} =\beta_{1,2}=\beta_{2,1}=\beta_{2,2}=\beta_{2,3}=\beta_{3,2}=0
\end{eqnarray}
\subsection{Compressible Mooney-Rivlin model}
Another example of compressible hyperelastic material model is Mooney-Rivlin model. The strain energy density per unit reference volume is given by:
\begin{align}\label{eq:59}
\psi&=C_1(J^{(-2/3)}I_C-3) + C_2(J^{(-2/3)}II_C-3) +{K\over 2}(J-1)^2\cr
&=C_1(J^{(-2/3)}I_1-3) + C_2\left({1\over2}J^{(-2/3)}(I_1^2-I_2)-3\right) +{K\over 2}(J-1)^2
\end{align}
where $C_1$ and $C_2$ are the material constants and $K$ is the bulk modulus. The $\beta$'s and its derivatives can be then given as follows:
\begin{eqnarray}\label{eq:60}
&&\beta_1= 2J^{(-2/3)}(C_1+C_2I_1),\quad \beta_2= -2C_2J^{(-2/3)},\cr\cr &&\beta_3=-{2\over3}J^{(-5/3)}\left(C_1I_1 + {1\over2}C_2(I_1^2-I_2)\right) +K(J-1),\cr
&&\beta_{1,1}= 4C_2J^{(-2/3)},\quad \beta_{1,3}= -{4\over3}J^{(-5/3)}(C_1+C_2I_1),\quad \beta_{2,3}= {4\over3}C_2 J^{(-5/3)}\cr
&&\beta_{3,1}=-{4\over3}J^{(-5/3)}\left(C_1 + C_2I_1\right),
\quad\beta_{3,2}={4\over3}C_2 J^{(-5/3)}\cr\cr
&&\beta_{3,3}={10\over9}J^{(-8/3)}\left(C_1I_1 + {1\over2}C_2(I_1^2-I_2)\right) +K\cr\cr
&&\beta_{1,2}=\beta_{2,1}=\beta_{2,2}=0
\end{eqnarray}
\section{Numerical examples}
In this section, we present a number of numerical examples illustrating the formulation of shell theory presented in this study. For all the examples, the surface coordinate $({\eta_1,\eta_2})$ are taken as $({\theta,s})$ as presented in Appendix A, and all the shell surfaces presented in the examples can be seen as a curved or straight pipe with constant or variable radius. Also, the mesh discretization $(n_\theta\times n_s)$ means $n_\theta$ and $n_s$ elements along the $\theta-$ and $s$- direction, respectively for all the numerical examples.
\subsection{Semi-cylindrical shell subjected to point load }
\begin{figure}[!htbp]
\centering
\begin{adjustbox}{center}
\begin{tabular}{cc}
\subfloat[ Natural configuration, load and boundary condition ]{\includegraphics[trim=2cm 1cm 0.5cm 1cm,clip, scale=0.31]{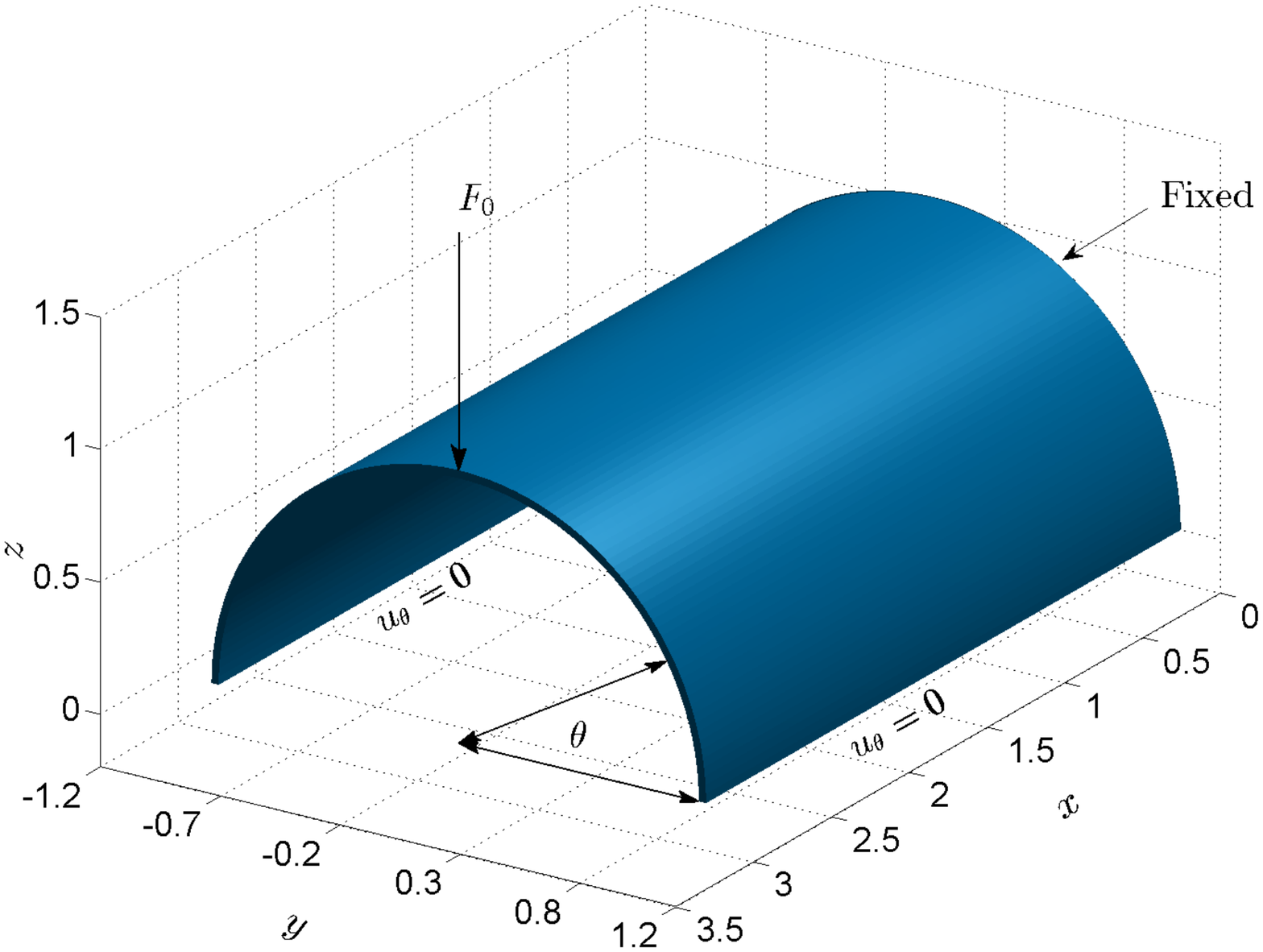}}
&\hskip-0.5cm \subfloat[ Deformed shape for $F_0 = 1614.5$ lb. ]{\includegraphics[trim=2cm 1cm 0.8cm 1cm,clip,scale=0.31]{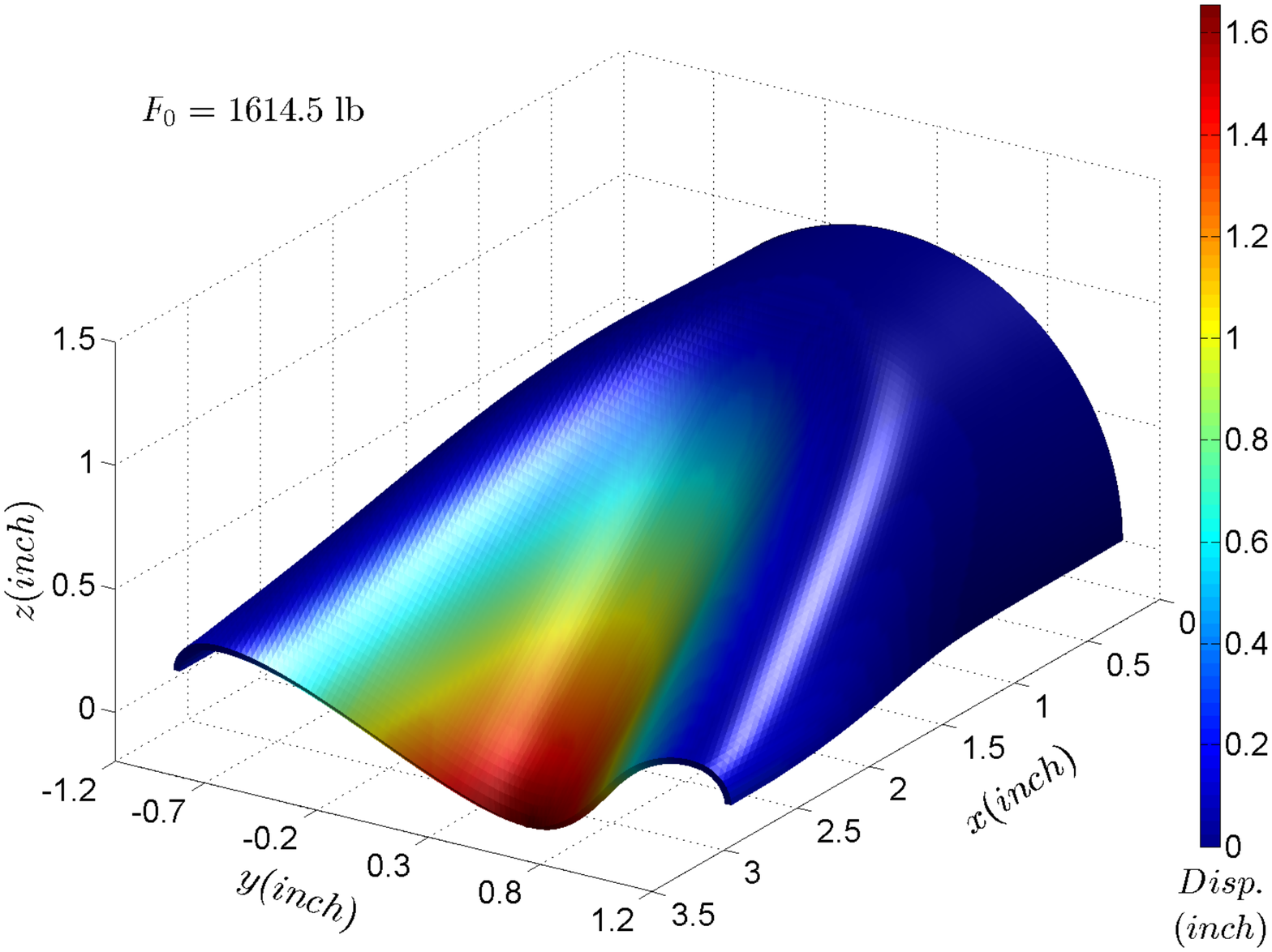}}
\end{tabular}
\end{adjustbox}
\caption{Natural configuration and the deformed geometry of the cylindrical shell panel subjected to a point load.}\label{fig:3}
\end{figure}

First, we consider an example of a very common benchmark problem of semi-cylindrical shell subjected to point load, as shown in Fig.~\ref{fig:3}(a) for Saint Venant–Kirchhoff nonlinear material model. The geometric and material properties of the semi-cylinder are the following:
\begin{equation}\label{eq:61}
L = 3.048\hbox{ in.} , \qquad R  = 1.016 \hbox{ in.} , \qquad h  = 0.03\hbox{ in.}
\end{equation}
\begin{equation}\label{eq:62}
E = 20.685 \times 10^6\hbox{ psi.}  , \qquad \nu = 0.3 ,\qquad
\end{equation}
where $L$, $R$, and $h$ are the length, mean radius, and thickness of the semicircular cylinder; $E$ and $\nu$ are the modulus of elasticity and Poisson's ratio, respectively. One end of the cylindrical panel is completely fixed, and the straight edges are constrained to have $u_\theta$ equal to zero.

Uniform $(16\times 6)$ cubic spectral Lagrange elements are used for the nonlinear finite element analysis. $(4\times4)$ Gauss points are used in an element to integrate the stiffness and tangent matrices. For solving the nonlinear finite element equation, arc-length method (see \cite{crisfield1981fast} and \cite{reddy2015introduction}) is employed to have program-controlled load increment with the error tolerance equal to $10^{-3}$.

\begin{figure}[!htbp]
\centerline{\includegraphics[trim=1cm 1cm 1cm 1cm,clip,scale=0.46]{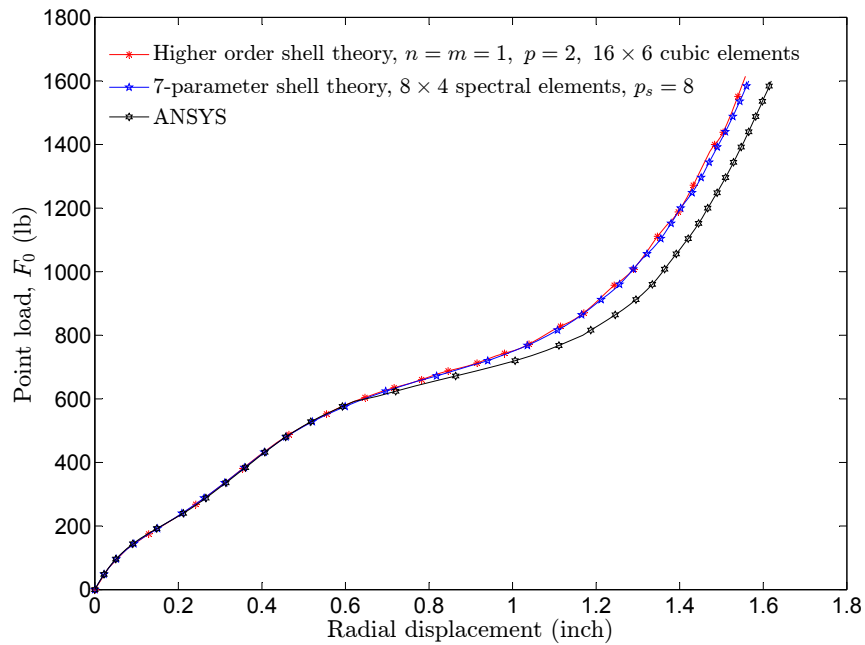}}
\caption{Comparison of radial displacement of mid surface at the point of application of load on the cylindrical panel by general higher order shell theory to 7-parameter shell theory (\cite{rivera2016stress}) and ANSYS.}
\label{fig:4}
\end{figure}

Figure~\ref{fig:3}(b) shows the deformed shape of the semi-cylinder at point load $1614.5$ lb. Also, Fig.~\ref{fig:4} shows the load Vs. radial displacement (at the point of load application) plot and the result has been compared with
$7$-parameter shell theory for linear material (see Rivera and Reddy \cite{rivera2016stress}) and ANSYS. The solution from the present study is in good agreement with the study of the 7-parameter shell theory of Rivera and Reddy \cite{rivera2016stress}, whereas ANSYS solution is not accurate for large deformation.

\subsection{Hyperboloidal shell subjected to point loads}
\begin{figure}[!htbp]
\centerline{\includegraphics[trim=1cm 1cm 1cm 1cm,clip,scale=0.5]{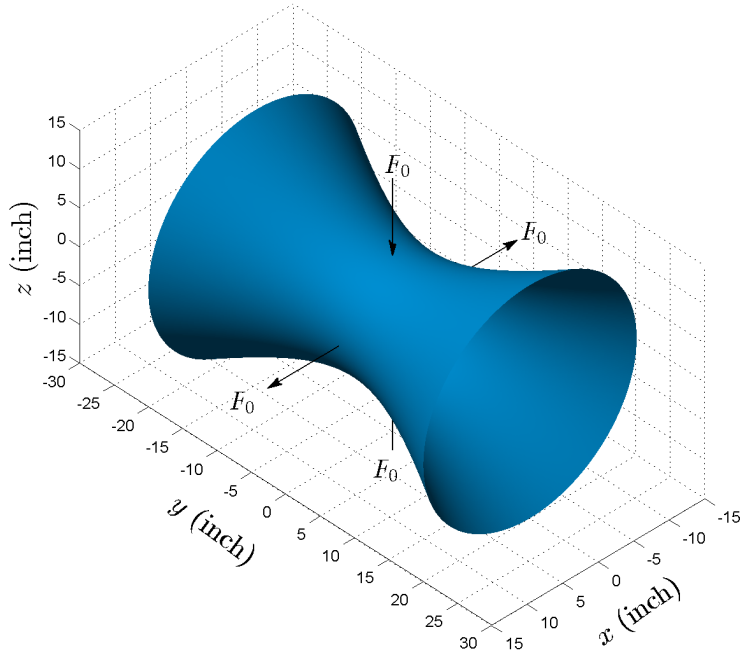}}
\caption{Undeformed geometry of the hyperboloidal shell. Both the end cross-section are completely free and four point loads are applied at the mid cross-section of the shell.}
\label{fig:5}
\end{figure}

In this example, we consider the deformation of a hyperboloidal shell under point loads, as shown in Fig.~\ref{fig:5}. Both the end cross-sections of the hyperboloidal shell are completely free. Two nonlinear material models, namely, compressible neo-Hookean model and Saint Venant–Kirchhoff nonlinear material model have been considered for this case. The geometric and material parameters of the cylindrical panel used are:
\begin{equation}\label{eq:63}
L = 40 \hbox{ in.} , \qquad R_1  = 7.5 \hbox{ in.}, \qquad R_2  = 15 \hbox{ in.} , \qquad h  = 0.04 \hbox{ in.}
\end{equation}
\begin{equation}\label{eq:64}
\mu = 1.6 \times 10^6  \hbox{ psi.}, \qquad \nu = 0.25
\end{equation}
where $L$, $R_1$, $R_2$, and $h$ are the length, mean radii (minimum and maximum radii, respectively), and thickness of the shell; $\mu$ and $\nu$ are the shear modulus and Poisson's ratio, respectively.

\noindent Exploiting the symmetry of the problem, only $1/8$th of the full domain is considered as the computational domain for the nonlinear finite element analysis. Three different meshes with different order ($p_s$) of spectral Lagrange elements have been used for the computational domain; (i) $(5 \times 5)$ spectral elements, $p_s = 8$, (ii) $(10 \times 10)$ spectral elements, $p_s = 4$, and (iii) $(20 \times 20)$ spectral elements, $p_s = 2$. All meshes contain $41\times41$ nodes. Full Gauss quadrature is used to integrate the stiffness and tangent matrices for all orders of spectral elements. Newton's method has been applied to solve the nonlinear finite element equations.

\begin{figure}[!htbp]
\centering
\begin{adjustbox}{center}
\begin{tabular}{cc}
\subfloat[ $F_0 = 124.7$ lb. ]{\includegraphics[trim=2cm 1cm 0.5cm 1cm,clip, scale=0.4]{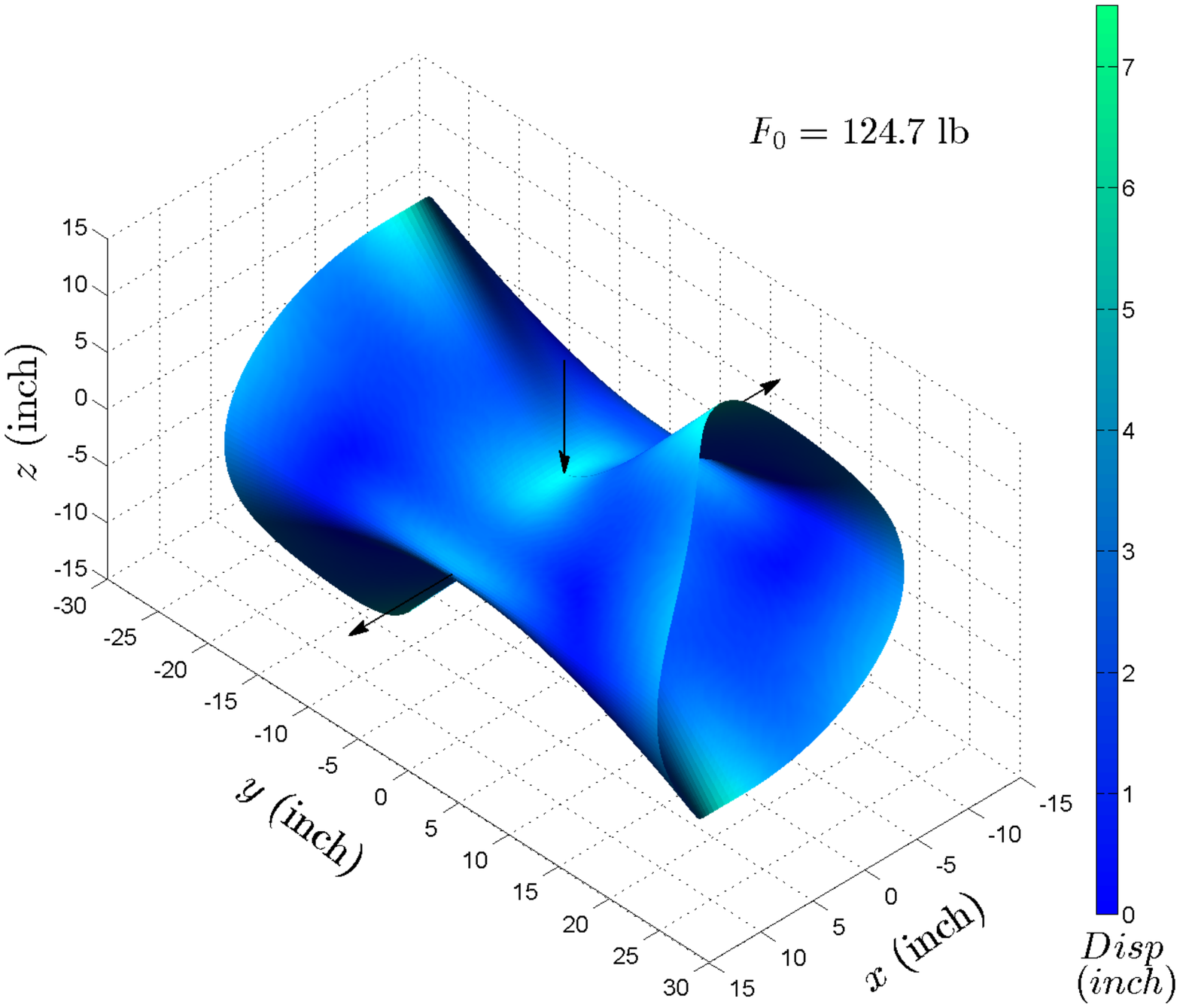}}
&\hskip-1cm \subfloat[ $F_0 = 563.0$ lb. ]{\includegraphics[trim=4cm 1cm 0.8cm 1cm,clip,scale=0.4]{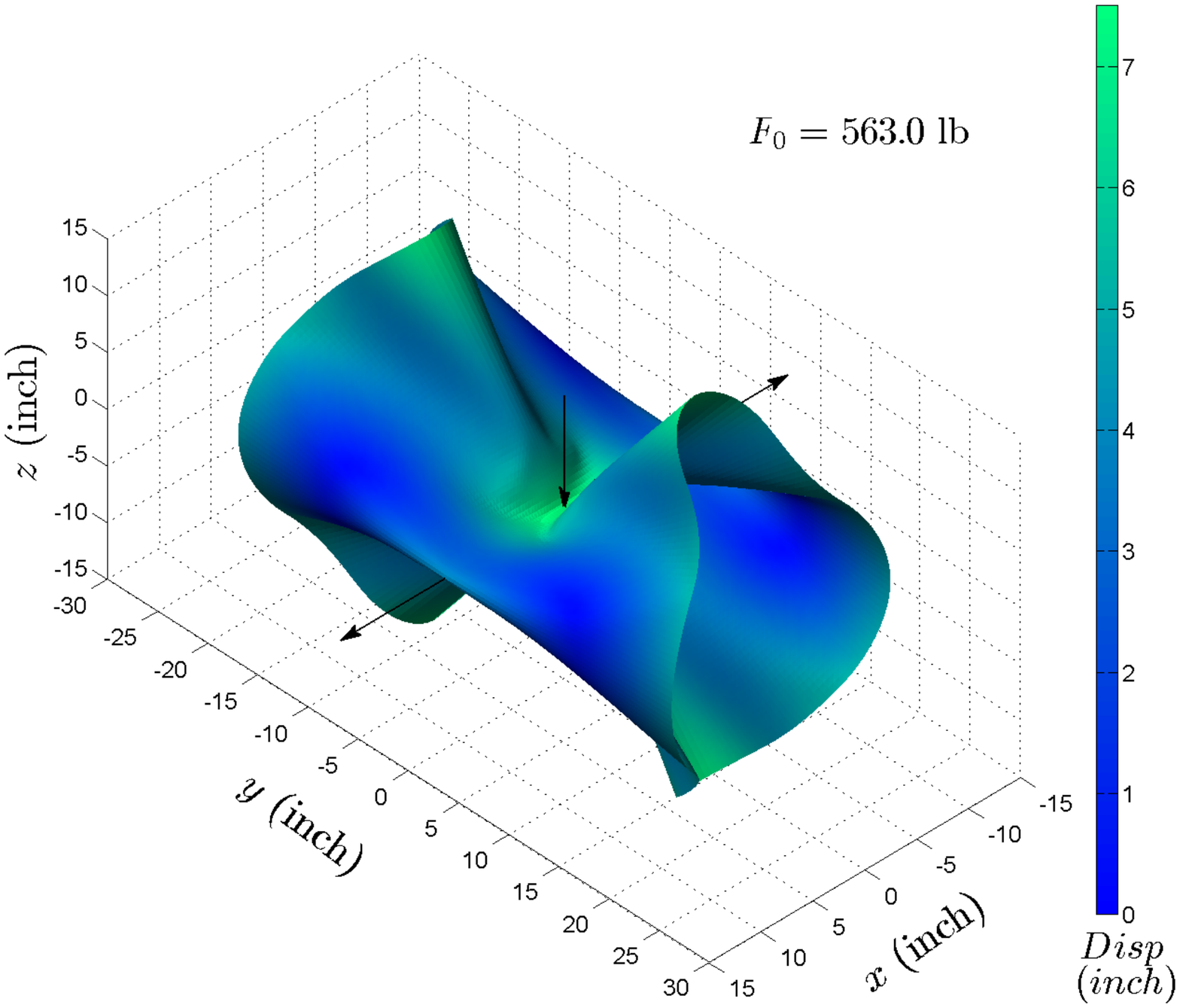}}
\end{tabular}
\end{adjustbox}
\caption{Deformed geometries of hyperboloidal shell subjected to point load considering neo-Hookean material model.}\label{fig:6}
\end{figure}

\begin{figure}[!htbp]
\centerline{\includegraphics[trim=0cm 1cm 0cm 0cm,clip,scale=0.45]{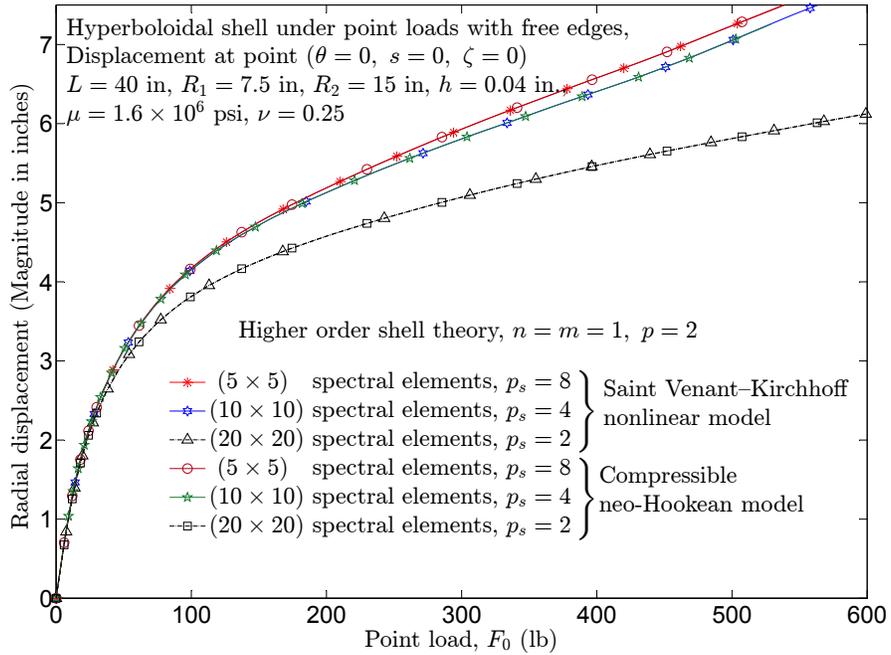}}
\caption{Load-displacement plot at the point of application of compressive point load. Here we note that higher-order spectral elements give the converged solution whereas the quadratic element gives an inaccurate solution for the large deformation and rotation. Also, the solution for compressible neo-Hookean material model and Saint Venant–Kirchhoff nonlinear material model are very similar as this problem is an example of small strain and large rotation.}
\label{fig:7}
\end{figure}

\noindent Figures~\ref{fig:6} (a) and (b) shows the deformed shape of the hyperboloidal shell for two different load cases, namely, $F_0=124.7$ lb and
$F_0=563.0$ lb, respectively, for the compressible neo-Hookean material model. Figure~\ref{fig:7} shows the magnitude of the radial displacement at the point of application of the compressive point load
versus the magnitude of the applied load $F_0$ for both compressible neo-Hookean and Saint Venant–Kirchhoff nonlinear material models. The displacements for both these nonlinear material models are almost the same. The reason for this could be that the considered shell being an example of large rotation and small strain. Moreover, for the small strain experienced, both compressible neo-Hookean and Saint Venant–Kirchhoff nonlinear material models tend to linearize to the same linear stress-strain relation for given material constants. Also, we can observe from Fig.~\ref{fig:7} that the quadratic element does not give accurate results for large deformation.

\subsection{Thin circular arc shaped shell-strip subjected to point load }
\begin{figure}[htbp]
\centerline{\includegraphics[trim=1cm 1cm 1cm 3cm,clip,scale=0.46]{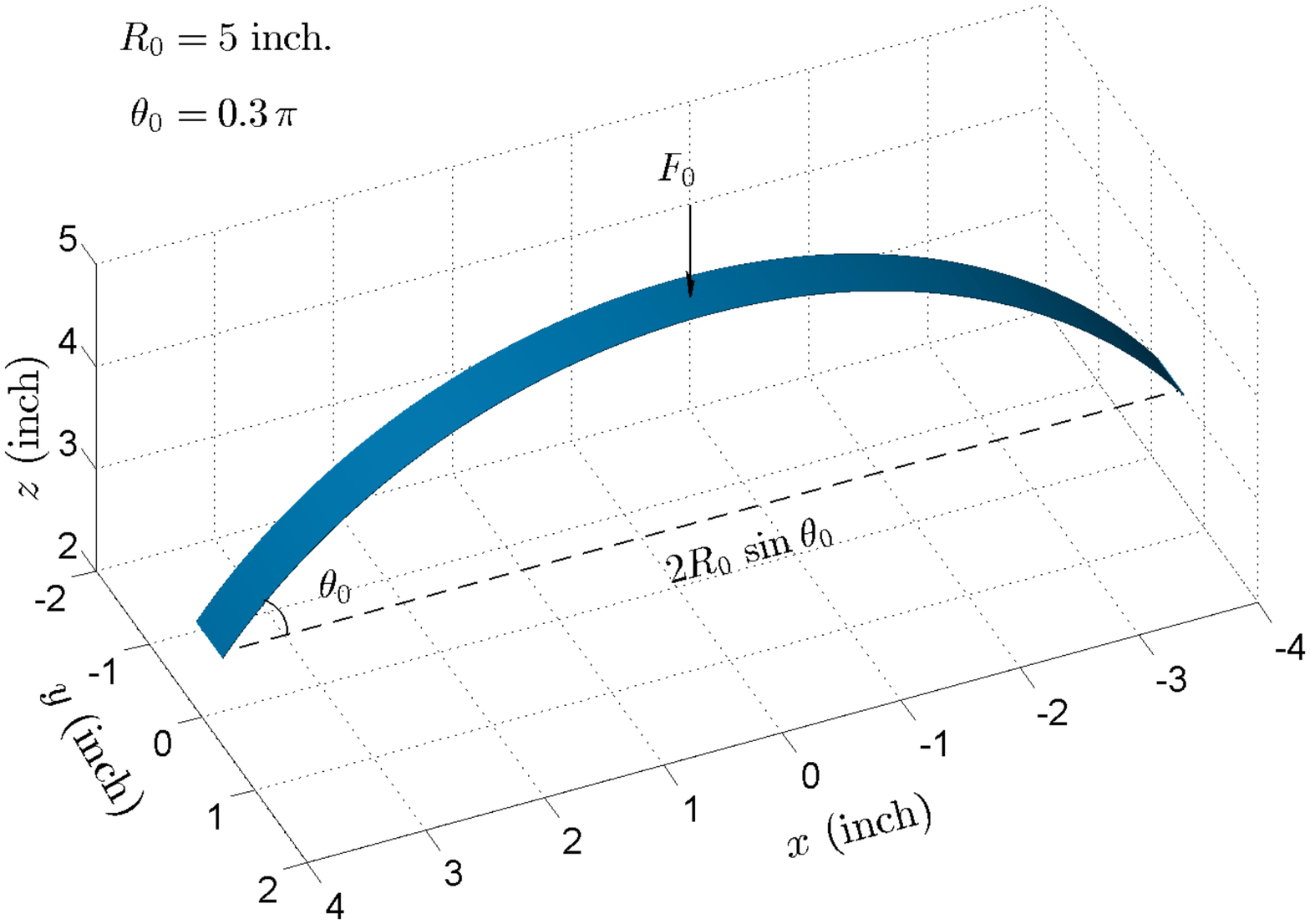}}
\caption{Original shape of thin circular arc shaped shell-strip subjected to point load.}
\label{fig:8}
\end{figure}

In this example, we analyze a circular arc-shaped thin shell strip under point load as shown in Fig.~\ref{fig:8} for two different boundary conditions; (i) Completely fixed (ii) hinged at the shorter edges of the shell stripe considering Saint Venant–Kirchhoff nonlinear material model.
This loading condition shows the snap-through motion of the structure, and hence the arc-length method has been applied to solve the nonlinear finite element equation. The geometric and material properties of the semi-cylinder are the following:
\begin{equation}\label{eq:65}
\theta_0 = 0.3\pi , \qquad R_0  = 5 \hbox{ in.} \qquad h  = 0.04 \hbox{ in.}
\end{equation}
\begin{equation}\label{eq:66}
E = 2.0685\times10^7  \hbox{ psi.}, \qquad \nu = 0.3
\end{equation}
where 2$\theta_0$ is the angle inscribed by the arc-length at the center of the circle of which this arc-length is part of and $R_0$ is the radius of the same circle. $h$ is the thickness of the shell-stripe. $E$ and $\nu$ are the modulus of elasticity and Poisson's ratio, respectively.

\begin{figure}[!htbp]
\centering
\begin{adjustbox}{center}
\begin{tabular}{cc}
\subfloat[ Fixed boundary condition ]{\includegraphics[trim=2cm 1cm 0.5cm 1cm,clip, scale=0.31]{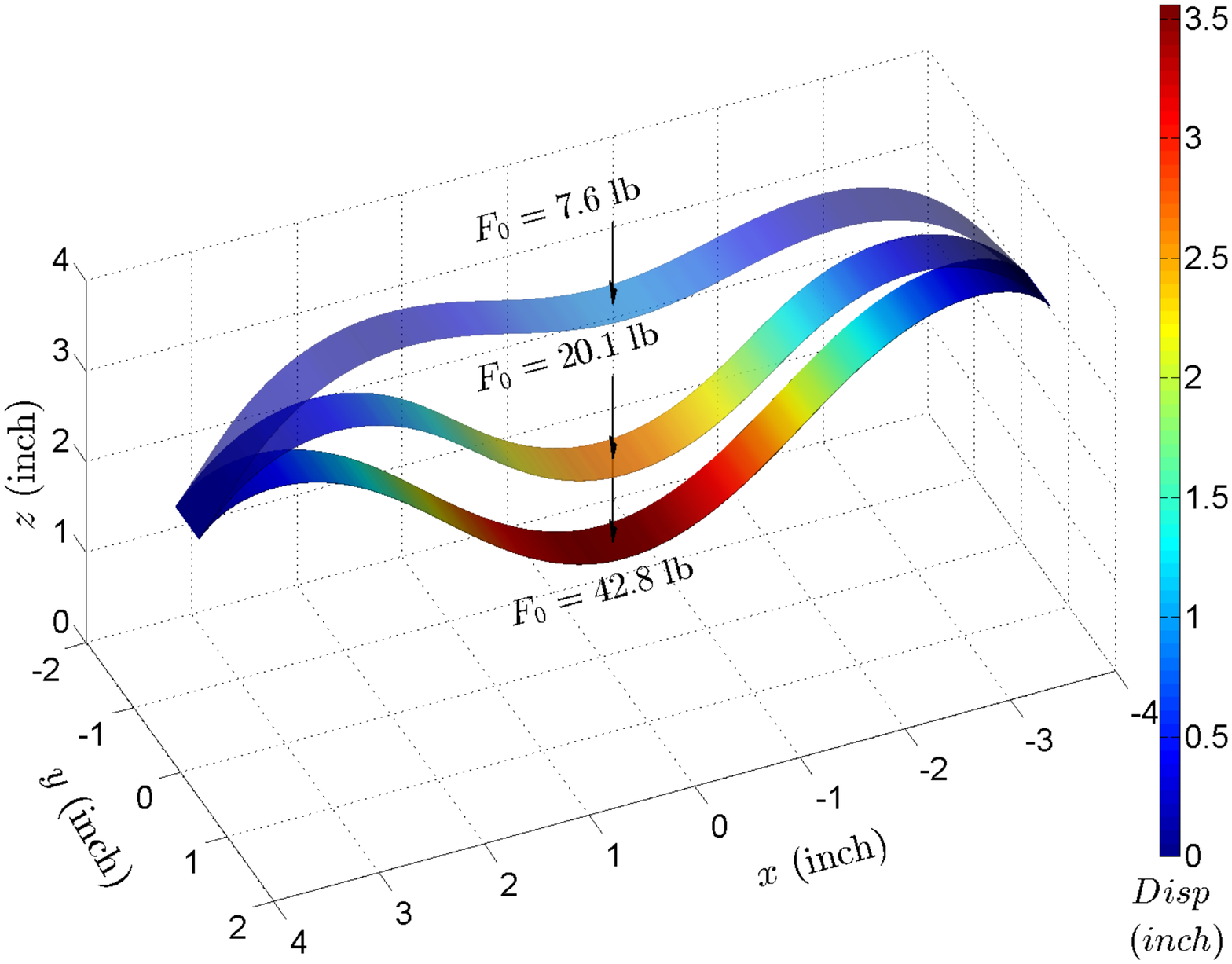}}
&\hskip-0.5cm \subfloat[ Hinged boundary condition ]{\includegraphics[trim=2cm 1cm 0.8cm 1cm,clip,scale=0.31]{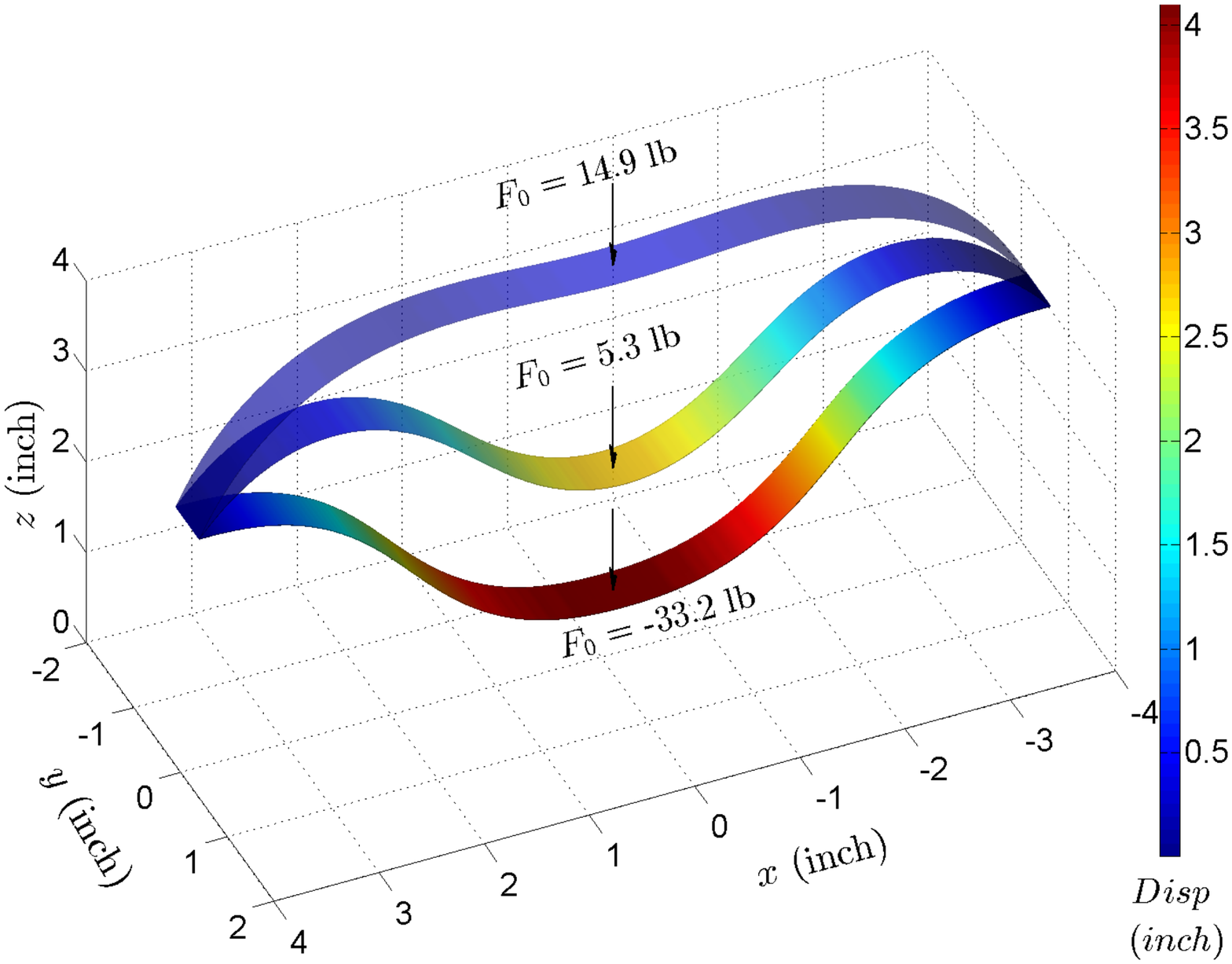}}
\end{tabular}
\end{adjustbox}
\caption{Deformed geometries of shell-strip subjected to point load for fixed and hinged boundary conditions.}\label{fig:9}
\end{figure}

\begin{figure}[!htbp]
\centering
\begin{adjustbox}{center}
\begin{tabular}{cc}
\subfloat[ Fixed boundaries]{\includegraphics[trim=0.8cm 1cm 0.5cm 1cm,clip, scale=0.31]{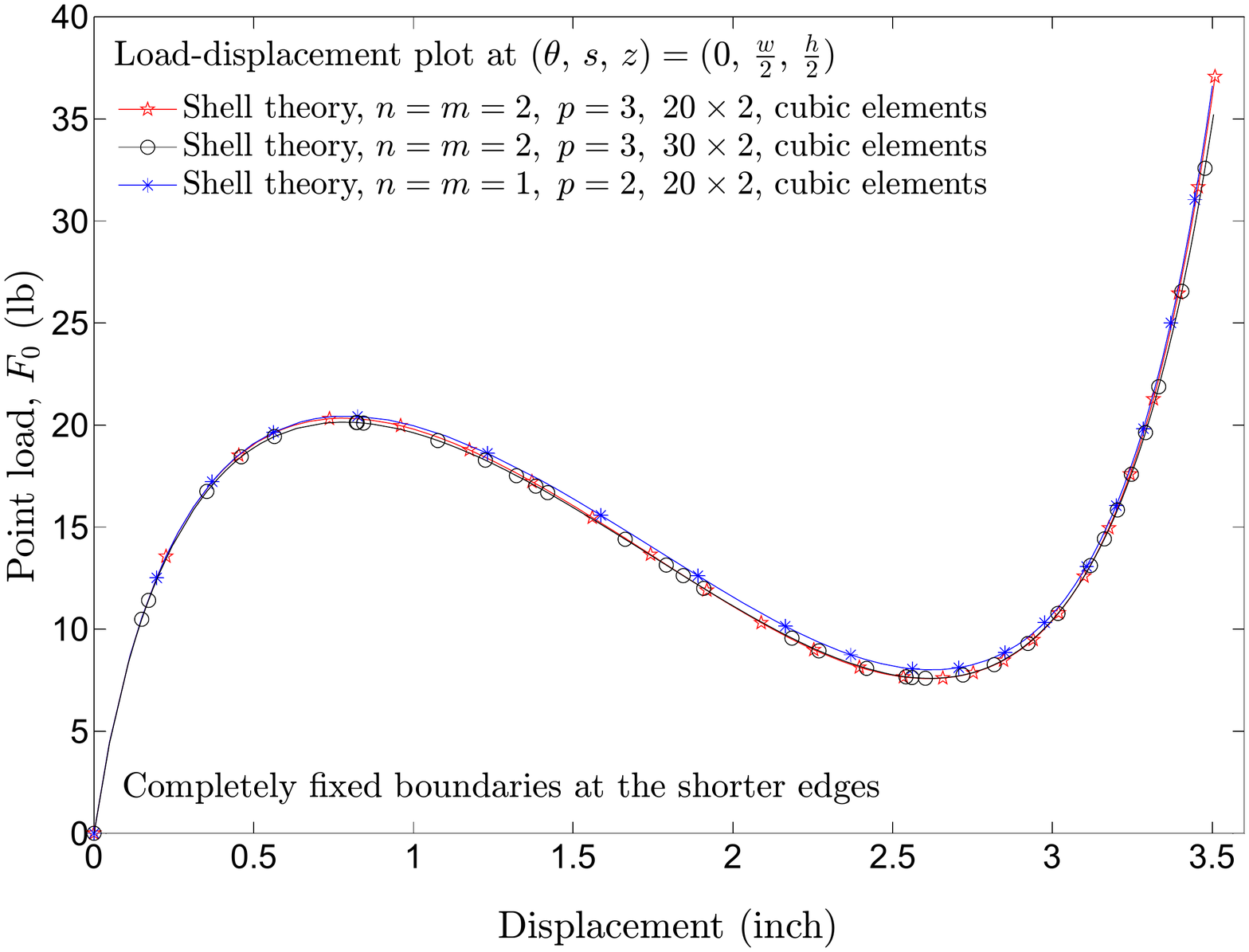}}
&\hskip-0.5cm \subfloat[ Hinged boundaries ]{\includegraphics[trim=0.8cm 1cm 0.8cm 1cm,clip,scale=0.31]{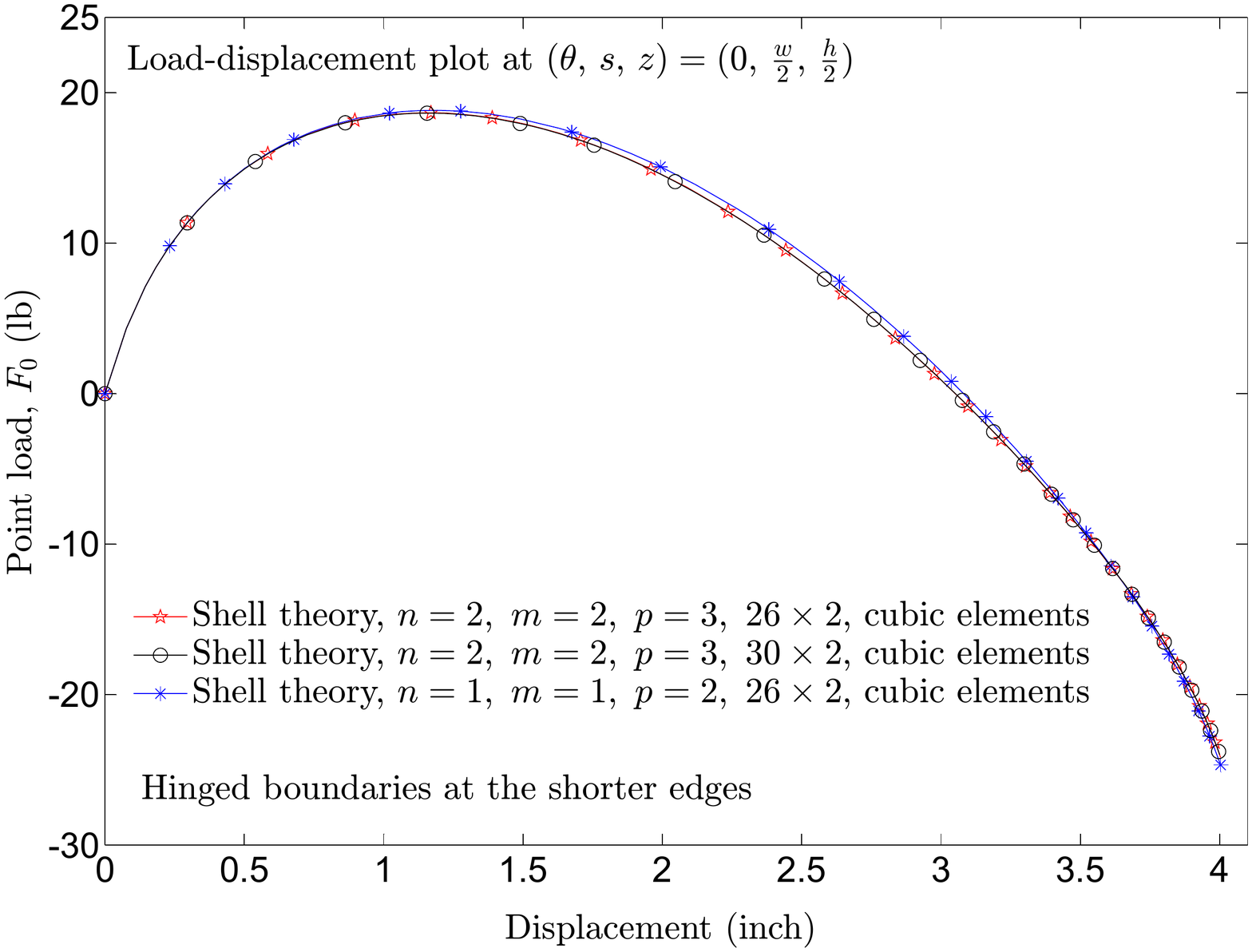}}
\end{tabular}
\end{adjustbox}
\caption{Load displacement curve (or the equilibrium path) at the point of load application for fixed and hinged boundaries, which shows the snap-through motion. We note here that the load-displacement plots are convergent for the successive refinement to reach the convergent equilibrium paths.}\label{fig:10}
\end{figure}

Three different combination of meshes and different orders of the shell theory have been used considering cubic spectral finite elements; (i) $(20 \times 2)$ (20 elements along arc-length and 2 elements along the width of the shell-stripe) cubic spectral elements, (ii) $(30 \times 10)$ cubic spectral elements with $n=m=1,\ p=2$ and $n=m=2,\ p=3$ as the order of approximation of higher-order shell theory. Full Gauss points ($4\times4$) are used to integrate the stiffness and tangent matrices. Figures \ref{fig:9}(a) and (b) shows the deformed shapes for three different loads along the equilibrium path for fixed and hinged boundary condition, respectively, whereas Figs. \ref{fig:10}(a) and (b) shows the load-displacement (at the point of load application) plot along the equilibrium path; here we note the snap-through motion for load-controlled scenario. Also, the solutions are convergent for refined mesh sizes and orders of the shell theory.

\subsection{Circular cylinder under internal pressure}
\begin{figure}[!htbp]
\centering
\begin{adjustbox}{center}
\begin{tabular}{cc}
\subfloat[ Natural configuration of the cylindrical shell ]{\includegraphics[trim=2cm 1cm 2.5cm 1cm,clip, scale=0.35]{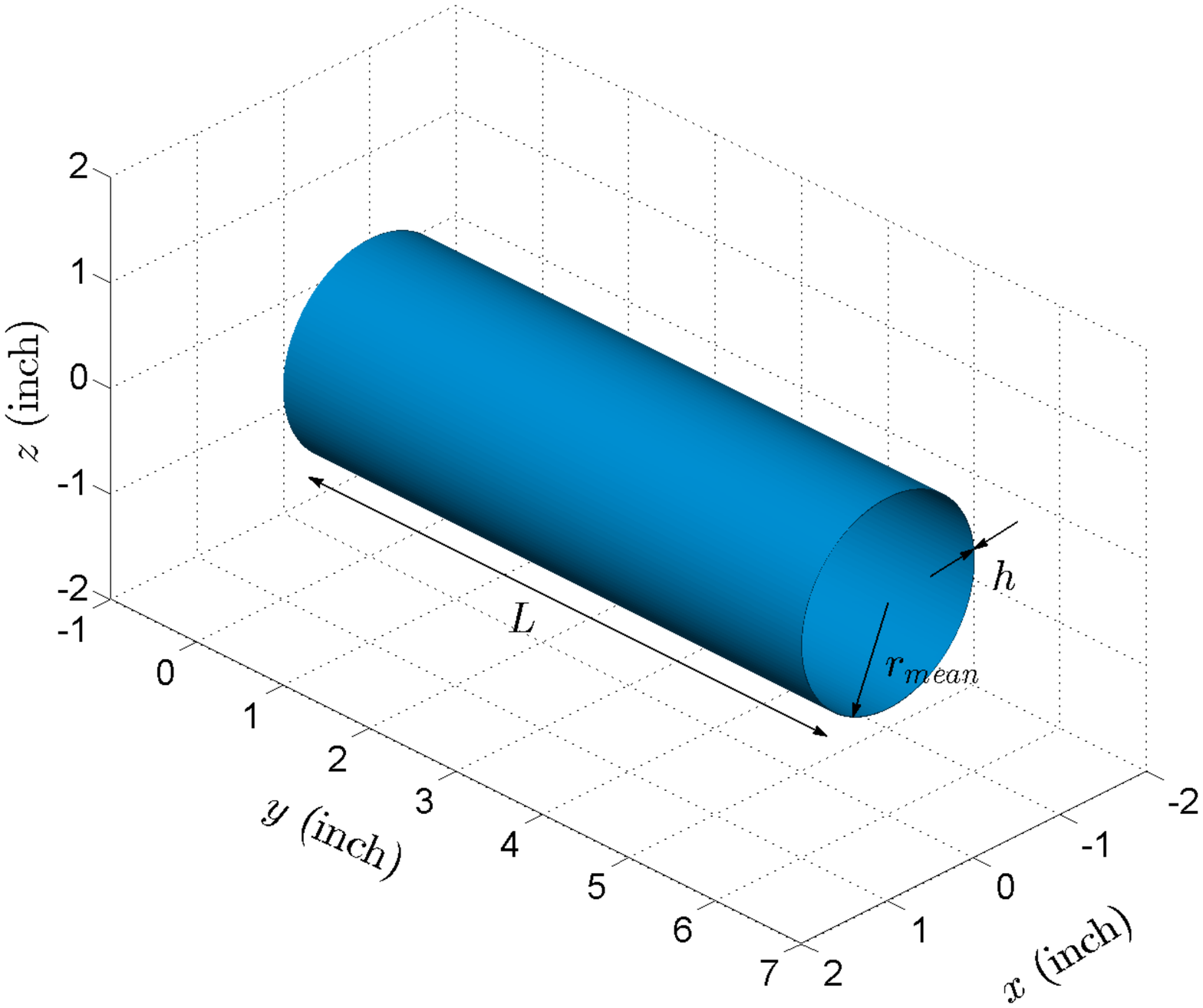}}
&\hskip-0.5cm \subfloat[ Deformed shape for $P_0 = 278.7 $ psi. ]{\includegraphics[trim=2cm 1cm 1cm 1cm,clip,scale=0.35]{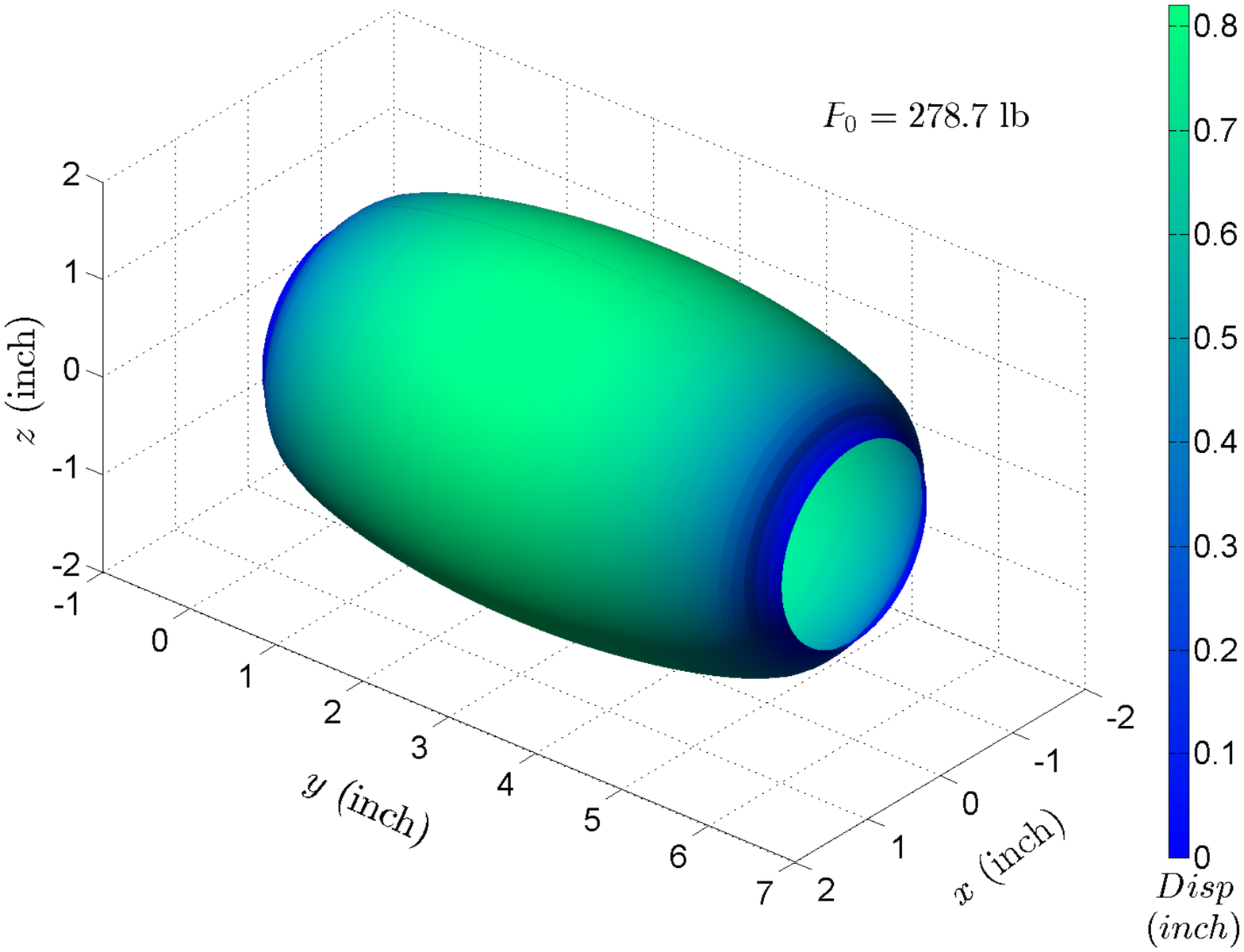}}
\end{tabular}
\end{adjustbox}
\caption{Natural and deformed geometries of cylindrical shell subjected to internal pressure.}\label{fig:11}
\end{figure}

In this example, a thin circular cylinder (see Fig. \ref{fig:11}(a)) subjected to internal pressure is analyzed considering compressible neo-Hookean and Saint
Venant–Kirchhoff nonlinear material models. Both ends of the cylinder are completely fixed. The geometric and material parameters  used are as follows:
\begin{equation}\label{eq:67}
L = 6 \hbox{ in.} , \qquad R  = 1\hbox{ in.} , \qquad h  = 0.01\hbox{ in.}
\end{equation}
\begin{equation}\label{eq:68}
\mu = 3.333 \times 10^4 \hbox{ psi.}, \qquad \nu = 0.3
\end{equation}
where $L$, $R$, and $h$ are the length, mean radius, and thickness, respectively, of the cylindrical shell. Moreover, $\mu$ and $\nu$ are the shear modulus and Poisson's ratio, respectively.

\noindent As in the case of hyperboloidal shell,  we exploit the symmetry and model only $1/8$th of the full domain as the computational domain for the analysis. Three different grid sizes are considered using quadratic elements for two different orders of approximation of the shell theory. Full Gauss quadrature is  used to integrate the stiffness and tangent matrices for elements.  In this example, we have used the arc-length method to solve the nonlinear finite element equations to have a program-controlled load increment.

\begin{figure}[!htbp]
\centering
{\includegraphics[trim=0.8cm 1cm 0.8cm 1cm,clip,scale=0.4]{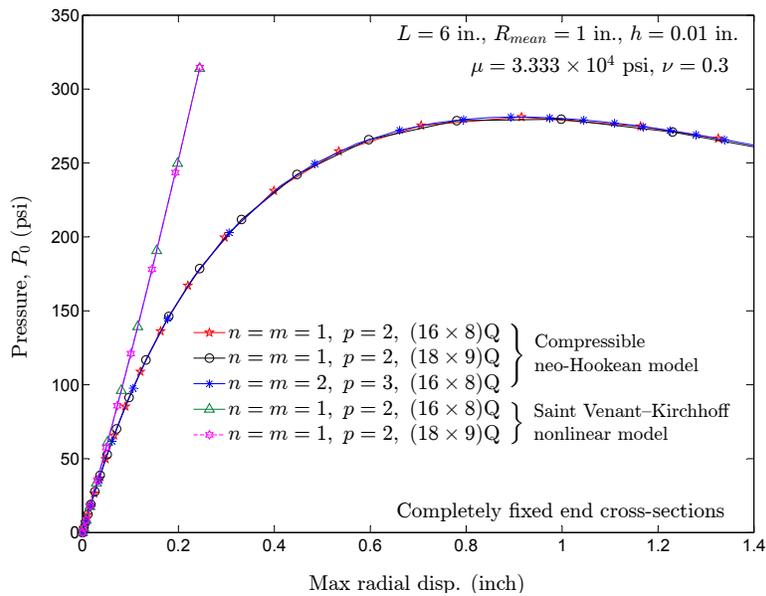}}
\caption{Load versus maximum radial displacement plot for the internally pressurized cylinder. Here, we notice that the plot coincides for small displacement {\it i.e.} small strain as both neo-Hookean and Saint Venant-Kirchhoff nonlinear models linearise to the same linear constitutive relation. Also, we note here that the load-displacement plots are convergent for the successive refinement of mesh and approximation order to reach the convergent equilibrium path. }\label{fig:12}
\end{figure}

Figure \ref{fig:11}(b) shows the deformed shapes for the internal pressure $P_0 = 278.7$ psi. Further, Fig. \ref{fig:12} shows the load vs. maximum radial displacement plot for both compressible neo-Hookean and Saint
Venant–Kirchhoff nonlinear material models. Also, in this case, the solutions are convergent for different mesh sizes and orders of shell theories for quadratic elements.

\subsection{Spiral tube under internal pressure}
\begin{figure}[htbp]
\centerline{\includegraphics[trim=1cm 1cm 1cm 2cm,clip,scale=0.46]{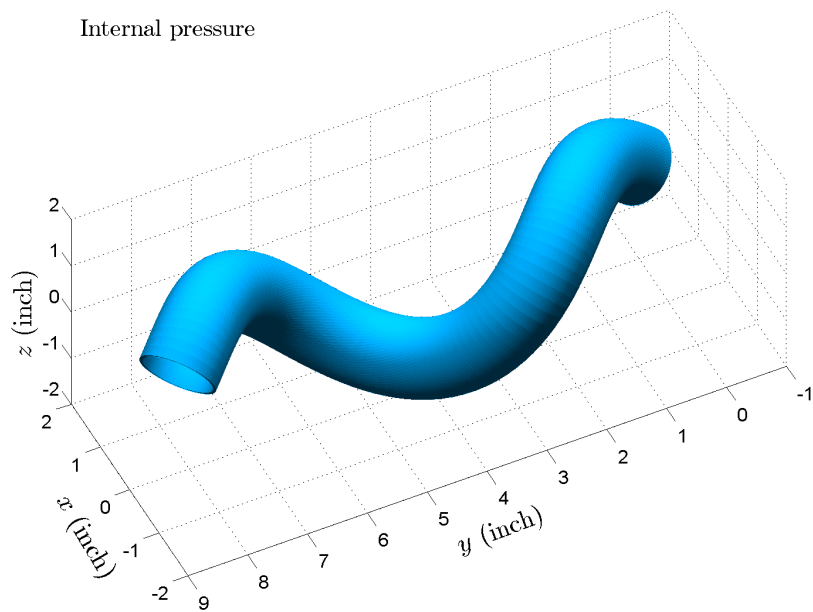}}
\caption{Original shape of spiral tube subjected to internal pressure.}
\label{fig:13}
\end{figure}

In this numerical example, we consider a spiral tube (see Fig.~\ref{fig:13}) under internal pressure with the central reference spiral curve given by following parametric equation in rectangular cartesian coordinate system:
\begin{equation}\label{eq:69}
 x =\cos\left({t\over \sqrt{2}}\right), \quad  y ={t\over \sqrt{2}}, \quad
 z =\sin\left({t\over \sqrt{2}}\right)
\end{equation}
where $t$ is an arbitrary parameter for this space curve and the arc length coordinate $s$ is same as the parameter $t$ for this particular case. Both the curvature ($\kappa$) and torsion ($\tau$) of the spiral curve are constant having the values 0.5 per inch and $-0.5$ per inch, respectively. The axis of the reference spiral curve lies along the $y$-axis. In this case also, two material models, namely, compressible neo-Hookean and  Saint Venant-Kirchhoff models are considered for the nonlinear finite element analysis. Other geometric and the material parameters are as follows:
\begin{equation}\label{eq:70}
L = 12 \hbox{ in.} , \qquad R_{mean}  = 0.3\hbox{ in.} , \qquad h  = 0.05\hbox{ in.}
\end{equation}
\begin{equation}\label{eq:71}
\mu = 3.333 \times 10^4 \hbox{ psi.}, \qquad \nu = 0.3
\end{equation}
where $L$, $R_{mean}$, and $h$ are the length, mean radius, and thickness, respectively, of the spiral tube. The curvilinear surface coordinate system are built over the underlying Frenet frame of the reference spiral curve of the tube, and hence the natural covariant basis of the surface coordinate will be non-orthogonal, but in this case we have considered the orthonormal basis for the analysis.

\begin{figure}[!htbp]
\centering
\begin{adjustbox}{center}
\begin{tabular}{cc}
\subfloat[ Deformed shape for $P_0 = 1762.8 $ psi ]{\includegraphics[trim=2cm 1cm 0cm 1cm,clip, scale=0.35]{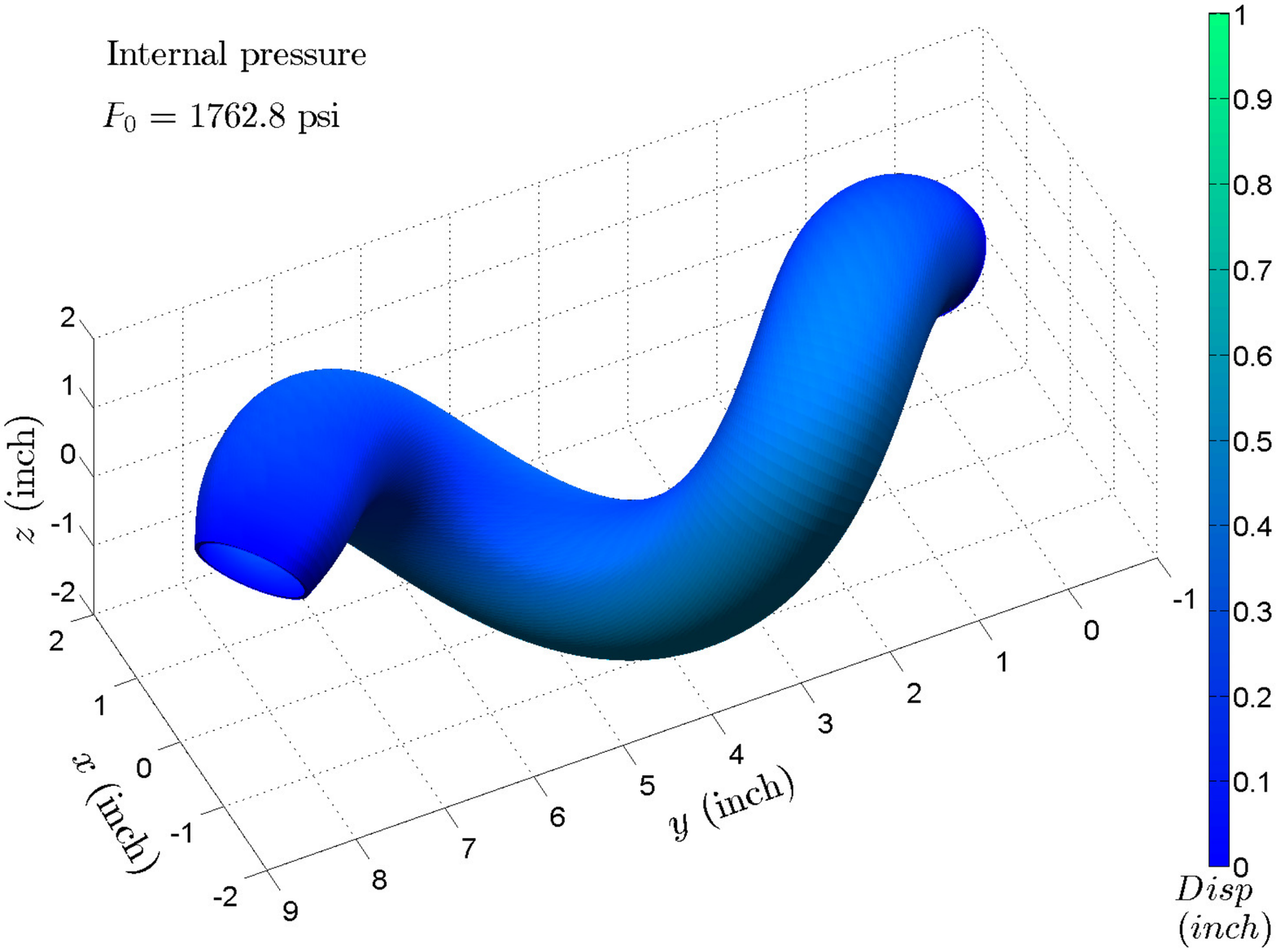}}
&\hskip-0.5cm \subfloat[ Deformed shape for $P_0 = 2026.1 $ psi. ]{\includegraphics[trim=2cm 1cm 0cm 1cm,clip,scale=0.35]{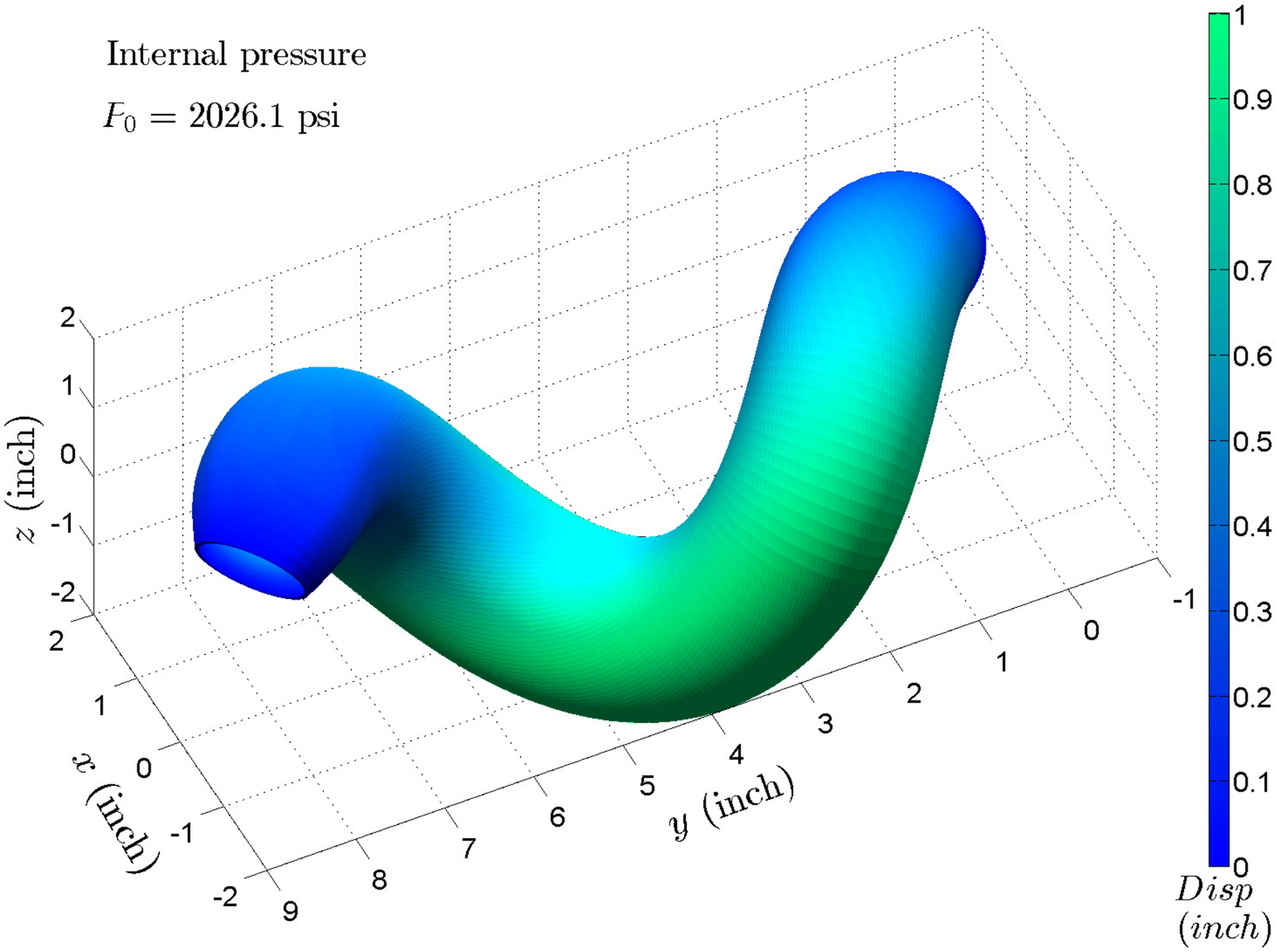}}
\end{tabular}
\end{adjustbox}
\caption{Deformed shape of the spiral tube subjected to internal pressure for compressible neo-Hookean material model.}\label{fig:14}
\end{figure}

\begin{figure}[!htbp]
\centering
\begin{adjustbox}{center}
\begin{tabular}{cc}
\subfloat[ Compressible neo-Hookean material model ]{\includegraphics[trim=2cm 1cm 0cm 1cm,clip, scale=0.35]{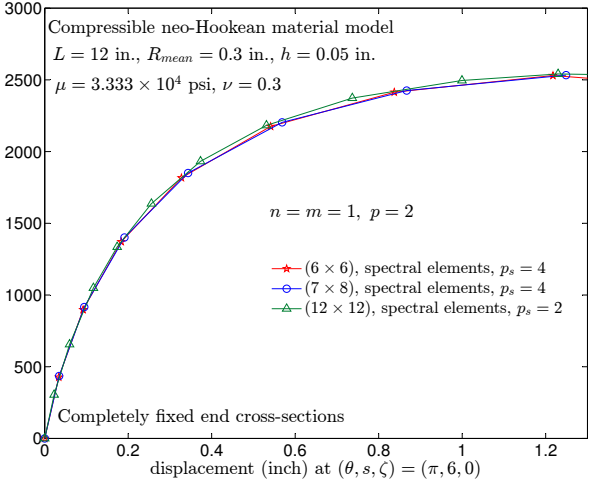}}
&\hskip-0.5cm \subfloat[ Saint Venant-Kirchhoff material model]{\includegraphics[trim=2cm 1cm 0cm 1cm,clip,scale=0.35]{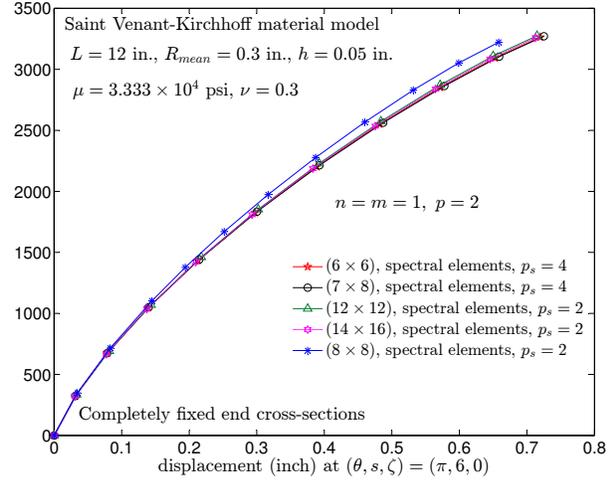}}
\end{tabular}
\end{adjustbox}
\caption{Load-displacement plot at $(\theta,s,\zeta)= (\pi,6,0)$ for the deformed spiral tube under internal pressure for the spectral elements of different order, mesh refinement and order of approximation of displacement field to see the convergence of solutions for different kinds of refinement in finite element implementations for both compressible neo-Hookean and saint Venant-Kirchhoff material models. It is observed that we obtain converged solution for successive refinement of the finite element mesh.}\label{fig:15}
\end{figure}

\begin{figure}[htbp]
\centerline{\includegraphics[trim=1cm 1cm 1cm 1cm,clip,scale=0.46]{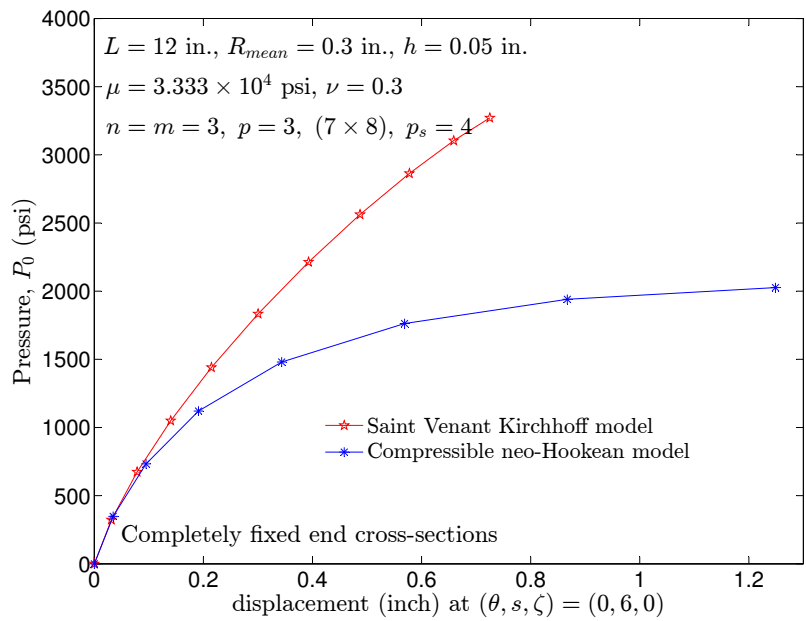}}
\caption{Load-displacement plot at $(\theta,s,\zeta)= (\pi,6,0)$ for the deformed spiral tube for compressible neo-Hookean and saint Venant-Kirchhoff material models. It is noted that for small deformation both material models undergo approximately similar deformation field.}
\label{fig:16}
\end{figure}

For the nonlinear finite element method, full domain of the shell mid-surface has been considered with spectral elements with order $p_s=2$ and $4$. $7$ -parameter shell theory has been considered to model the tube's thickness stretch. To solve the nonlinear finite element equation, the arc-length method has been considered; this way, we allow the program to decide the next load step while tracing the deformed tube's equilibrium path. Figures \ref{fig:14}(a) and (b) shows two different deformed configuration for neo-Hookean material model for internal pressure $P_0=1776.8$ psi. and $P_0=2026.1$ psi, respectively, considering $(7\times8)$ spectral element with $p_s=4$ for 7-parameter shell theory. Here, $7$ and $8$ elements are along the circumferential and longitudinal directions, respectively. In all cases, full Gauss points are used for integrating the stiffness, tangent matrix, and force vector. Figure \ref{fig:15} shows the load-displacement plot at $(\theta,s,\zeta) = (\pi,6,0)$ for the deformed tube for the spectral elements of a different order, mesh refinement and order of approximation of displacement field to see the convergence of solutions for different kinds of refinement in finite element implementations for both compressible neo-Hookean and saint Venant-Kirchhoff material models. We observe that spectral elements of order $p_s = 4$ give better convergent solution than the quadratic elements for a given number of nodes. Figure \ref{fig:16} shows a similar plot comparing the two nonlinear material models, namely, compressible neo-Hookean and Saint Venant-Kirchhoff model. For small deformation, both the models converge to similar deformations, as observed in the previous examples. The numerical values of the pressure and the displacement at $(\theta,s,\zeta) = (\pi,6,0)$ are tabulated in Table 2 which is calculated using the arc-length method.

\begin{table}[!htbp]\small
\centering
\caption{Magnitude of displacement at the point $(\theta,s,\zeta)= (\pi,6,0)$ or $(x,y,z)= (-0.724,4.242,-1.427)$ at the mid surface of the tube for the given internal pressure considering $(7 \times 8)$ spectral element of order $p_s = 4$ calculated via arc-length method for 7-parameter shell theory.}\label{Table:2}
\vskip 1mm
\begin{adjustbox}{center}
\begin{tabular}{|cc|cc|}
\hline
\multicolumn{2}{|c|}{Compressible neo-Hookean}&\multicolumn{2}{c|}{Saint Venant-Kirchhoff}\\
\multicolumn{2}{|c|}{material model}&\multicolumn{2}{c|}{material model}\\ \hline
Pressure (ksi)  &displacement (inch)&Pressure (ksi)&displacement (inch) \\
\hline
0.0000	&	0.0000	&	0.0000	&	0.0000	\\
0.3479	&	0.0346	&	0.3230	&	0.0314	\\
0.7328	&	0.0949	&	0.6738	&	0.0785	\\
1.1206	&	0.1910	&	1.0498	&	0.1401	\\
1.4803	&	0.3435	&	1.4407	&	0.2148	\\
1.7628	&	0.5688	&	1.8338	&	0.3003	\\
1.9403	&	0.8671	&	2.2135	&	0.3928	\\
2.0261	&	1.2488	&	2.5624	&	0.4872	\\
\hline
\end{tabular}
\end{adjustbox}
\end{table}
\section{Summary and conclusions}
In this study, we introduced a novel general higher-order shell theory for the general compressible hyperelastic material model using an orthonormal moving frame.  The geometry of the shell surface has been represented exactly, and the governing equation and its finite element model have been derived in terms of the surface coordinates and orthonormal bases. A novel method to obtain the kinematic invariants on the curvilinear coordinates with the orthonormal moving frame is presented in the appendix of this study.  Also, the displacement field along the normal line of the shell reference surface has been approximated by general Taylor series or Legendre polynomial in normal coordinate $\zeta$, which makes the general higher-order shell theory suitable for analysis of both thin or thick shell structures by just changing the approximation orders of the displacement field. This way, we can model thickness stretch to desired order to avoid any numerical locking in the shell thickness direction. In the finite element model, we have used the higher-order spectral element to avoid membrane locking.   Further, we also present the specialization of various compressible hyperelastic material models such as compressible neo-Hookean, Saint-Venant -Kirchhoff model, and Mooney-Rivlin model. Various numerical examples have been presented to illustrate the use of various orders of shell theories using Newton's or arc-length methods. We have shown that the solutions converge for successively refined meshes or orders of displacement approximations. Thus, in this study, it is shown how such an approach could result in a more efficient and accurate computational shell theory compared to the existing computational shell theories.

The general higher-order shell theory presented herein has the potential for extension (due to its orthonormal basis framework) in various other complex constitutive relations of the incompressible hyperelastic material model, plasticity, stretch dependent hyperelasticity, viscoelasticity, viscoplasticity or other implicit constitutive relation which can be carried out as future work. Applications of the present model to biological systems to gain insights into their functionalities are awaiting attention.

\section*{Acknowledgments}
The authors are grateful for the financial support provided by the Oscar S. Wyatt Endowed Chair in the Department of Mechanical Engineering at Texas A\&M University.

\renewcommand{\theequation}{A.\arabic{equation}}
\renewcommand{\thefigure}{A.\arabic{figure}}
\renewcommand{\thesection}{A}
\setcounter{equation}{0}
\setcounter{figure}{0}
\setcounter{section}{0}
\section*{Appendix A: Gradient for various curvilinear coordinate systems}
\subsection{Gradient for curvilinear cylindrical coordinate system with orthonormal moving frame on the reference surface of the shell-structure}
\begin{figure}[htbp]
\centering{\includegraphics[trim=0cm 1cm 0cm 1cm,clip,scale=0.50]{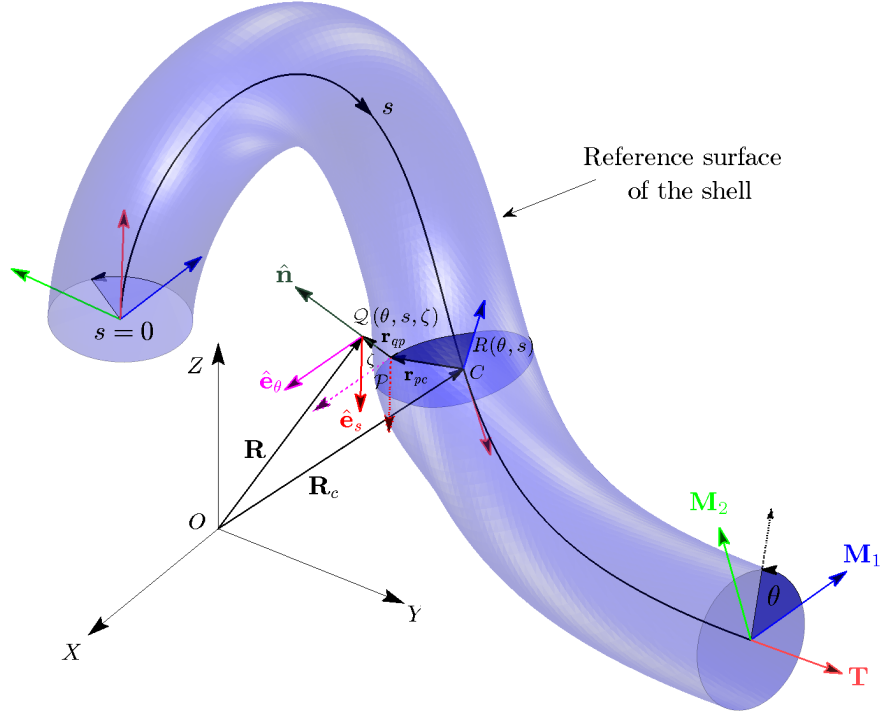}}
\caption{Arbitrary surface in curvilinear cylindrical coordinate system.}
\label{fig:A1}
\end{figure}

Let us consider a curved cylindrical pipe-like reference surface with arbitrary varying cross-section embedded in three-dimensional Euclidean space, $\mathbb{R}^3$, as shown in Fig.~\ref{fig:A1}.
We consider a space curve, $\mathcal{C}$, as the axis of the reference surface. We frame the reference curve by a general hybrid frame $({\bf T},\,{\bf M}_1,\, {\bf M}_2)$ which varies along the curve according to the following rule:
\begin{eqnarray}\label{eq:A1}
{d\over ds}\left\{\begin{matrix}
         {\bf T} \\
         {\bf M}_1 \\
         {\bf M}_2
       \end{matrix}\right\}=\left[\begin{matrix}
                                     0 & \kappa_1 & \kappa_2 \\
                                     -\kappa_1 & 0 & \kappa_3 \\
                                     -\kappa_2 & -\kappa_3 & 0
                                   \end{matrix}\right]\left\{\begin{matrix}
         {\bf T} \\
         {\bf M}_1 \\
         {\bf M}_2
       \end{matrix}\right\}
\end{eqnarray}
Here $\kappa_1,\, \kappa_2$ and $\kappa_3$ are the components of the Cartan matrix of the moving frame (see Arbind {\it et al.} \cite{arbind2018curved} for details). Now we consider a curvilinear cylindrical coordinate system with surface coordinates, $\eta_1=\theta$ and $\eta_2=s$ and normal coordinate $\zeta$; $s$ is the arc-length coordinate measured along the reference curve $\mathcal{C}$, $\theta$ is the angle measure from the ${\bf M}_1$ towards ${\bf M}_2$ (see Fig.\ref{fig:A1}) and $\zeta$ is the normal coordinate in the direction of the normal to the reference surface. At any cross-section, the distance of any arbitrary point $\mathcal{P}$ on the reference surface from the point $C$ at the reference curve $\mathcal{C}$ is defined as $R(\theta,s)$ which is given for the geometry of the reference surface. The basis vectors at any arbitrary point $\mathcal{Q}$, whose position vector is ${\bf R}$, are defined as follows:
\begin{eqnarray}\label{eq:A2}
\hat{\bf e}_{\zeta} =\hat{\bf n}={{\bf R}_{,\zeta}\over||{\bf R}_{,\zeta}||},\quad
\hat{\bf e}_{\theta}=\hat{\bf e}_{\eta_1} ={{\bf R}_{,\eta_1}\over||{\bf R}_{,\eta_1}||},\quad \hbox{and}\quad \hat{\bf e}_{s}=\hat{\bf e}_{\eta_2}=\hat{\bf e}_{\zeta}\times\hat{\bf e}_{\eta_1}
\end{eqnarray}
where $\hat{\bf n}$ is the unit normal to the reference surface. For any arbitrary tubular reference surface. The orthonormal bases can be given as following (see Appendix B of \cite{arbind2018Neo-Hookean} for detail derivation of $\hat{\bf n}$ and $\hat{\bf e}_\theta$):
\begin{align}\label{eq:A3}
\hat{\bf e}_\theta &={1\over\alpha_5}\left[\alpha_3{\bf M}_1+\alpha_2{\bf M}_2\right]\cr
\hat{\bf e}_s &={1\over\alpha_5\alpha_4}\left[\xi \alpha_5^2{\bf T}-R\alpha_1\alpha_2{\bf M}_1+R\alpha_1\alpha_3{\bf M}_2 \right]\cr
{\hat {\bf n}} &= {1\over\alpha_4}\left[R\alpha_1{\bf T} + \xi\alpha_2{\bf M}_1 -\xi\alpha_3{\bf M}_2\right]
\end{align}
where
\begin{align}\label{eq:A4}
\alpha_1 &= \left(\kappa_3 {\partial R \over \partial \theta}- {\partial R \over \partial s} \right),\quad \alpha_2 = \left({\partial R \over \partial \theta}\, \sin\theta+ R\,\cos\theta\right),\quad\alpha_3=\left({\partial R \over \partial \theta}\, \cos\theta- R\,\sin\theta\right)\cr
\alpha_4 &= \sqrt{R^2\alpha_1^2+\xi^2(\alpha_2^2+\alpha_3^2)},\quad \alpha_5=\sqrt{\alpha_2^2+\alpha_3^2}=\sqrt{\left({\partial R \over \partial \theta}\right)^2+R^2}
\end{align}
and
\begin{equation}\label{eq:A5}
\xi = 1-R\hat{\kappa}_1,\quad \hat{\kappa}_1 = \kappa_1 \cos\theta\,+\kappa_2\sin\theta,\quad \hat{\kappa}_2 = \kappa_1 \sin\theta\,-\kappa_2\cos\theta
\end{equation}
Then the curved frame $({\bf T},\, {\bf M}_1,\, {\bf M}_2)$ can be given in term of the assumed orthonormal bases as following:
\begin{align}\label{eq:A6}
{\bf T}   &= \xi{\alpha_5\over \alpha_4}\hat{\bf e}_s+R{\alpha_1\over \alpha_4}\hat{\bf n}\cr
{\bf M}_1 &={\alpha_3\over \alpha_5}\hat{\bf e}_\theta-R{\alpha_1\alpha_2\over \alpha_4\alpha_5}\hat{\bf e}_s+\xi{\alpha_2\over \alpha_4}\hat{\bf n} \cr
{\bf M}_2 &={\alpha_2\over \alpha_5}\hat{\bf e}_\theta+R{\alpha_1\alpha_3\over \alpha_4\alpha_5}\hat{\bf e}_s-\xi{\alpha_3\over \alpha_4}\hat{\bf n}
\end{align}
Next, we obtain the exterior derivative of the considered basis vectors  as follows:
\begin{equation}\label{eq:A7}
\begin{aligned}
&d\hat{\bf e}_\theta &=&&&&(a_{1}\,d\theta +b_{1}\,ds)\,\hat{\bf e}_s&& +(a_{2}\,d\theta +b_{2}\,ds)\,\hat{\bf n}&\cr
&d\hat{\bf e}_s &=&&-(a_{1}\,d\theta+ b_{1}\,ds)\,\hat{\bf e}_\theta &&&&+(a_{3}\,d\theta+ b_{3}\,ds)\,\hat{\bf n}&\cr
&d\hat{\bf n} &=&&-(a_{2}\,d\theta +b_{2}\,ds)\,\hat{\bf e}_\theta&&-(a_{3}\,d\theta +b_{3}\,ds)\,\hat{\bf e}_s &&&
\end{aligned}
\end{equation}
where
\begin{align}\label{eq:A8}
a_{1} & = {R\alpha_1\over \alpha_4\alpha_5^2}\left(R^2+2\left({\partial R\over \partial \theta}\right)^2-R{\partial^2 R\over \partial \theta^2}\right)\cr
a_{2} & = {\xi\over \alpha_4\alpha_5}\left(\alpha_2{\partial \alpha_3\over \partial \theta}-\alpha_3{\partial\alpha_2\over \partial \theta}\right)\cr
a_{3} & = {\alpha_1\alpha_5\over\alpha_4^2}\left(R^2\hat{\kappa}_2-{\partial R\over\partial\theta}\right)+{R\xi\over\alpha_4^2}\left(\alpha_1{\partial \alpha_5 \over\partial\theta}-\alpha_5{\partial \alpha_1 \over\partial\theta}\right)\cr
b_{1} & ={R\kappa_3\alpha_1\over\alpha_4}-{\xi\over \alpha_4}(\hat{\kappa}_1{\partial R\over\partial\theta}-R\hat{\kappa}_2)+{R\alpha_1\over \alpha_4\alpha_5^2}\left({\partial R\over \partial s}{\partial R\over \partial \theta}-R{\partial^2 R\over \partial s\partial \theta}\right) \cr
b_{2} & =  -{\xi\kappa_3\alpha_5\over\alpha_4}-{R\alpha_1\over \alpha_4\alpha_5}(\hat{\kappa}_1{\partial R\over\partial\theta}-R\hat{\kappa}_2)+{\xi\over \alpha_4\alpha_5}\left(\alpha_2{\partial \alpha_3\over \partial s}-\alpha_3{\partial\alpha_2\over \partial s}\right)\cr
b_{3} & = {1\over \alpha_5}\left(R\hat{\kappa}_1+{\partial R\over \partial \theta}\hat{\kappa}_2\right)-{\alpha_1\alpha_5\over\alpha_4^2}\left(\xi{\partial R\over\partial s}+R{\partial \xi\over\partial s}\right)-{R\xi\over\alpha_4^2}\left(\alpha_5{\partial \alpha_1 \over\partial s}-\alpha_1{\partial \alpha_5 \over\partial s}\right)
\end{align}
The position vector of any arbitrary point $\mathcal{Q}$ is
\begin{eqnarray}\label{eq:A9}
{\bf R}={\bf R}_c+{\bf r}_{pc}+{\bf r}_{qp}={\bf R}_c + R\, \cos\theta\, {\bf M}_1 + R\, \sin\theta\, {\bf M}_2 + \zeta\,\hat{\bf n}
\end{eqnarray}
Further, the length element, which is a vector valued one form, is
given as
\begin{align}\label{eq:A10}
d{\bf R}&=\left[\left(\alpha_5-a_2\zeta\right)d\theta+\left({1\over \alpha_5} {\partial R\over\partial s}{\partial R\over\partial \theta}+\kappa_3{R^2\over\alpha_5}-b_2\zeta\right)ds\right]\hat{\bf e}_\theta\cr&+\left[-a_3\zeta\,d\theta+\left({\alpha_4\over \alpha_5}-b_3\zeta\right)ds\right]\hat{\bf e}_s+d\zeta\, \hat{\bf n}
\end{align}
Let us write the length element $d{\bf R}$ as the following column vector\footnote{In the column vector, the first element is component along $\hat{\bf e}_\theta$ and the second and third elements are the components along $\hat{\bf e}_s$ and $\hat{\bf n}$, respectively. }:
\begin{align}\label{eq:A11}
\{d{\bf R}\}&=\left\{\begin{matrix}
\left(\alpha_5-a_2\zeta\right)d\theta+\left({1\over \alpha_5} {\partial R\over\partial s}{\partial R\over\partial \theta}+\kappa_3{R^2\over\alpha_5}-b_2\zeta\right)ds  \\
-a_3\zeta\,d\theta+\left({\alpha_4\over \alpha_5}-b_3\zeta\right)ds \\
                   d\zeta
                 \end{matrix}\right\}
\end{align}
Further the volume element $dV$, which is scaler valued three form, can be obtained as following:
\begin{align}\label{eq:A12}
dV&=\left(\left(\alpha_5-a_2\zeta\right)d\theta+\left({1\over \alpha_5} {\partial R\over\partial s}{\partial R\over\partial \theta}+\kappa_3{R^2\over\alpha_5}-b_2\zeta\right)ds\right)\wedge\left(-a_3\zeta\,d\theta+\left({\alpha_4\over \alpha_5}-b_3\zeta\right)ds\right)\wedge d\zeta\cr
&=\left[\left(\alpha_5-a_2\zeta\right)\left({\alpha_4\over \alpha_5}-b_3\zeta\right)
+a_3\zeta\left({1\over \alpha_5} {\partial R\over\partial s}{\partial R\over\partial \theta}+\kappa_3{R^2\over\alpha_5}-b_2\zeta\right)\right]d\theta\wedge ds\wedge d\zeta
\end{align}
where $(d\theta \wedge ds)$ means wedge product (see
Flanders \cite{flanders1963differential}) or exterior product between $d\theta$ and $ds$ and so on.
Dropping the wedge sign (and maintaining the order $d\theta$, $ds$, and $d\zeta$), we can rewrite the volume element as
\begin{align}\label{eq:A13}
dV&=g\, d\theta\, ds\, d\zeta
\end{align}
where
\begin{align}\label{eq:A14}
g = \left(\alpha_5-a_2\zeta\right)\left({\alpha_4\over \alpha_5}-b_3\zeta\right)
+a_3\zeta\left({1\over \alpha_5} {\partial R\over\partial s}{\partial R\over\partial \theta}+\kappa_3{R^2\over\alpha_5}-b_2\zeta\right)
\end{align}
Next, from Eq. (A.12), we have the following:
\begin{align}\label{eq:A15}
\left\{\begin{matrix}
                   d\theta\\
                   ds\\
                   d\zeta
                 \end{matrix}\right\}&=\left[\begin{matrix}
                   c_1&c_2&0 \\
                   c_3&c_4&0  \\
                   0&0&1
                 \end{matrix}\right]\left\{d{\bf R}\right\}
\end{align}
where
\begin{equation}\label{eq:A16}
\begin{aligned}
& c_1&=&{1\over g}\left({\alpha_4\over \alpha_5}-b_3\zeta\right), & &c_3&=&{1\over g}a_3\zeta &\cr
& c_2&=&\,-{1\over g}\left({1\over \alpha_5} {\partial R\over\partial s}{\partial R\over\partial \theta}+\kappa_3{R^2\over\alpha_5}-b_2\zeta\right), & &c_4&=& \,{1\over g}\left(\alpha_5-a_2\zeta\right)&
\end{aligned}
\end{equation}

Now, let us consider the displacement vector as ${\bf u}=u_\theta\hat{\bf e}_\theta+u_s \hat{\bf e}_s+u_\zeta\hat{\bf n}$, then we have the differential, $d{\bf u}:=({\boldsymbol\nabla}{\bf u})d{\bf R}$, where ${\boldsymbol\nabla}{\bf u}$ is the gradient of displacement vector. Further, using Eq.~\eqref{eq:A1}, $d{\bf u}$ can be given as
\begin{align}\label{eq:A17}
 d{\bf u} &= d{\bf u}_{rel}+u_\theta\,d\hat{\bf e}_\theta+u_s\, d\hat{\bf e}_s+u_\zeta\,d\hat{\bf n}\cr
 &= d{\bf u}_{rel}+u_\theta\,((a_{1}\,d\theta+b_{1}\,ds)\,\hat{\bf e}_s +(a_{2}\,d\theta+b_{2}\,ds)\,\hat{\bf n})\cr
& +u_s\, (-(a_{1}\,d\theta +b_{1}\,ds)\,\hat{\bf e}_\theta +(a_{3}\,d\theta +b_{3}\,ds)\,\hat{\bf n})\cr
 &+u_\zeta\,(-(a_{2}\,d\theta+ b_{2}\,ds)\,\hat{\bf e}_\theta-(a_{3}\,d\theta + b_{3}\,ds)\,\hat{\bf e}_s)
\end{align}
which can be again written as a column vector as
\begin{align}\label{eq:A18}
 d{\bf u} &= d{\bf u}_{rel}+\left[\begin{matrix}
                                    -a_1u_s-a_2u_\zeta & -b_1u_s-b_2u_\zeta&0 \\
                                    a_1u_\theta-a_3u_\zeta & b_1u_\theta-b_3u_\zeta&0 \\
                                    a_2u_\theta+a_3u_s & b_2u_\theta+b_3u_s&0
                                  \end{matrix}\right]\left\{\begin{matrix}
                   d\theta\\
                    ds \\
                    d\zeta
                  \end{matrix}\right\}
\end{align}
and
\begin{align}\label{eq:A19}
d{\bf u}_{rel} &=\left[\begin{matrix}
u_{\theta,\theta}&u_{\theta,s}&u_{\theta,\zeta}\\
u_{s,\theta}&u_{s,s}&u_{s,\zeta}\\
u_{\zeta,\theta}&u_{\zeta,s}&u_{\zeta,\zeta}
                                  \end{matrix}\right]\left\{\begin{matrix}
                   d\theta\\
                    ds\\
                    d\zeta
                  \end{matrix}\right\}
\end{align}
Hence
\begin{align}\label{eq:A20}
 d{\bf u} &= \left[\begin{matrix}
                                    u_{\theta,\theta}-a_1u_s-a_2u_\zeta & u_{\theta,s}-b_1u_s-b_2u_\zeta&u_{\theta,\zeta} \\
                                    u_{s,\theta}+a_1u_\theta-a_3u_\zeta & u_{s,s}+b_1u_\theta-b_3u_\zeta&u_{s,\zeta} \\
                                    u_{\zeta,\theta}+a_2u_\theta+a_3u_s & u_{\zeta,s}+b_2u_\theta+b_3u_s&u_{\zeta,\zeta}
                                  \end{matrix}\right]\left\{\begin{matrix}
                   d\theta\\
                    d s\\
                    d\zeta
                  \end{matrix}\right\}\cr
&=\left[\begin{matrix}
                                    u_{\theta,\theta}-a_1u_s-a_2u_\zeta & u_{\theta,s}-b_1u_s-b_2u_\zeta&u_{\theta,\zeta} \\
                                    u_{s,\theta}+a_1u_\theta-a_3u_\zeta & u_{s,s}+b_1u_\theta-b_3u_\zeta&u_{s,\zeta} \\
                                    u_{\zeta,\theta}+a_2u_\theta+a_3u_s & u_{\zeta,s}+b_2u_\theta+b_3u_s&u_{\zeta,\zeta}
                                  \end{matrix}\right]\left[\begin{matrix}
                   c_1&c_2&0 \\
                   c_3&c_4&0  \\
                   0&0&1
                 \end{matrix}\right]\{d{\bf R}\}\cr
&=\left[\begin{matrix}
          \begin{matrix}
            [c_1(u_{\theta,\theta}-a_1u_s-a_2u_\zeta) \\
            +c_3(u_{\theta,s}-b_1u_s-b_2u_\zeta)]
          \end{matrix} & \begin{matrix}
            [c_2(u_{\theta,\theta}-a_1u_s-a_2u_\zeta) \\
            +c_4(u_{\theta,s}-b_1u_s-b_2u_\zeta)]
          \end{matrix} & u_{\theta,\zeta} \\ \\
          \begin{matrix}
            [c_1(u_{s,\theta}+a_1u_\theta-a_3u_\zeta) \\
            +c_3(u_{s,s}+b_1u_\theta-b_3u_\zeta)]
          \end{matrix} & \begin{matrix}
            [c_2(u_{s,\theta}+a_1u_\theta-a_3u_\zeta) \\
            +c_4(u_{s,s}+b_1u_\theta-b_3u_\zeta)]
          \end{matrix} & u_{s,\zeta}\\ \\
          \begin{matrix}
            [c_1(u_{\zeta,\theta}+a_2u_\theta+a_3u_s) \\
            +c_3(u_{\zeta,s}+b_2u_\theta+b_3u_s)]
          \end{matrix}  & \begin{matrix}
            [c_2(u_{\zeta,\theta}+a_2u_\theta+a_3u_s) \\
            +c_4(u_{\zeta,s}+b_2u_\theta+b_3u_s)]
          \end{matrix} & u_{\zeta,\zeta}
        \end{matrix}\right]\{d{\bf R}\}
\end{align}
Now, comparing with $d{\bf u}=({\boldsymbol\nabla}{\bf u})\,d{\bf R}$, the gradient of the displacement vector ${\bf u}$ can be expressed as
\begin{align}\label{eq:A21}
{\boldsymbol\nabla}{\bf u} &=\left[\begin{matrix}
          \begin{matrix}
            [c_1(u_{\theta,\theta}-a_1u_s-a_2u_\zeta) \\
            +c_3(u_{\theta,s}-b_1u_s-b_2u_\zeta)]
          \end{matrix} & \begin{matrix}
            [c_2(u_{\theta,\theta}-a_1u_s-a_2u_\zeta) \\
            +c_4(u_{\theta,s}-b_1u_s-b_2u_\zeta)]
          \end{matrix} & u_{\theta,\zeta} \\ \\
          \begin{matrix}
            [c_1(u_{s,\theta}+a_1u_\theta-a_3u_\zeta) \\
            +c_3(u_{s,s}+b_1u_\theta-b_3u_\zeta)]
          \end{matrix} & \begin{matrix}
            [c_2(u_{s,\theta}+a_1u_\theta-a_3u_\zeta) \\
            +c_4(u_{s,s}+b_1u_\theta-b_3u_\zeta)]
          \end{matrix} & u_{s,\zeta}\\ \\
          \begin{matrix}
            [c_1(u_{\zeta,\theta}+a_2u_\theta+a_3u_s) \\
            +c_3(u_{\zeta,s}+b_2u_\theta+b_3u_s)]
          \end{matrix}  & \begin{matrix}
            [c_2(u_{\zeta,\theta}+a_2u_\theta+a_3u_s) \\
            +c_4(u_{\zeta,s}+b_2u_\theta+b_3u_s)]
          \end{matrix} & u_{\zeta,\zeta}
        \end{matrix}\right]
\end{align}
In tensor notation the gradient can be written as
\begin{align}\label{eq:A22}
{\boldsymbol \nabla}{\bf u}&=\left(c_1(u_{\theta,\theta}-a_1u_s-a_2u_\zeta)
            +c_3(u_{\theta,s}-b_1u_s-b_2u_\zeta)\right) \hat{\bf e}_\theta\otimes\hat{\bf e}_\theta\cr
            &+\left(c_2(u_{\theta,\theta}-a_1u_s-a_2u_\zeta)
            +c_4(u_{\theta,s}-b_1u_s-b_2u_\zeta)\right)\hat{\bf e}_\theta\otimes\hat{\bf e}_s+u_{\theta,\zeta}\,\hat{\bf e}_\theta\otimes\hat{\bf n}\cr
            &+\left(c_1(u_{s,\theta}+a_1u_\theta-a_3u_\zeta)
            +c_3(u_{s,s}+b_1u_\theta-b_3u_\zeta)\right) \hat{\bf e}_s\otimes\hat{\bf e}_\theta\cr
            &+\left(c_2(u_{s,\theta}+a_1u_\theta-a_3u_\zeta)
            +c_4(u_{s,s}+b_1u_\theta-b_3u_\zeta)\right)\hat{\bf e}_s\otimes\hat{\bf e}_s+u_{s,\zeta}\,\hat{\bf e}_s\otimes\hat{\bf n}\cr
            &+\left(c_1(u_{\zeta,\theta}+a_2u_\theta+a_3u_s)
            +c_3(u_{\zeta,s}+b_2u_\theta+b_3u_s)\right) \hat{\bf n}\otimes\hat{\bf e}_\theta\cr
            &+\left(c_2(u_{\zeta,\theta}+a_2u_\theta+a_3u_s)
            +c_4(u_{\zeta,s}+b_2u_\theta+b_3u_s)\right)\hat{\bf n}\otimes\hat{\bf e}_s+u_{\zeta,\zeta}\,\hat{\bf n}\otimes\hat{\bf n}
\end{align}
The deformation gradient can be obtained as ${\bf F}={\bf I}+{\boldsymbol \nabla}{\bf u}$. The components of the deformation gradient can be written as the column vector, $\tilde{\bf F}$ as defined in Eq.~\eqref{eq:22}. Then the coefficients ${\bf G}_1,\, {\bf G}_2$ and ${\bf G}_3$ can be expressed as
\begin{align}\label{eq:A23}
{\bf G}_1&= \left[\begin{matrix}
                   c_1{\bf A}_{\theta,\theta}+c_3{\bf A}_{\theta,s} & -(b_1c_3+a_1 c_1){\bf A}_{s}& -(b_2c_3+a_2 c_1){\bf A}_{\zeta}\\
                   c_2{\bf A}_{\theta,\theta}+c_4{\bf A}_{\theta,s} & -(b_1c_4+a_1 c_2){\bf A}_{s}& -(b_2c_4+a_2 c_2){\bf A}_{\zeta} \\
                   {\bf A}_{\theta,\zeta}  & {\bf 0} & {\bf 0} \\
                   (b_1c_3+a_1 c_1){\bf A}_{\theta} & c_1{\bf A}_{s,\theta}+c_3{\bf A}_{s,s} & -(b_3c_3+a_3 c_1){\bf A}_{\zeta} \\
                   (b_1c_4+a_1 c_2){\bf A}_{\theta} & c_2{\bf A}_{s,\theta}+c_4{\bf A}_{s,s} & -(b_3c_4+a_3 c_2){\bf A}_{\zeta} \\
                   {\bf 0} &{\bf A}_{s,\zeta}& {\bf 0}\\
                   (b_2c_3+a_2 c_1){\bf A}_{\theta} & (b_3c_3+a_3 c_1){\bf A}_{s}& c_1{\bf A}_{\zeta,\theta}+c_3{\bf A}_{\zeta,s}  \\
                   (b_2c_4+a_2 c_2){\bf A}_{\theta} & (b_3c_4+a_3 c_2){\bf A}_{s}& c_2{\bf A}_{\zeta,\theta}+c_4{\bf A}_{\zeta,s}\\
                   {\bf 0} &{\bf 0}&{\bf A}_{\zeta,\zeta}
                 \end{matrix}\right]\cr\cr
{\bf G}_2 &=\left[\begin{matrix}
                       c_1{\bf A}_\theta & {\bf 0} & {\bf 0} \\
                       c_2{\bf A}_\theta & {\bf 0} & {\bf 0} \\
                       {\bf 0} & {\bf 0} & {\bf 0} \\
                       {\bf 0} & c_1{\bf A}_{s} & {\bf 0} \\
                       {\bf 0} & c_2{\bf A}_{s} & {\bf 0} \\
                       {\bf 0} & {\bf 0} & {\bf 0} \\
                       {\bf 0} & {\bf 0} & c_1{\bf A}_{\zeta} \\
                       {\bf 0} & {\bf 0} & c_2{\bf A}_{\zeta} \\
                       {\bf 0} & {\bf 0} & {\bf 0}
                     \end{matrix}\right],\quad
{\bf G}_3 =\left[\begin{matrix}
                       c_3{\bf A}_\theta & {\bf 0} & {\bf 0} \\
                       c_4{\bf A}_\theta & {\bf 0} & {\bf 0} \\
                       {\bf 0} & {\bf 0} & {\bf 0} \\
                       {\bf 0} & c_3{\bf A}_{s} & {\bf 0} \\
                       {\bf 0} & c_4{\bf A}_{s} & {\bf 0} \\
                       {\bf 0} & {\bf 0} & {\bf 0} \\
                       {\bf 0} & {\bf 0} & c_3{\bf A}_{\zeta} \\
                       {\bf 0} & {\bf 0} & c_4{\bf A}_{\zeta} \\
                       {\bf 0} & {\bf 0} & {\bf 0}
                     \end{matrix}\right],\quad
\tilde{\bf I} =\left\{\begin{matrix}
                      1 \\
                      0 \\
                      0 \\
                      0 \\
                      1 \\
                      0 \\
                      0 \\
                      0 \\
                      1
                    \end{matrix}\right\}
\end{align}
where $(\ )_{,\zeta}$ represents the partial derivative with respect to $\zeta$ and so on.
\subsubsection{Specialization for tubular curved shell with circular cross-section}
The expression for the gradient given Eq.~\eqref{eq:A22} can be specialized for curved tubular shell with circular cross-section. For this case, the radius $R$ would be a function of the coordinate $s$ only. In this case the $\alpha$'s in Eq.~\eqref{eq:A4} becomes:
\begin{align}\label{eq:A24}
\alpha_1 &= - {\partial R \over \partial s} ,\quad \alpha_2 = R\,\cos\theta,\quad\alpha_3=- R\,\sin\theta\cr
\alpha_4 &= R\sqrt{\alpha_1^2+\xi^2},\quad \alpha_5=R
\end{align}
and $a$'s defined in Eq.~\eqref{eq:A8} can be written as follows:
\begin{align}\label{eq:A25}
a_{1} & = -{R\over \alpha_4} {\partial R \over \partial s},\quad
a_{2}  = {-\xi\over \sqrt{\alpha_1^2+\xi^2}},\quad
a_{3}  = {\alpha_1 R\hat{\kappa}_2\over(\alpha_1^2+\xi^2)}\cr
b_{1} & ={R\kappa_3\alpha_1\over\alpha_4}+{\xi R\hat{\kappa}_2\over \alpha_4} \cr
b_{2} & =  -{R\xi\kappa_3\over\alpha_4}-{R\hat{\kappa}_2\over \alpha_4}{\partial R\over\partial s}\cr
b_{3} & = \hat{\kappa}_1+{R\over\alpha_4^2}{\partial R\over\partial s}{\partial R\xi\over\partial s}
\end{align}
and the $c$'s defined in Eq.~\eqref{eq:A16} becomes:
\begin{equation}\label{eq:A26}
\begin{aligned}
& c_1&=&{1\over g}\left({\alpha_4\over \alpha_5}-b_3\zeta\right), & &c_3&=&{1\over g}a_3\zeta &\cr
& c_2&=&\,-{1\over g}\left(\kappa_3R-b_2\zeta\right), & &c_4&=& \,{1\over g}\left(R-a_2\zeta\right)&
\end{aligned}
\end{equation}
with
\begin{align}\label{eq:A27}
g = \left(R-a_2\zeta\right)\left({\alpha_4\over \alpha_5}-b_3\zeta\right)
+a_3\zeta\left(\kappa_3R-b_2\zeta\right)
\end{align}
\subsubsection{Specialization for curved tubular shell with constant radius}
\begin{figure}[!htbp]
\begin{adjustbox}{center}\hskip20pt
\begin{tabular}{ccc}
\subfloat[ ]{\includegraphics[trim=2cm 0.5cm 5.5cm 1cm,clip, scale=0.3]{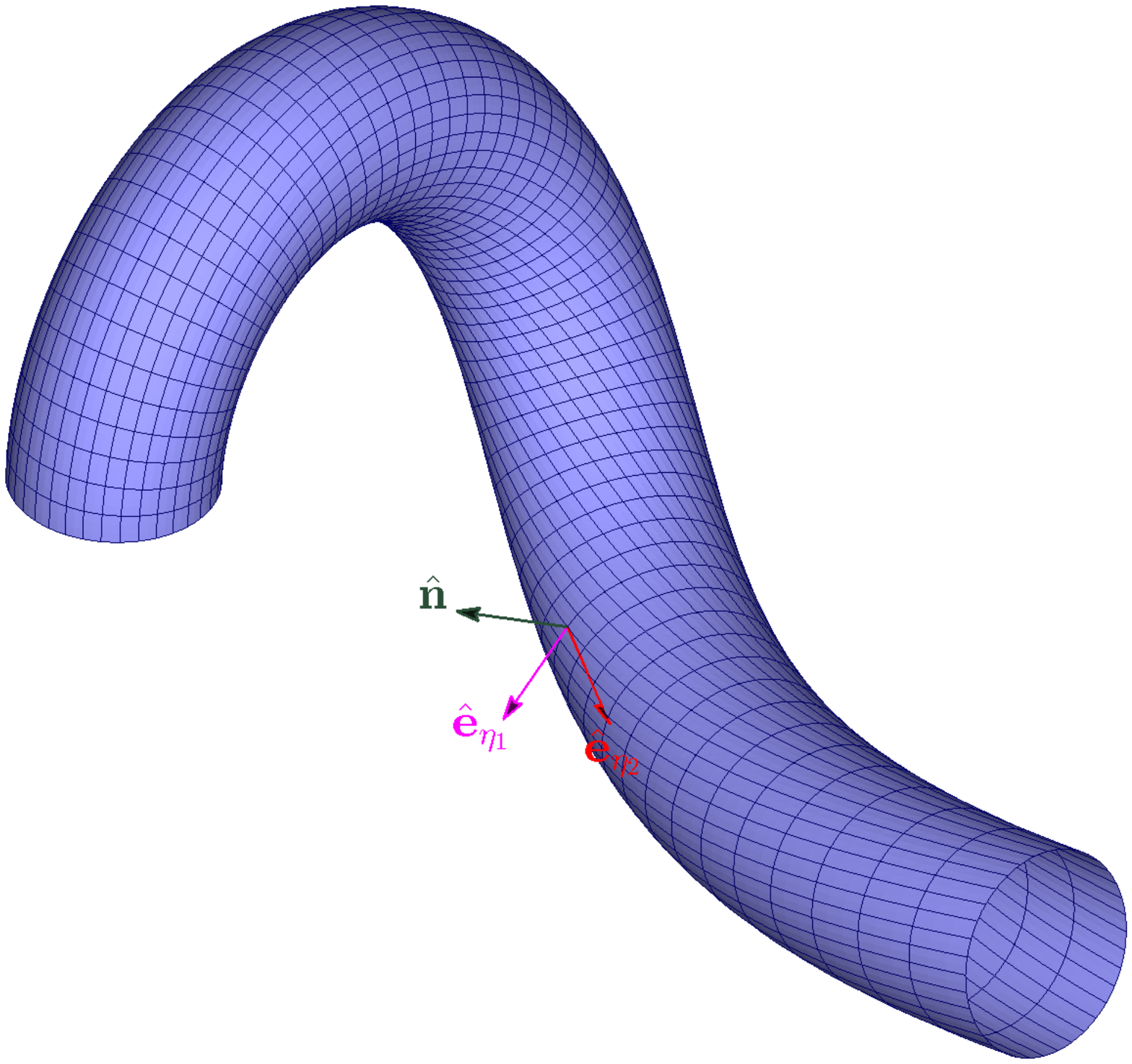}}
& \subfloat[]{\includegraphics[trim=3cm 0.5cm 5.5cm 1cm,clip,scale=0.3]{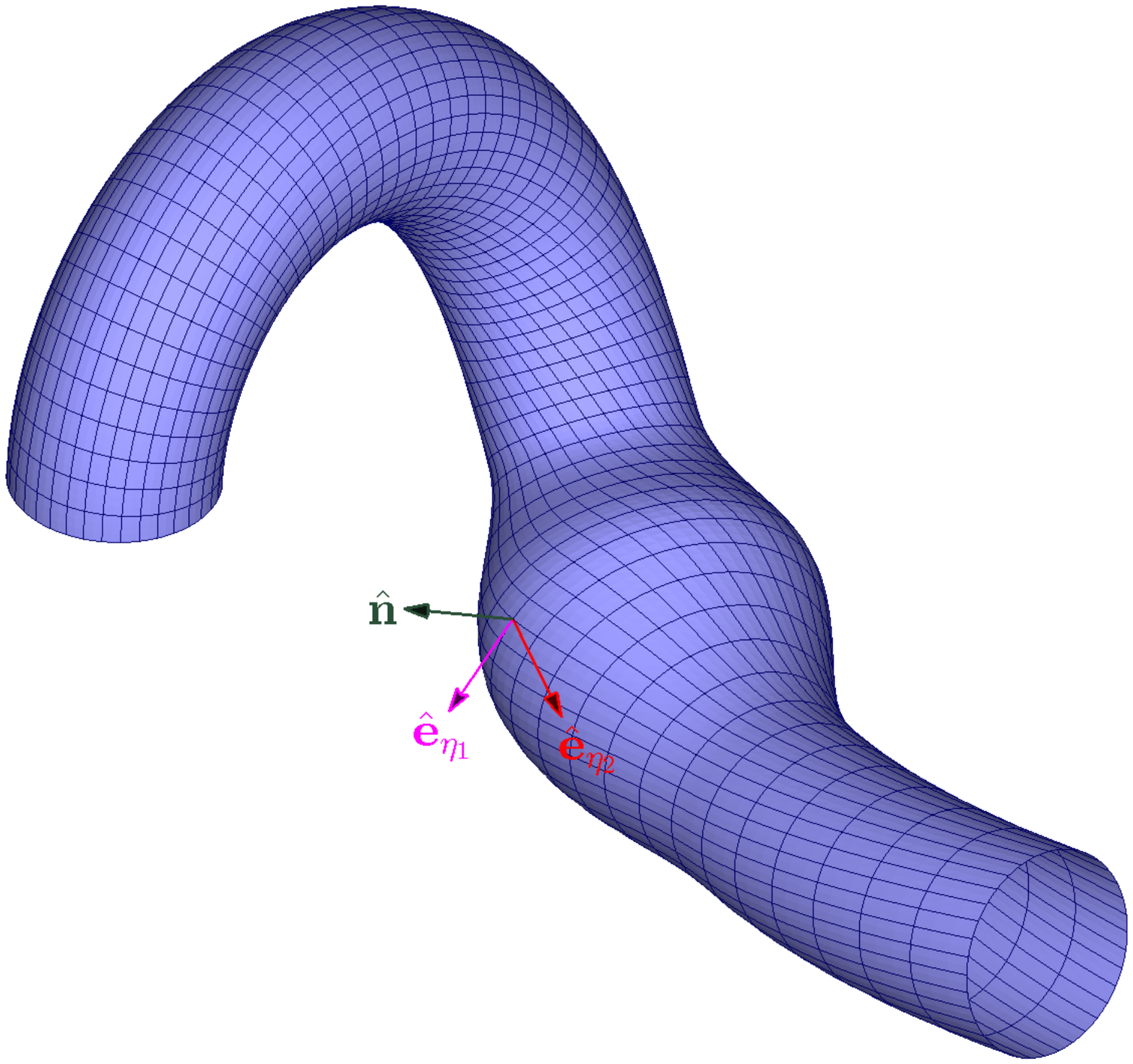}}
& \subfloat[]{\includegraphics[trim=2.5cm 0.5cm 5.5cm 1cm,clip,scale=0.3]{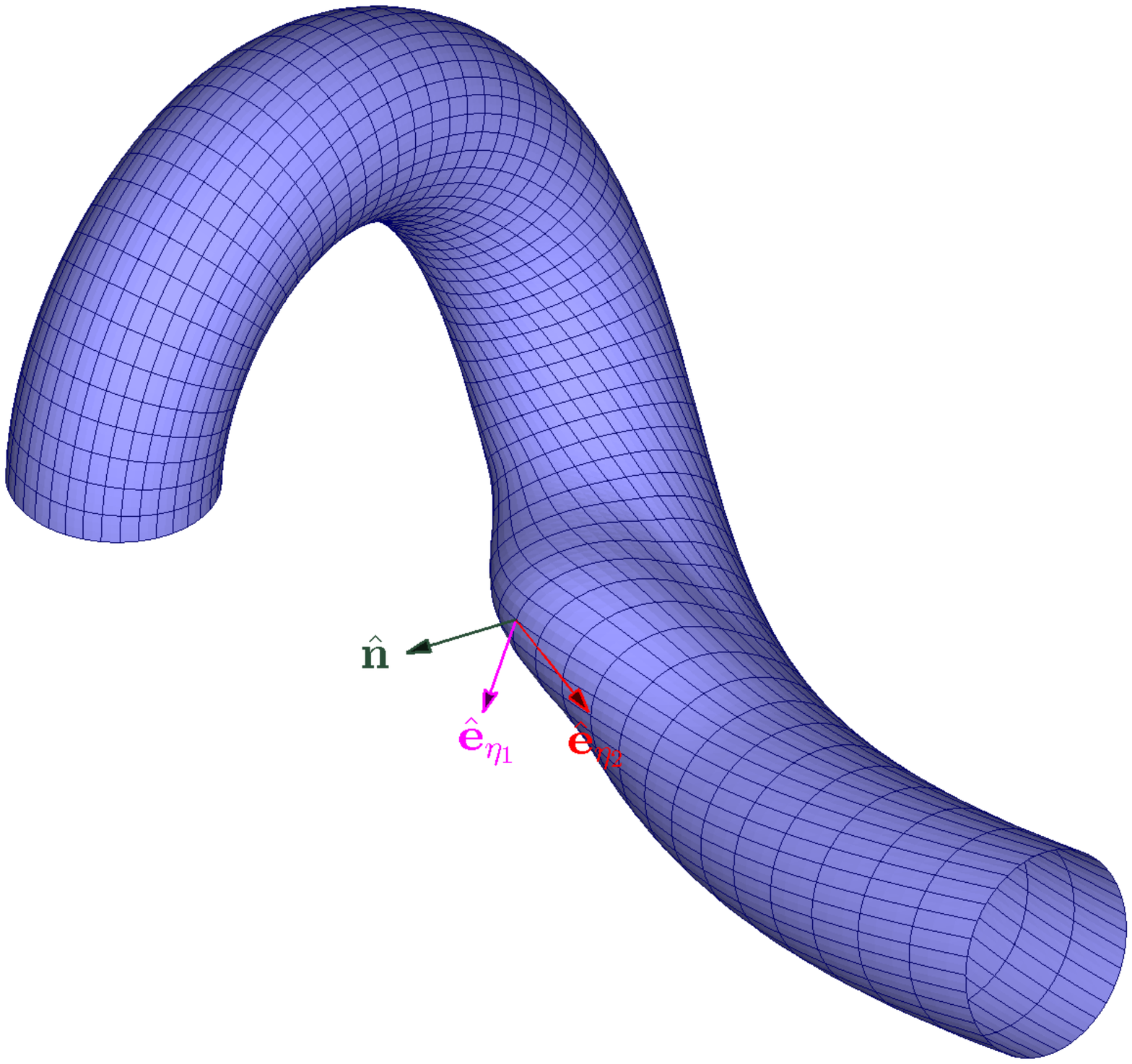}}
\end{tabular}
\end{adjustbox}
\caption{Curved tubular shells with constant or varying radial distance and Cartan's moving frame on their surfaces.}\label{fig:A2}
\end{figure}
The expression for the gradient given Eq.~\eqref{eq:A20} can be specialized for the curved tubular shell with a circular cross-section. For this, the radius $R$ will be a function of the coordinate $s$ only. In this case, the $\alpha$'s in Eq.~\eqref{eq:A4} become:
\begin{align}\label{eq:A28}
\alpha_1 &= 0 ,\quad \alpha_2 = R\,\cos\theta,\quad\alpha_3=- R\,\sin\theta\cr
\alpha_4 &= R\xi,\quad \alpha_5=R
\end{align}
and $a$'s defined in Eq.~\eqref{eq:A8} can be given as follows:
\begin{align}\label{eq:A29}
a_{1} & = 0,\quad
a_{2}  = -1,\quad
a_{3}  = 0\cr
b_{1} & =\hat{\kappa}_2,\quad
b_{2}  =  -\kappa_3,\quad
b_{3}  = \hat{\kappa}_1
\end{align}
and the $c$'s defined in Eq.~\eqref{eq:A16} become:
\begin{equation}\label{eq:A30}
\begin{aligned}
& c_1&=&{1\over \left(R+\zeta\right)}, & &c_3&=&0 &\cr
& c_2&=&\,-{\kappa_3\over \left(\xi-\hat{\kappa}_1\zeta\right)}, & &c_4&=& \,{1\over \left(\xi-\hat{\kappa}_1\zeta\right)}&
\end{aligned}
\end{equation}
and
\begin{align}\label{eq:A31}
g = \left(R+\zeta\right)\left(\xi-\hat{\kappa}_1\zeta\right)
\end{align}
\subsubsection{Specialization for surface of revolution}
The above gradient can be specialized for shells with surface of revolution as well. For the surface of revolution the radius of cross-section again would be function of $s$, that is, $R=R(s)$ along with the components of the Cartan matrix defined in Eq.~\eqref{eq:A1} would also be equal to zero because in that case the reference curve is a straight line.Hence, in this case, we have
\begin{equation}\label{eq:A32}
\kappa_1=\kappa_2=\kappa_3=\hat{\kappa}_1=\hat{\kappa}_2=0,\quad \xi=1,
\end{equation}
Also, the $\alpha$'s in Eq.~\eqref{eq:A4} become:
\begin{align}\label{eq:A33}
\alpha_1 &= - {\partial R \over \partial s} ,\quad \alpha_2 = R\,\cos\theta,\quad\alpha_3=- R\,\sin\theta\cr
\alpha_4 &= R\sqrt{\alpha_1^2+1},\quad \alpha_5=R
\end{align}
and $a$'s defined in Eq.~\eqref{eq:A8} can be given as follows:
\begin{align}\label{eq:A34}
a_{1} & = -{1\over \sqrt{\alpha_1^2+1}} {\partial R \over \partial s},\quad
a_{2}  = {-1\over \sqrt{\alpha_1^2+1}},\quad
a_{3}  = 0\cr
b_{1} & =0 ,\quad
b_{2}  =  0,\quad
b_{3}  = {R^2\over\alpha_4^2}\left({\partial^2 R\over\partial s^2}\right)
\end{align}
and the $c$'s defined in Eq.~\eqref{eq:A16} become:
\begin{equation}\label{eq:A35}
\begin{aligned}
& c_1&=&{1\over \left(R-a_2\zeta\right)}, & &c_3&=&\,0 &\cr
& c_2&=&\,0, & &c_4&=& \,{1\over \left(\sqrt{\alpha_1^2+1}-b_3\zeta\right)}&
\end{aligned}
\end{equation}
and
\begin{align}\label{eq:A36}
g = \left(R-a_2\zeta\right)\left(\sqrt{\alpha_1^2+1}-b_3\zeta\right)
\end{align}

\begin{figure}[!htbp]
\begin{adjustbox}{center}
\begin{tabular}{ccc}
\subfloat[ ]{\includegraphics[trim=6cm 0.5cm 5cm 1.8cm,clip, scale=0.3]{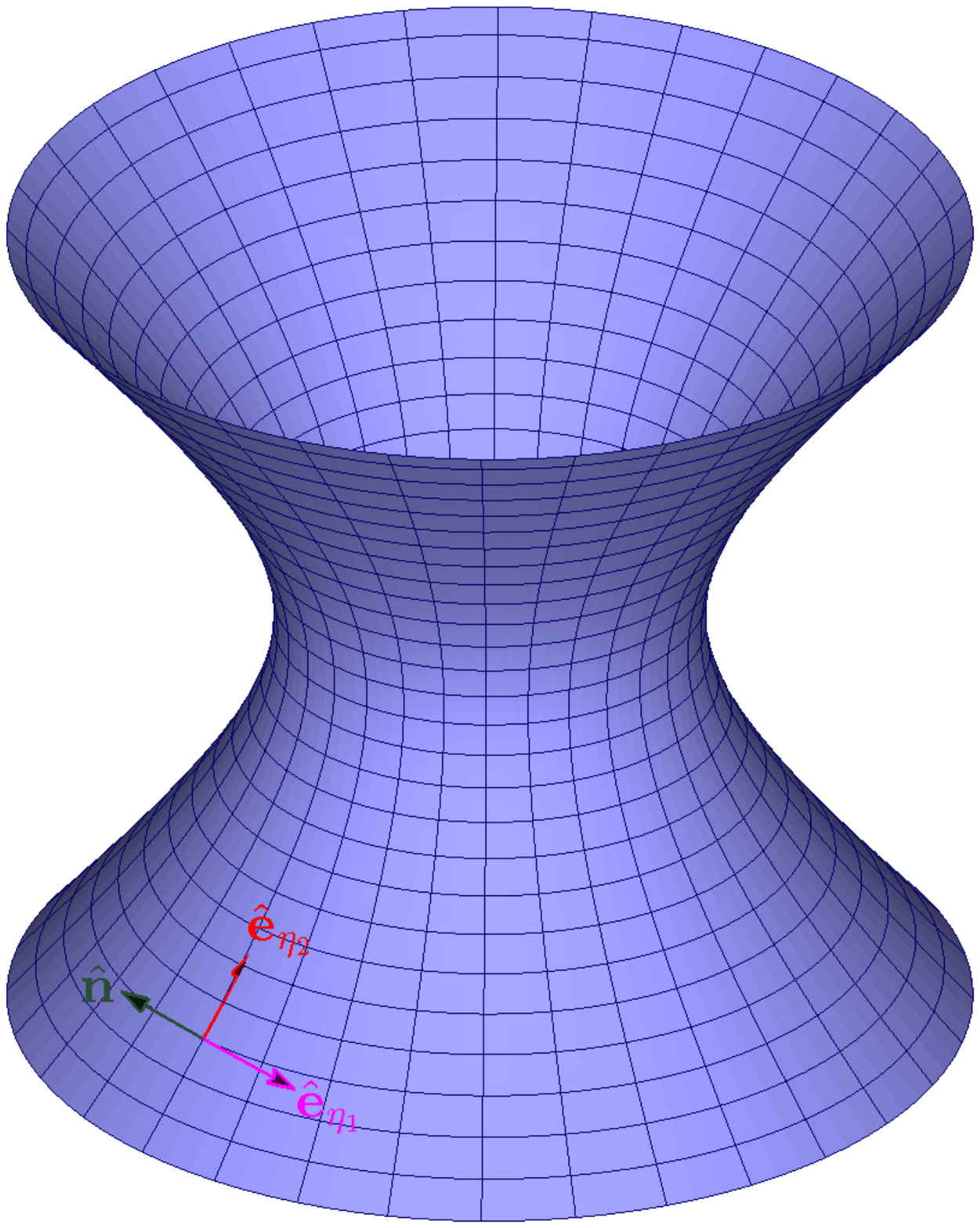}}
& \subfloat[]{\includegraphics[trim=7cm 2.5cm 5.5cm 2cm,clip,scale=0.36]{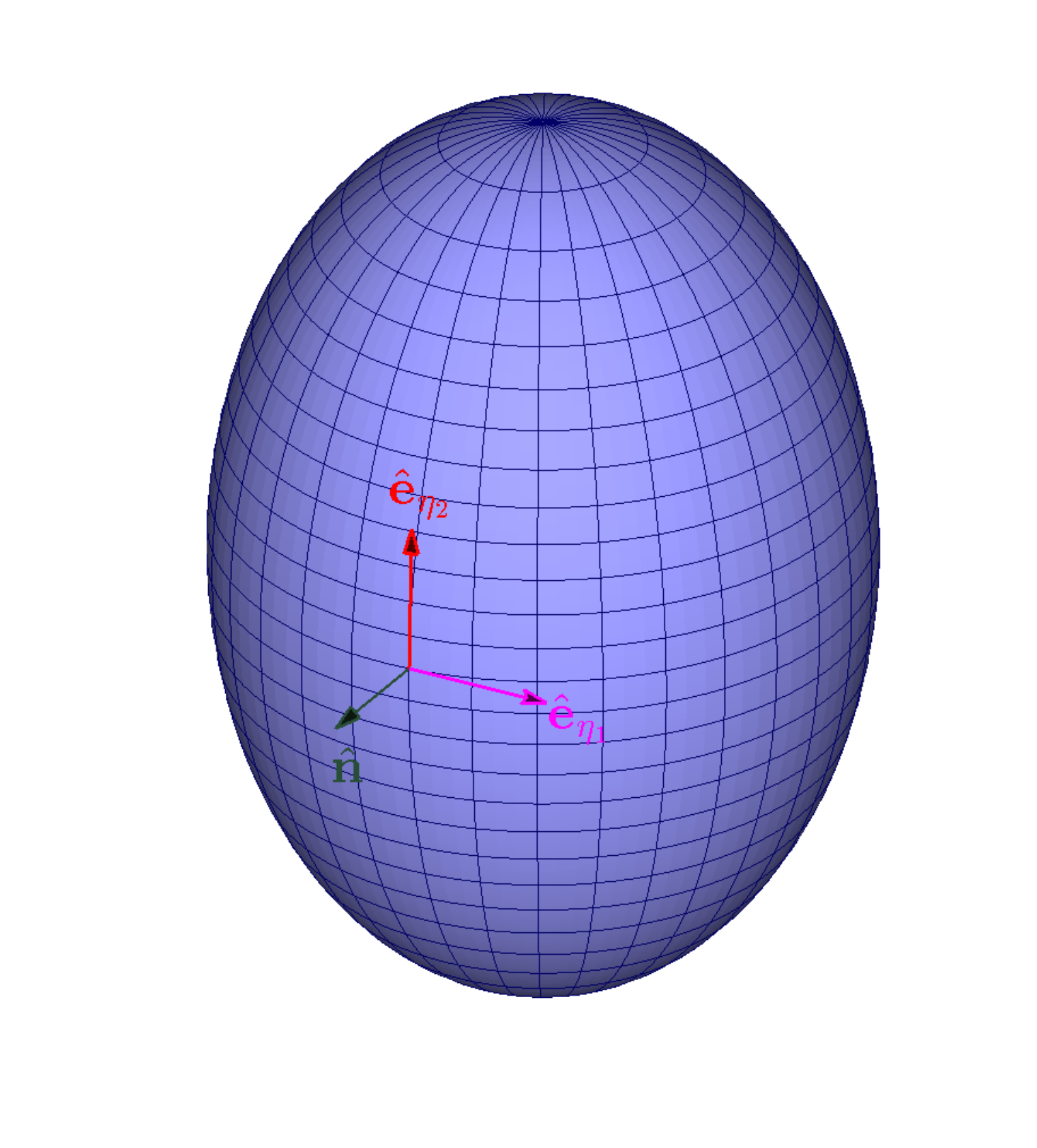}}&
\subfloat[]{\includegraphics[trim=6cm 2cm 5cm 2cm,clip,scale=0.36]{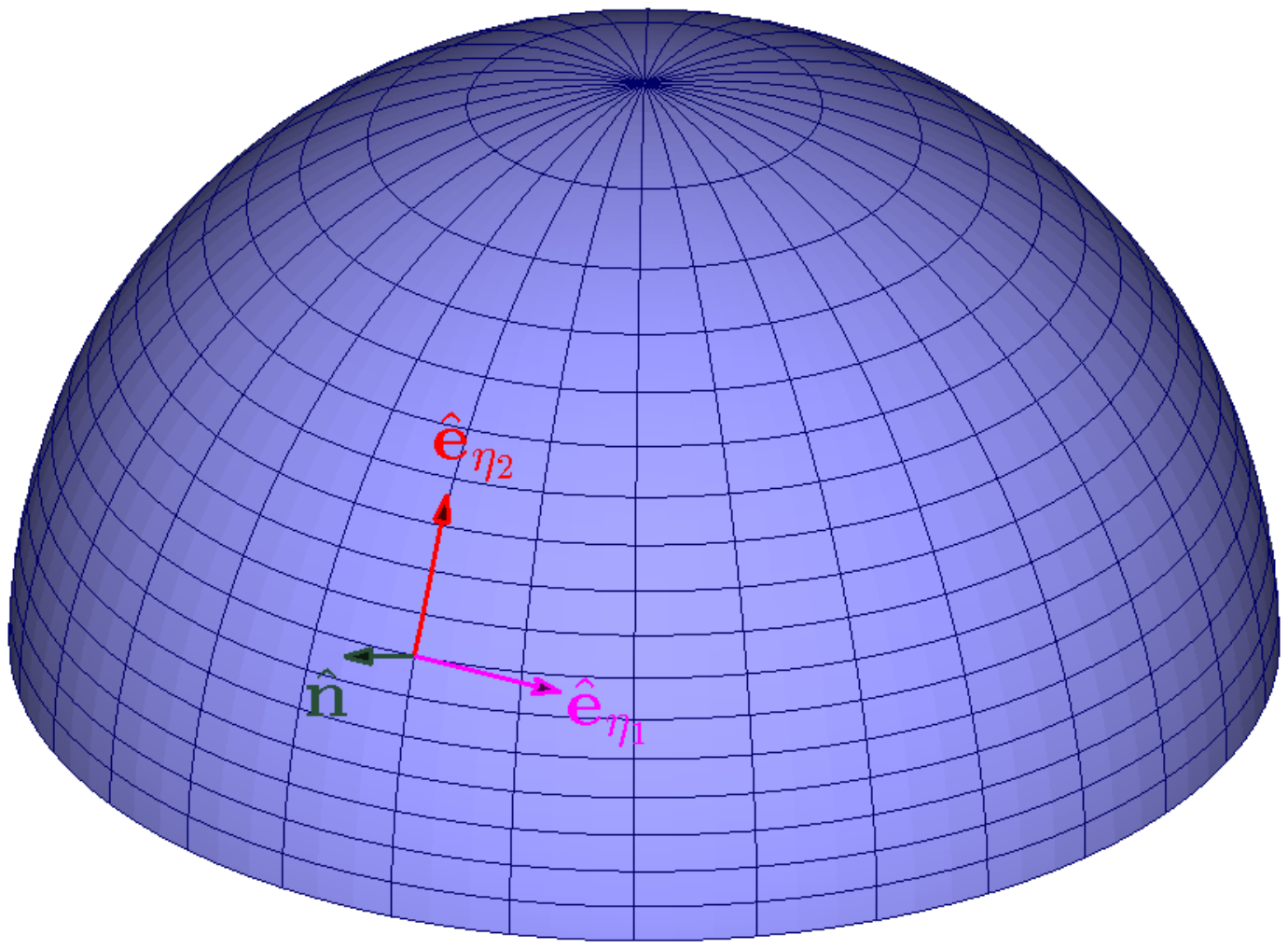}}\\
\subfloat[]{\includegraphics[trim=6cm 4cm 5.5cm 1cm,clip,scale=0.34]{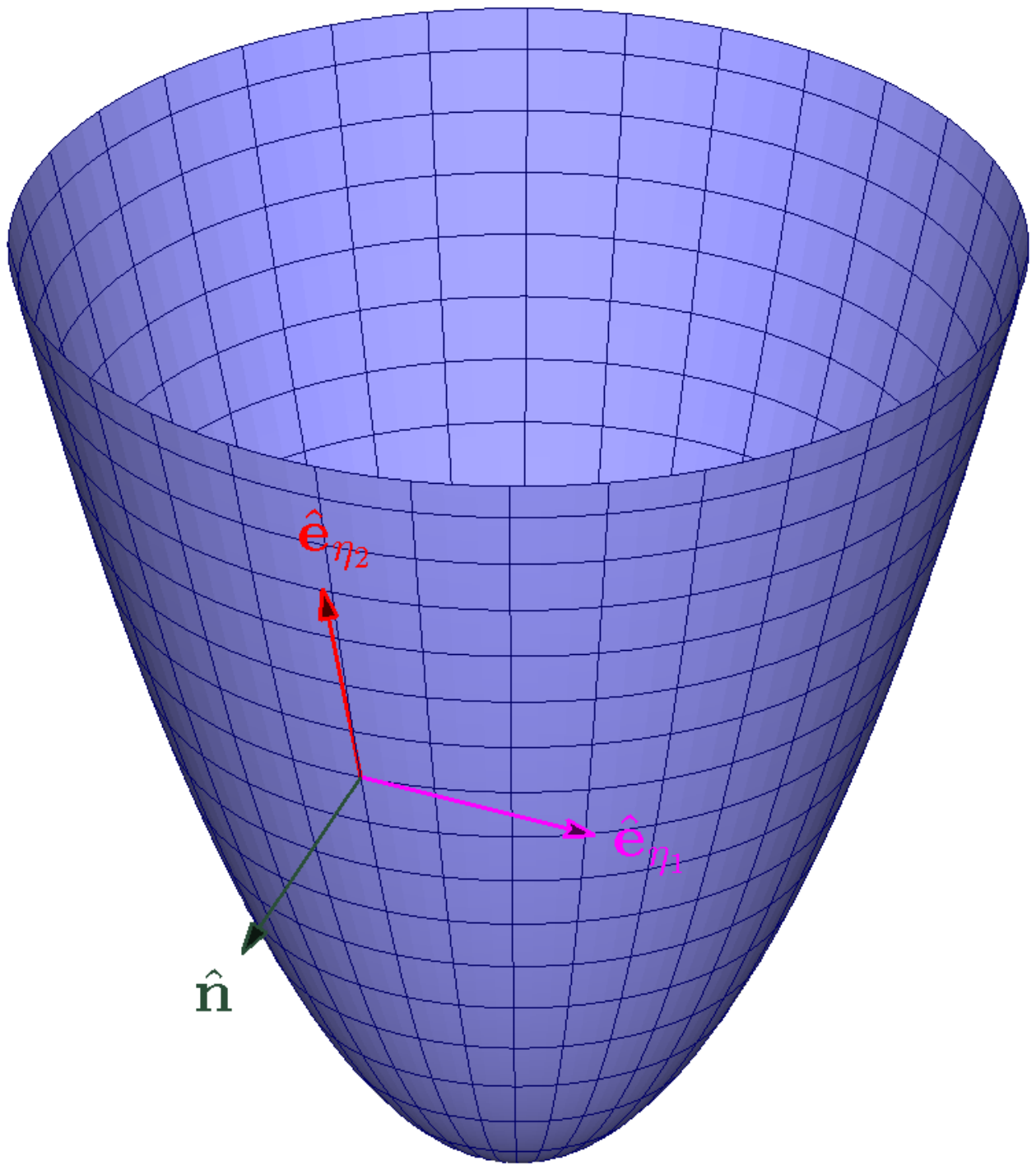}}
&\subfloat[]{\includegraphics[trim=6cm 2.2cm 5.5cm 1cm,clip,scale=0.34]{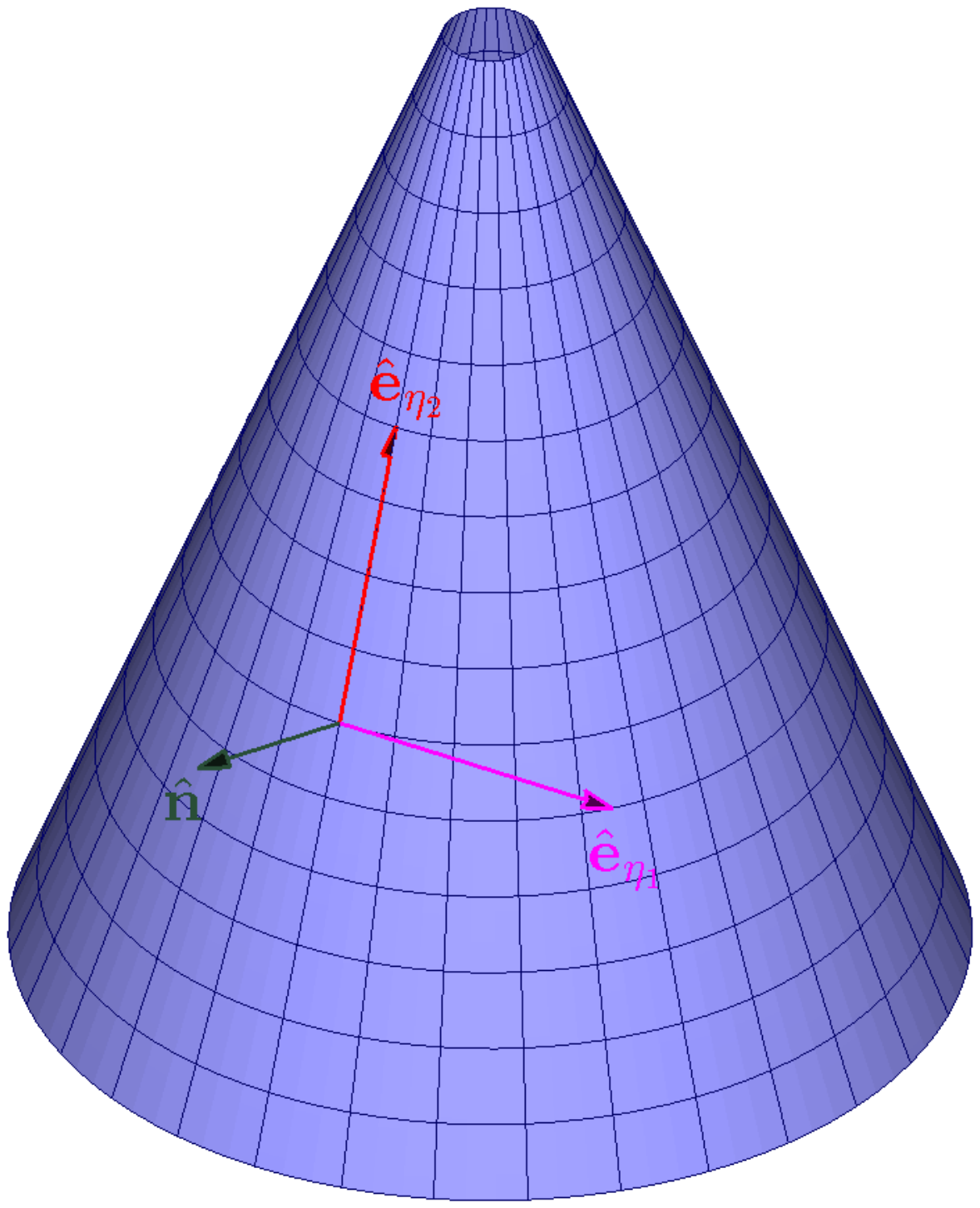}}
& \subfloat[]{\includegraphics[trim=6cm 2cm 5cm 1cm,clip,scale=0.29]{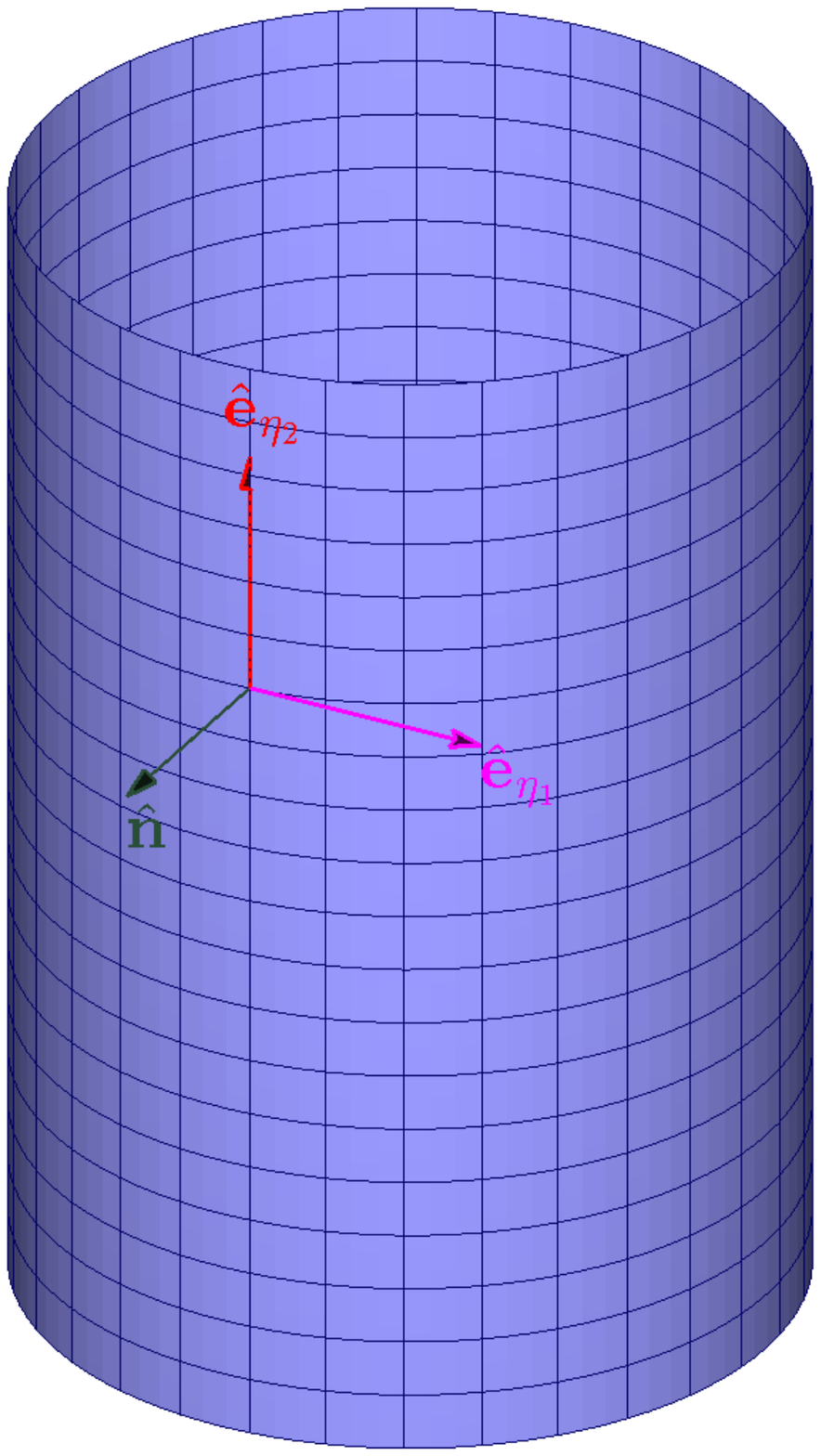}}
\end{tabular}
\end{adjustbox}
\caption{Surface of revolution}\label{fig:A3}
\end{figure}
\subsubsection{Specialization for cylindrical shell}
The cylindrical surface is a very special kind of surface of revolution  where $R$ is constant. In this case the $\alpha$'s in Eq.~\eqref{eq:A4} become:
\begin{align}\label{eq:A37}
\alpha_1 &= 0 ,\quad \alpha_2 = R\,\cos\theta,\quad\alpha_3=- R\,\sin\theta\cr
\alpha_4 &= R,\quad \alpha_5=R
\end{align}
and $a$'s defined in Eq.~\eqref{eq:A8} can be given as follows:
\begin{align}\label{eq:A38}
a_{1} & = 0,\quad
a_{2}  = -1,\quad
a_{3}  = 0\cr
b_{1} & =0 ,\quad
b_{2}  =  0,\quad
b_{3}  = 0
\end{align}
and the $c$'s defined in Eq.~\eqref{eq:A16} become:
\begin{equation}\label{eq:A39}
\begin{aligned}
& c_1&=&{1\over \left(R+\zeta\right)}, & &c_3&=&\,0 &\cr
& c_2&=&\,0, & &c_4&=& \,1&
\end{aligned}
\end{equation}
and
\begin{align}\label{eq:A40}
g = \left(R+\zeta\right)
\end{align}

\subsection{Spherical shell}
Here we also derive the gradient of the spherical shell following similar procedure for the sake of completeness. Let us consider, a spherical shell with radius $R$ and the surface coordinates $\eta_1=\phi$ and $\eta_2=\theta$, where $\theta$ is the azimuthal angle in the $xy$-plane, measured from the $x$-axis and $\phi$ is the polar angle (also known as the zenith angle or colatitude) as shown in the Fig.~\ref{fig:A1}. Also, $\zeta$ is the coordinate along the normal direction $\hat{\bf n}$, which is essentially the radial direction in the case of the spherical shell.
\begin{figure}[htbp]
\centering{\includegraphics[trim=3.7cm 4.6cm 3.7cm 1.5cm,clip,scale=0.63]{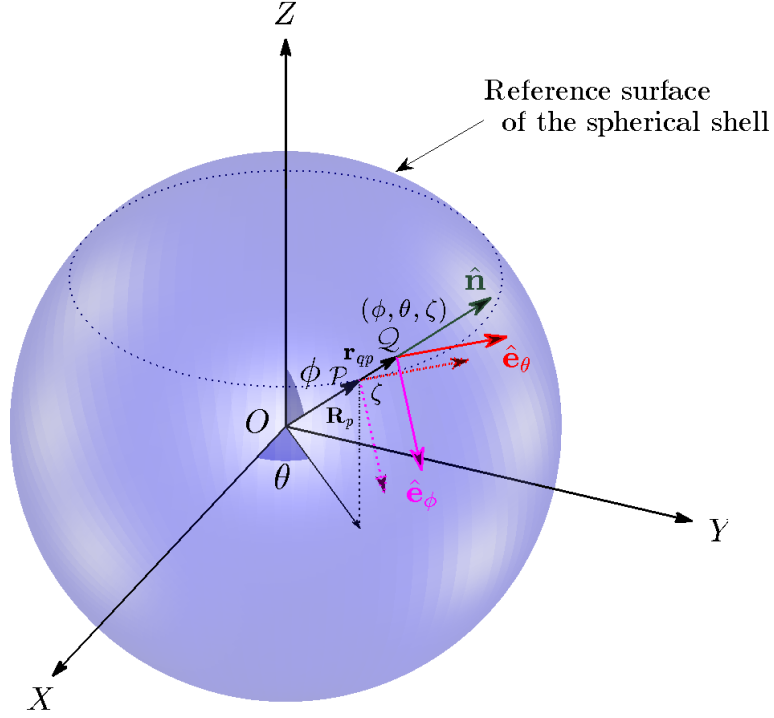}}
\caption{Arbitrary surface in Spherical coordinate system.}
\label{fig:A4}
\end{figure}
\begin{align}\label{eq:A41}
\hat{\bf e}_{\eta_1}&=\hat{\bf e}_{\phi} = \cos\theta\,\cos\phi\,\hat{\bf e}_x+\sin\theta\,\cos\phi\,\hat{\bf e}_y-\sin\phi\, \hat{\bf e}_z\cr
\hat{\bf e}_{\eta_2}&=\hat{\bf e}_{\theta} = -\sin\theta\,\hat{\bf e}_x +\cos\theta\,\hat{\bf e}_y\cr
\hat{\bf e}_{\zeta} &=\hat{\bf n}=\cos\theta\,\sin\phi\,\hat{\bf e}_x+\sin\theta\,\sin\phi\,\hat{\bf e}_y+\cos\phi\, \hat{\bf e}_z
\end{align}
Next, we obtain the exterior derivative of the above considered orthonormal basis vectors can be given as follows:
\begin{equation}\label{eq:A42}
\begin{aligned}
&d\hat{\bf e}_\phi &=&&&&\cos\phi\,d\theta\,\,\hat{\bf e}_\theta&& -d\phi\,\hat{\bf n}&\cr
&d\hat{\bf e}_\theta &=&&-\cos\phi\,d\theta\,\hat{\bf e}_\phi &&&&-\sin\phi\,d\theta\,\hat{\bf n}&\cr
&d\hat{\bf n} &=&&d\phi\,\hat{\bf e}_\phi && +\sin\phi\,d\theta\,\,\hat{\bf e}_\theta&&&
\end{aligned}
\end{equation}
The position vector of any arbitrary point $\mathcal{Q}$ can be given as follows:
\begin{eqnarray}\label{eq:A43}
{\bf R}=(R+\zeta)\,\hat{\bf n}
\end{eqnarray}
Further, the length element, which is a vector valued one form, is:
\begin{align}\label{eq:A44}
d{\bf R}&=(R+\zeta)d\phi\,\hat{\bf e}_\phi+(R+\zeta)\sin\phi\,d\theta\,\,\hat{\bf e}_\theta +d\zeta\,\hat{\bf n}
\end{align}
Let us write the length element $d{\bf R}$ as the following column vector\footnote{Here, in the column vector, the first element is component along $\hat{\bf e}_\phi$ and the second and third elements are the components along $\hat{\bf e}_\theta$ and $\hat{\bf n}$, respectively. }:
\begin{align}\label{eq:A45}
\{d{\bf R}\}&=\left\{\begin{matrix}
                    (R+\zeta)d\phi\\
                   (R+\zeta)\sin\phi\,d\theta\\
                   d\zeta
                 \end{matrix}\right\}=\left[\begin{matrix}
                   (R+\zeta)&0&0 \\
                   0&(R+\zeta)\sin\phi&0  \\
                   0&0&1
                 \end{matrix}\right]\left\{\begin{matrix}
                   d\phi \\
                   d\theta\\
                   d\zeta
                 \end{matrix}\right\}
\end{align}
Further, the volume element $dV$ in this case can be obtained as follows:
\begin{align}\label{eq:A46}
dV&=g\, d\phi\,d\theta\, d\zeta,\quad \text{where}\quad g = (R+\zeta)^2\sin\phi
\end{align}
Now, let us consider the displacement vector as ${\bf u}=u_\phi\hat{\bf e}_\phi+u_\theta \hat{\bf e}_\theta+u_\zeta\hat{\bf n}$, then we have the differential, $d{\bf u}:=({\boldsymbol\nabla}{\bf u})d{\bf R}$, where ${\boldsymbol\nabla}{\bf u}$ is the gradient of displacement vector. Now, $d{\bf u}$ can be given as following:
\begin{align}\label{eq:A47}
d{\bf u} &= d{\bf u}_{rel}+u_\phi\,d\hat{\bf e}_\phi +u_\theta\,d\hat{\bf e}_\theta+u_\zeta\,d\hat{\bf n}\cr
&= d{\bf u}_{rel} +(u_\zeta\,d\phi-u_\theta\,\cos\phi\,d\theta)\,\hat{\bf e}_\phi \cr&+(u_\phi\,\cos\phi+u_\zeta\,\sin\phi)\,d\theta\,\hat{\bf e}_\theta -(u_\theta\,\sin\phi\,d\theta+u_\phi\,d\phi)\,\hat{\bf n}
\end{align}
which can be again written as a column vector as follows:
\begin{align}\label{eq:A48}
 \{d{\bf u}\} &= \left(\left[\begin{matrix}
                                    u_{\phi,\phi}&u_{\phi,\theta}&u_{\phi,\zeta}\\
                                    u_{\theta,\phi}&u_{\theta,\theta}&u_{\theta,\zeta}\\
                                    u_{\zeta,\phi}&u_{\zeta,\theta}&u_{\zeta,\zeta}
                                  \end{matrix}\right]+\left[\begin{matrix}
                                  u_\zeta &-u_\theta\,\cos\phi &  0\\
                                  0&u_\phi\,\cos\phi+u_\zeta\,\sin\phi&0 \\
                                   -u_\phi &-u_\theta\,\sin\phi & 0
                                  \end{matrix}\right]\right)\left\{\begin{matrix}
                   d\phi \\
                   d\theta\\
                   d\zeta
                 \end{matrix}\right\}
\end{align}
and
\begin{align}\label{eq:A49}
 \{d{\bf u}\} &= \left[\begin{matrix}
                                    u_{\phi,\phi}+u_\zeta&u_{\phi,\theta}-u_\theta\,\cos\phi&u_{\phi,\zeta}\\
                                    u_{\theta,\phi}&u_{\theta,\theta}+\cos\phi\, u_\phi+u_\zeta\,\sin\phi&u_{\theta,\zeta}\\
                                    u_{\zeta,\phi}-u_\phi&u_{\zeta,\theta}-u_\theta\,\sin\phi&u_{\zeta,\zeta}
                                  \end{matrix}\right]\left\{\begin{matrix}
                   d\phi\\
                   d\theta \\
                   d\zeta
                  \end{matrix}\right\}\cr\cr
                  &=\left[\begin{matrix}
                                    u_{\phi,\phi}+u_\zeta&u_{\phi,\theta}-u_\theta\,\cos\phi&u_{\phi,\zeta}\\
                                    u_{\theta,\phi}&u_{\theta,\theta}+\cos\phi\, u_\phi+u_\zeta\,\sin\phi&u_{\theta,\zeta}\\
                                    u_{\zeta,\phi}-u_\phi&u_{\zeta,\theta}-u_\theta\,\sin\phi&u_{\zeta,\zeta}
                                  \end{matrix}\right]\left[\begin{matrix}
                   {1\over(R+\zeta)}&0&0 \\
                   0&{1\over(R+\zeta)\sin\phi}&0  \\
                   0&0&1
                 \end{matrix}\right]\{d{\bf R}\}\cr\cr
                  &=\left[\begin{matrix}
                                    {1\over(R+\zeta)}(u_{\phi,\phi}+u_\zeta)&{1\over(R+\zeta)\sin\phi}(u_{\phi,\theta}-u_\theta\,\cos\phi)&u_{\phi,\zeta}\\[8pt]
                                    {1\over(R+\zeta)}u_{\theta,\phi}&{1\over(R+\zeta)\sin\phi}(u_{\theta,\theta}+\cos\phi\, u_\phi+u_\zeta\,\sin\phi)&u_{\theta,\zeta}\\[8pt]
                                    {1\over(R+\zeta)}(u_{\zeta,\phi}-u_\phi)&{1\over(R+\zeta)\sin\phi}(u_{\zeta,\theta}-u_\theta\,\sin\phi)&u_{\zeta,\zeta}
                                  \end{matrix}\right]\{d{\bf R}\}
\end{align}
Now, comparing with $d{\bf u}=({\boldsymbol\nabla}{\bf u})\,d{\bf R}$, the gradient of the displacement vector ${\bf u}$ can be given as follows:
\begin{align}\label{eq:A50}
{\boldsymbol\nabla}{\bf u} &=\left[\begin{matrix}
                                    {1\over(R+\zeta)}(u_{\phi,\phi}+u_\zeta)&{1\over(R+\zeta)\sin\phi}(u_{\phi,\theta}-u_\theta\,\cos\phi)&u_{\phi,\zeta}\\[8pt]
                                    {1\over(R+\zeta)}u_{\theta,\phi}&{1\over(R+\zeta)\sin\phi}(u_{\theta,\theta}+\cos\phi\, u_\phi+u_\zeta\,\sin\phi)&u_{\theta,\zeta}\\[8pt]
                                    {1\over(R+\zeta)}(u_{\zeta,\phi}-u_\phi)&{1\over(R+\zeta)\sin\phi}(u_{\zeta,\theta}-u_\theta\,\sin\phi)&u_{\zeta,\zeta}
                                  \end{matrix}\right]
\end{align}
In tensor notation the gradient can be written as follows:
\begin{align}\label{eq:A51}
{\boldsymbol \nabla}{\bf u}&= {1\over(R+\zeta)}(u_{\phi,\phi}+u_\zeta)\,\hat{\bf e}_\phi\otimes\hat{\bf e}_\phi+{1\over(R+\zeta)\sin\phi}(u_{\phi,\theta}-u_\theta\,\cos\phi)\hat{\bf e}_\phi\otimes\hat{\bf e}_\theta+u_{\phi,\zeta}\,\hat{\bf e}_\phi\otimes\hat{\bf n}\cr
            &+{1\over(R+\zeta)}u_{\theta,\phi} \hat{\bf e}_\theta\otimes\hat{\bf e}_\phi+{1\over(R+\zeta)\sin\phi}(u_{\theta,\theta}+\cos\phi\, u_\phi+u_\zeta\,\sin\phi)\hat{\bf e}_\theta\otimes\hat{\bf e}_\theta+u_{\theta,\zeta}\,\hat{\bf e}_\theta\otimes\hat{\bf n}\cr
            &+{1\over(R+\zeta)}(u_{\zeta,\phi}-u_\phi) \hat{\bf n}\otimes\hat{\bf e}_\phi+{1\over(R+\zeta)\sin\phi}(u_{\zeta,\theta}-u_\theta\,\sin\phi)\hat{\bf n}\otimes\hat{\bf e}_\theta+u_{\zeta,\zeta}\,\hat{\bf n}\otimes\hat{\bf n}
\end{align}
The coefficients ${\bf G}_1,\, {\bf G}_2$ and ${\bf G}_3$ used in Eq.~\eqref{eq:23} for the gradient can be given as
\begin{align}\label{eq:A52}
{\bf G}_1&= \left[\begin{matrix}
             {1\over(R+\zeta)}{\bf A}_{\phi,\phi}&{\bf 0}&{1\over(R+\zeta)}{\bf A}_\zeta\\[5pt]
             {1\over(R+\zeta)\sin\phi}{\bf A}_{\phi,\theta}&-{\cos\phi\over(R+\zeta)\sin\phi}{\bf A}_\theta&{\bf 0}\\[5pt]
             {\bf A}_{\phi,\zeta}&{\bf 0}&{\bf 0}\\[5pt]
             {\bf 0}&{1\over(R+\zeta)}{\bf A}_{\theta,\phi}&{\bf 0}\\[5pt]
             {\cos\phi\over(R+\zeta)\sin\phi}\,{\bf A}_\phi&{1\over(R+\zeta)\sin\phi}{\bf A}_{\theta,\theta}&{1\over(R+\zeta)}{\bf A}_\zeta\\[5pt]
             {\bf 0}&{\bf A}_{\theta,\zeta}&{\bf 0}\\[5pt]
             -{1\over(R+\zeta)}{\bf A}_\phi&{\bf 0}&{1\over(R+\zeta)}{\bf A}_{\zeta,\phi}\\[5pt]
             {\bf 0}&-{1\over(R+\zeta)}{\bf A}_\theta &{1\over(R+\zeta)\sin\phi}{\bf A}_{\zeta,\theta}\\[5pt]
              {\bf 0}& {\bf 0}& {\bf A}_{\zeta,\zeta}
             \end{matrix}\right],\quad\tilde{\bf I} =\left\{\begin{matrix}
                      1 \\[5pt]
                      0 \\[5pt]
                      0 \\[5pt]
                      0 \\[5pt]
                      1 \\[5pt]
                      0 \\[5pt]
                      0 \\[5pt]
                      0 \\[5pt]
                      1
                    \end{matrix}\right\}\cr\cr\cr
{\bf G}_2& ={1\over(R+\zeta)}\left[\begin{matrix}
{\bf A}_{\phi}&{\bf 0}&{\bf 0}\\[5pt]
             {\bf 0}&{\bf 0}&{\bf 0}\\[5pt]
             {\bf 0}&{\bf 0}&{\bf 0}\\[5pt]
             {\bf 0}&{\bf A}_{\theta}&{\bf 0}\\[5pt]
             {\bf 0}&{\bf 0}&{\bf 0}\\[5pt]
             {\bf 0}&{\bf 0}&{\bf 0}\\[5pt]
             {\bf 0}&{\bf 0}&{\bf A}_{\zeta}\\[5pt]
             {\bf 0}&{\bf 0} &{\bf 0}\\[5pt]
              {\bf 0}& {\bf 0}& {\bf 0}
             \end{matrix}\right],\quad
{\bf G}_3 ={1\over(R+\zeta)\sin\phi}\left[\begin{matrix}
             {\bf 0}&{\bf 0}&{\bf 0}\\[5pt]
             {\bf A}_{\phi}&{\bf 0}&{\bf 0}\\[5pt]
             {\bf 0}&{\bf 0}&{\bf 0}\\[5pt]
             {\bf 0}&{\bf 0}&{\bf 0}\\[5pt]
             {\bf 0}&{\bf A}_{\theta}&{\bf 0}\\[5pt]
             {\bf 0}&{\bf 0}&{\bf 0}\\[5pt]
             {\bf 0}&{\bf 0}&{\bf 0}\\[5pt]
             {\bf 0}&{\bf 0} &{\bf A}_{\zeta}\\[5pt]
              {\bf 0}& {\bf 0}& {\bf 0}
             \end{matrix}\right]\cr&&
\end{align}
\subsection{Specialization to Plates}
The simplest case of the shell is the plates, which can be analyzed using rectangular cartesian coordinates $(x,y)$ in the plane of the plate and $z$ coordinate in the normal thickness direction. In this case, we will consider $\eta_1=x$, $\eta_2=y$ and $\zeta = z$. The orthonormal basis vector would be simply the cartesian basis vectors $\hat{\bf \eta}_1=\hat{\bf e}_x,\,\hat{\bf \eta}_2=\hat{\bf e}_y$ and $\hat{\bf n}=\hat{\bf e}_z$ and in this coordinate system the gradient can simply be given as follows:
\begin{align}\label{eq:A53}
{\boldsymbol \nabla}{\bf u}&=u_{x,x} \hat{\bf e}_x\otimes\hat{\bf e}_x+u_{x,y}\hat{\bf e}_x\otimes\hat{\bf e}_y+u_{x,z}\,\hat{\bf e}_x\otimes\hat{\bf e}_z\cr
            &+u_{y,x}\hat{\bf e}_y\otimes\hat{\bf e}_x+u_{y,y}\hat{\bf e}_y\otimes\hat{\bf e}_y+u_{y,z}\,\hat{\bf e}_y\otimes\hat{\bf e}_z\cr
            &+u_{z,x}\hat{\bf e}_z\otimes\hat{\bf e}_x+u_{z,y}\hat{\bf e}_z\otimes\hat{\bf e}_y+u_{x,z}\,\hat{\bf e}_z\otimes\hat{\bf e}_z
\end{align}
The coefficients ${\bf G}_1,\, {\bf G}_2$ and ${\bf G}_3$ used in Eq.~\eqref{eq:23}, for the deformation gradient, become:
\begin{eqnarray}\label{eq:A54}
{\bf G}_1= \left[\begin{matrix}
                   {\bf A}_{x,x} & {\bf 0}& {\bf 0}\\
                   {\bf A}_{x,y} & {\bf 0}& {\bf 0}\\
                   {\bf A}_{x,z} & {\bf 0}& {\bf 0}\\
                   {\bf 0}& {\bf A}_{y,x} & {\bf 0}\\
                   {\bf 0}& {\bf A}_{y,y} & {\bf 0}\\
                   {\bf 0}& {\bf A}_{y,z} & {\bf 0}\\
                   {\bf 0}& {\bf 0} & {\bf A}_{z,x}\\
                   {\bf 0}& {\bf 0} & {\bf A}_{z,y}\\
                   {\bf 0}& {\bf 0} & {\bf A}_{z,z}\\
                 \end{matrix}\right],\
{\bf G}_2 =\left[\begin{matrix}
                       {\bf A}_{x} & {\bf 0} & {\bf 0} \\
                       {\bf 0} & {\bf 0} & {\bf 0} \\
                       {\bf 0} & {\bf 0} & {\bf 0} \\
                       {\bf 0} & {\bf A}_{y} & {\bf 0} \\
                       {\bf 0} & {\bf 0} & {\bf 0} \\
                       {\bf 0} & {\bf 0} & {\bf 0} \\
                       {\bf 0} & {\bf 0} & {\bf A}_{z} \\
                       {\bf 0} & {\bf 0} & {\bf 0} \\
                       {\bf 0} & {\bf 0} & {\bf 0}
                     \end{matrix}\right],\
 {\bf G}_3 =\left[\begin{matrix}
                       {\bf 0}  & {\bf 0} & {\bf 0} \\
                       {\bf A}_{x} & {\bf 0} & {\bf 0} \\
                       {\bf 0} & {\bf 0} & {\bf 0} \\
                       {\bf 0} & {\bf 0}  & {\bf 0} \\
                       {\bf 0} & {\bf A}_{y} & {\bf 0} \\
                       {\bf 0} & {\bf 0} & {\bf 0} \\
                       {\bf 0} & {\bf 0} & {\bf 0}  \\
                       {\bf 0} & {\bf 0} & {\bf A}_{z} \\
                       {\bf 0} & {\bf 0} & {\bf 0}
                     \end{matrix}\right],\
\tilde{\bf I} =\left\{\begin{matrix}
                      1 \\
                      0 \\
                      0 \\
                      0 \\
                      1 \\
                      0 \\
                      0 \\
                      0 \\
                      1
                    \end{matrix}\right\}\cr
\end{eqnarray}
\renewcommand{\theequation}{B.\arabic{equation}}
\renewcommand{\thefigure}{B.\arabic{figure}}
\renewcommand{\thesection}{B}
\setcounter{equation}{0}
\setcounter{figure}{0}
\setcounter{section}{0}
\setcounter{subsection}{0}
\section*{Appendix B: The derivative of the invariants $I_1$, $I_2$ and $J$ with respect to displacement gradient }
The matrices ${\bf G}_0$ and ${\bf G}_{cof}$ used in the derivative of determinant $J$ (see Eqs.~\eqref{eq:26}--\eqref{eq:29}) are defined as follows:
\begin{align}\label{eq:B1}
{\bf G}_{0}&=\left[\ \begin{matrix*}[r]
                       0 & 0 & 0 & 0 & 1 & 0 & 0 & 0 & 1 \\
                       0 & 0 & 0 & -1 & 0 & 0 & 0 & 0 & 0 \\
                       0 & 0 & 0 & 0 & 0 & 0 & -1 & 0 & 0 \\
                       0 & -1 & 0 & 0 & 0 & 0 & 0 & 0 & 0 \\
                       1 & 0 & 0 & 0 & 0 & 0 & 0 & 0 & 1 \\
                       0 & 0 & 0 & 0 & 0 & 0 & 0 & -1 & 0 \\
                       0 & 0 & -1 & 0 & 0 & 0 & 0 & 0 & 0 \\
                       0 & 0 & 0 & 0 & 0 & -1 & 0 & 0 & 0 \\
                       1 & 0 & 0 & 0 & 1 & 0 & 0 & 0 & 0
                     \end{matrix*}\ \right]
\cr
{\bf G}_{cof} &=\left[\begin{matrix}
                        0 & 0 & 0 & 0 & L_{\zeta  \zeta}& -L_{\zeta  \eta_2} & 0 & -L_{\eta_2 \zeta} & L_{\eta_2\eta_2} \\
                        0 & 0 & 0 & -L_{\zeta  \zeta}  & 0 & L_{\zeta  \eta_1} & L_{\eta_2 \zeta}  & 0 & -L_{\eta_2\eta_1} \\
                        0 & 0 & 0 & L_{\zeta  \eta_2} & -L_{\zeta  \eta_1} & 0 & -L_{\eta_2\eta_2} & L_{\eta_2\eta_1} & 0 \\
                        0 & -L_{\zeta  \zeta}  & L_{\zeta  \eta_2} & 0 & 0 & 0 & 0 & L_{\eta_1  \zeta}  & -L_{\eta_1 \eta_2} \\
                        L_{\zeta  \zeta}  & 0 & -L_{\zeta  \eta_1} & 0 & 0 & 0 & -L_{\eta_1  \zeta}  & 0 & L_{\eta_1 \eta_1} \\
                        -L_{\zeta  \eta_2} & L_{\zeta  \eta_1} & 0 & 0 & 0 & 0 & L_{\eta_1 \eta_2} & -L_{\eta_1 \eta_1} & 0 \\
                        0 & L_{\eta_2 \zeta}  & -L_{\eta_2\eta_2} & 0 & -L_{\eta_1  \zeta}  & L_{\eta_1 \eta_2} & 0 & 0 & 0 \\
                        -L_{\eta_2 \zeta}  & 0 & L_{\eta_2\eta_1} & L_{\eta_1  \zeta}  & 0 & -L_{\eta_1 \eta_1} & 0 & 0 & 0 \\
                        L_{\eta_2\eta_2} & -L_{\eta_2\eta_1} & 0 & -L_{\eta_1 \eta_2} & L_{\eta_1 \eta_1} & 0 & 0 & 0 & 0
                      \end{matrix}\right]\cr
\end{align}
where $L_{ij}$ with $i,j = \eta_1,\eta_2,\zeta$ are the components of the displacement gradient tensor ${\bf L}={\boldsymbol\nabla}{\bf u}=L_{ij}\hat{\bf e}_i\otimes\hat{\bf e}_j$.

Further, the derivative of the invariants $I_1$ and $I_2$ with respect to ${\bf L}$ are:
\begin{eqnarray}\label{eq:B2}
{\partial I_1\over \partial{\bf L}} &=& 2({\bf I}+{\bf L})\cr\cr
{\partial I_2\over \partial{\bf L}} &=& 4{\bf F}{\bf C}=4({\bf I}+{\bf L})({\bf I}+{\bf L}+{\bf L}^{\rm T}+{\bf L}^{\rm T}{\bf L})\cr
&=&4({\bf I}+\underbrace{2{\bf L}+{\bf L}^{\rm T}}_{{\bf L}_1}+\underbrace{{\bf L}^{\rm T}{\bf L}+{\bf L}^2+{\bf L}{\bf L}^{\rm T}}_{{\bf L}_2}+\underbrace{{\bf L}{\bf L}^{\rm T}{\bf L}}_{{\bf L}_3})
\end{eqnarray}
where the components of tensors ${\bf L}_1,\, {\bf L}_2$ and ${\bf L}_3$ can be written in the column vector form in a similar fashion as described in Eq.~\eqref{eq:22} for the considered orthonormal coordinate system as
\begin{eqnarray}\label{eq:B3}
\tilde{\bf L}_1={\bf B}_1\tilde{\bf L},\quad \tilde{\bf L}_2={1\over 2}{\bf B}_2\tilde{\bf L},\quad \tilde{\bf L}_3={1\over 3}{\bf B}_3\tilde{\bf L}
\end{eqnarray}
where
{\small
\begin{align}\label{eq:B4}
{\bf B}_1& = \left[\begin{matrix}
                    3 & 0 & 0 & 0 & 0 & 0 & 0 & 0 & 0 \\
                    0 & 2 & 0 & 1 & 0 & 0 & 0 & 0 & 0 \\
                    0 & 0 & 2 & 0 & 0 & 0 & 1 & 0 & 0 \\
                    0 & 1 & 0 & 2 & 0 & 0 & 0 & 0 & 0 \\
                    0 & 0 & 0 & 0 & 3 & 0 & 0 & 0 & 0 \\
                    0 & 0 & 0 & 0 & 0 & 2 & 0 & 1 & 0 \\
                    0 & 0 & 1 & 0 & 0 & 0 & 2 & 0 & 0 \\
                    0 & 0 & 0 & 0 & 0 & 1 & 0 & 2 & 0 \\
                    0 & 0 & 0 & 0 & 0 & 0 & 0 & 0 & 3
                  \end{matrix}\right]_{(9\times9)}
\end{align}
}
and
{\small
\begin{align}\label{eq:B5}
\hskip-100pt{\bf B}_2&=\left[\begin{matrix}
                  6L_{\eta_1 \eta_1}  & 2L_{\eta_1 \eta_2}+L_{\eta_2 \eta_1}  & 2L_{\eta_1 \zeta }+L_{\zeta \eta_1}  & 2L_{\eta_2 \eta_1} +L_{\eta_1 \eta_2} & 0 \\
                  2L_{\eta_1 \eta_2}+L_{\eta_2 \eta_1}  & 2L_{\eta_1 \eta_1} +2L_{\eta_2 \eta_2} &  L_{\zeta \eta_2}+ L_{\eta_2 \zeta} & L_{\eta_1 \eta_1} +L_{\eta_2 \eta_2} & 2L_{\eta_1 \eta_2}+L_{\eta_2 \eta_1}   \\
                  2L_{\eta_1 \zeta}+L_{\zeta \eta_1}  & L_{\eta_2 \zeta}+L_{\zeta \eta_2} & 2L_{\eta_1 \eta_1} +2L_{\zeta \zeta} & L_{\eta_2 \zeta} & 0 \\
                  L_{\eta_1 \eta_2}+2L_{\eta_2 \eta_1}  & L_{\eta_1 \eta_1} +L_{\eta_2 \eta_2} & L_{\eta_2 \zeta} & 2L_{\eta_1 \eta_1} +2L_{\eta_2 \eta_2} & 2L_{\eta_2 \eta_1} +L_{\eta_1 \eta_2}  \\
                  0 & 2L_{\eta_1 \eta_2}+L_{\eta_2 \eta_1}  & 0 & 2L_{\eta_2 \eta_1} +L_{\eta_1 \eta_2} & 6L_{\eta_2 \eta_2} \\
                  0 & L_{\eta_1 \zeta} & L_{\eta_1 \eta_2}+L_{\eta_2 \eta_1}  & L_{\eta_1 \zeta}+L_{\zeta \eta_1}  & 2L_{\eta_2 \zeta}+L_{\zeta \eta_2} \\
                  2L_{\zeta \eta_1} +L_{\eta_1 \zeta} & L_{\zeta \eta_2} & L_{\eta_1 \eta_1} +L_{\zeta\zeta} & L_{\eta_2 \zeta}+L_{\zeta \eta_2} & 0  \\
                  0 & L_{\eta_1 \zeta}+L_{\zeta \eta_1}  & L_{\eta_1 \eta_2} & L_{\zeta \eta_1}  & 2L_{\zeta \eta_2}+L_{\eta_2 \zeta}  \\
                  0 & 0 & 2L_{\eta_1 \zeta}+L_{\zeta \eta_1}  & 0 & 0
                \end{matrix}\right.\cr\cr
                &\hskip100pt\left.\begin{matrix}
                 0 & 2L_{\zeta \eta_1} +L_{\eta_1 \zeta } & 0 & 0 \\
                  L_{\eta_1 \zeta} &  L_{\zeta \eta_2} &  L_{\zeta \eta_1} + L_{\eta_1 \zeta}& 0 \\
                  L_{\eta_2 \eta_1} +L_{\eta_1 \eta_2} & L_{\eta_1 \eta_1} +L_{\zeta \zeta} & L_{\eta_1 \eta_2} & L_{\zeta \eta_1} +2L_{\eta_1 \zeta } \\
                  L_{\zeta \eta_1} +L_{\eta_1 \zeta} & L_{\zeta \eta_2}+L_{\eta_2 \zeta} & L_{\zeta \eta_1}  & 0 \\
                  2L_{\eta_2 \zeta}+L_{\zeta \eta_2} & 0 & 2L_{\zeta \eta_2}+L_{\eta_2 \zeta} & 0 \\
                  2L_{\eta_2 \eta_2}+2L_{\zeta\zeta} & L_{\eta_2 \eta_1}  & L_{\eta_2 \eta_2}+L_{\zeta\zeta} & 2L_{\eta_2 \zeta}+L_{\zeta \eta_2} \\
                  L_{\eta_2 \eta_1}  & 2L_{\eta_1 \eta_1} +2L_{\zeta\zeta} & L_{\eta_2 \eta_1} +L_{\eta_1 \eta_2} & 2L_{\zeta \eta_1} +L_{\eta_1 \zeta} \\
                  L_{\eta_2 \eta_2}+L_{\zeta \zeta}  & L_{\eta_1 \eta_2}+L_{\eta_2 \eta_1}  & 2L_{\eta_2 \eta_2}+2L_{\zeta \zeta} & 2L_{\zeta \eta_2}+L_{\eta_2 \zeta} \\
                  2L_{\eta_2 \zeta}+L_{\zeta \eta_2} & 2L_{\zeta \eta_1} +L_{\eta_1 \zeta} & 2L_{\zeta \eta_2}+L_{\eta_2 \zeta} & 6L_{\zeta \zeta }
                \end{matrix}\right]_{(9\times9)}\cr\cr
\end{align}
}
{\small
\begin{align}\label{eq:B6}
\hskip-50pt{\bf B}_3 &= \left[\ \ \begin{matrix}
                          L_{\eta_1k}L_{\eta_1k}+L_{n\eta_1}L_{n\eta_1}+L_{\eta_1\eta_1}L_{\eta_1\eta_1} & L_{n\eta_2}L_{n\eta_1}+L_{\eta_1\eta_2}L_{\eta_1\eta_1} & L_{n\zeta}L_{n\eta_1}+L_{\eta_1\zeta}L_{\eta_1\eta_1} \\
                          L_{n\eta_1}L_{n\eta_2}+L_{\eta_1\eta_1}L_{\eta_1\eta_2} & L_{n\eta_2}L_{n\eta_2}+L_{\eta_1k}L_{\eta_1k}+L_{\eta_1\eta_2}L_{\eta_1\eta_2} & L_{n\zeta}L_{n\eta_2}+L_{\eta_1\zeta}L_{\eta_1\eta_2} \\
                          L_{n\eta_1}L_{n\zeta}+L_{\eta_1\eta_1}L_{\eta_1\zeta} & L_{n\eta_2}L_{n\zeta}+L_{\eta_1\eta_2}L_{\eta_1\zeta} & L_{\eta_1k}L_{\eta_1k}+L_{n\zeta}L_{n\zeta}+L_{\eta_1\zeta}L_{\eta_1\zeta} \\
                          L_{\eta_2k}L_{\eta_1k}+L_{\eta_2\eta_1}L_{\eta_1\eta_1} & L_{\eta_2\eta_2}L_{\eta_1\eta_1} & L_{\eta_2\zeta}L_{\eta_1\eta_1} \\
                          L_{\eta_2\eta_1}L_{\eta_1\eta_2} & L_{\eta_2k}L_{\eta_1k}+L_{\eta_2\eta_2}L_{\eta_1\eta_2} & L_{\eta_2\zeta}L_{\eta_1\eta_2} \\
                          L_{\eta_2\eta_1}L_{\eta_1\zeta} & L_{\eta_2\eta_2}L_{\eta_1\zeta} & L_{\eta_2k}L_{\eta_1k}+L_{\eta_2\zeta}L_{\eta_1\zeta}\\
                          L_{\zeta k}L_{\eta_1k}+L_{\zeta\eta_1}L_{\eta_1\eta_1} & L_{\zeta\eta_2}L_{\eta_1\eta_1} & L_{\zeta\zeta}L_{\eta_1\eta_1} \\
                          L_{\zeta\eta_1}L_{\eta_1\eta_2} & L_{\zeta k}L_{\eta_1k}+L_{\zeta\eta_2}L_{\eta_1\eta_2} & L_{\zeta\zeta}L_{\eta_1\eta_2} \\
                          L_{\zeta\eta_1}L_{\eta_1\zeta}& L_{\zeta\eta_2}L_{\eta_1\zeta} & L_{\zeta k}L_{\eta_1k}+L_{\zeta\zeta}L_{\eta_1\zeta}
                        \end{matrix}\right.\cr\cr
                        & \left.\begin{matrix}
                          L_{\eta_1k}L_{\eta_2k}+L_{\eta_1\eta_1}L_{\eta_2\eta_1} & L_{\eta_1\eta_2}L_{\eta_2\eta_1} & L_{\eta_1\zeta}L_{\eta_2\eta_1} \\
                          L_{\eta_1\eta_1}L_{\eta_2\eta_2} & L_{\eta_1k}L_{\eta_2k}+L_{\eta_1\eta_2}L_{\eta_2\eta_2} & L_{\eta_1\zeta}L_{\eta_2\eta_2} \\
                          L_{\eta_1\eta_1}L_{\eta_2\zeta} & L_{\eta_1\eta_2}L_{\eta_2\zeta} & L_{\eta_1k}L_{\eta_2k}+L_{\eta_1\zeta}L_{\eta_2\zeta} \\
                          L_{\eta_2k}L_{\eta_2k}+L_{n\eta_1}L_{n\eta_1}+L_{\eta_2\eta_1}L_{\eta_2\eta_1} & L_{n\eta_2}L_{n\eta_1}+L_{\eta_2\eta_2}L_{\eta_2\eta_1} & L_{n\zeta}L_{n\eta_1}+L_{\eta_2\zeta}L_{\eta_2\eta_1} \\
                          L_{n\eta_1}L_{n\eta_2}+L_{\eta_2\eta_1}L_{\eta_2\eta_2} & L_{\eta_2k}L_{\eta_2k}+L_{n\eta_2}L_{n\eta_2}+L_{\eta_2\eta_2}L_{\eta_2\eta_2}& L_{n\zeta}L_{n\eta_2}+L_{\eta_2\zeta}L_{\eta_2\eta_2} \\
                          L_{n\eta_1}L_{n\zeta}+L_{\eta_2\eta_1}L_{\eta_2\zeta} & L_{n\eta_2}L_{n\zeta}+L_{\eta_2\eta_2}L_{\eta_2\zeta} & L_{\eta_2k}L_{\eta_2k}+L_{n\zeta}L_{n\zeta}+L_{\eta_2\zeta}L_{\eta_2\zeta}\\
                          L_{\zeta k}L_{\eta_2k}+L_{\zeta\eta_1}L_{\eta_2\eta_1} & L_{\zeta\eta_2}L_{\eta_2\eta_1} & L_{\zeta\zeta}L_{\eta_2\eta_1} \\
                          L_{\zeta\eta_1}L_{\eta_2\eta_2} & L_{\zeta k}L_{\eta_2k}+L_{\zeta\eta_2}L_{\eta_2\eta_2} & L_{\zeta\zeta}L_{\eta_2\eta_2} \\
                          L_{\zeta\eta_1}L_{\eta_2\zeta} & L_{\zeta\eta_2}L_{\eta_2\zeta} & L_{\zeta k}L_{\eta_2k}+L_{\zeta\zeta}L_{\eta_2\zeta}
                        \end{matrix}\right.\cr\cr
                        &\left.\begin{matrix}
                          L_{\eta_1k}L_{\zeta k}+L_{\eta_1\eta_1}L_{\zeta\eta_1} & L_{\eta_1\eta_2}L_{\zeta\eta_1} & L_{\eta_1\zeta}L_{\zeta\eta_1} \\
                          L_{\eta_1\eta_1}L_{\zeta\eta_2} & L_{\eta_1k}L_{\zeta k}+L_{\eta_1\eta_2}L_{\zeta\eta_2} & L_{\eta_1\zeta}L_{\zeta\eta_2} \\
                          L_{\eta_1\eta_1}L_{\zeta\zeta} & L_{\eta_1\eta_2}L_{\zeta\zeta} & L_{\eta_1k}L_{\zeta k}+L_{\eta_1\zeta}L_{\zeta\zeta} \\
                          L_{\eta_2k}L_{\zeta k}+L_{\eta_2\eta_1}L_{\zeta\eta_1} & L_{\eta_2\eta_2}L_{\zeta\eta_1} & L_{\eta_2\zeta}L_{\zeta\eta_1} \\
                          L_{\eta_2\eta_1}L_{\zeta\eta_2} & L_{\eta_2k}L_{\zeta k}+L_{\eta_2\eta_2}L_{\zeta\eta_2} & L_{\eta_2\zeta}L_{\zeta\eta_2} \\
                          L_{\eta_2\eta_1}L_{\zeta\zeta} & L_{\eta_2\eta_2}L_{\zeta\zeta} & L_{\eta_2k}L_{\zeta k}+L_{\eta_2\zeta}L_{\zeta\zeta}\\
                          L_{\zeta k}L_{\zeta k}+L_{n\eta_1}L_{n\eta_1}+L_{\zeta\eta_1}L_{\zeta\eta_1} & L_{n\eta_2}L_{n\eta_1}+L_{\zeta\eta_2}L_{\zeta\eta_1} & L_{n\zeta}L_{n\eta_1}+L_{\zeta\zeta}L_{\zeta\eta_1} \\
                          L_{n\eta_1}L_{n\eta_2}+L_{\zeta\eta_1}L_{\zeta\eta_2} & L_{\zeta k}L_{\zeta k}+L_{n\eta_2}L_{n\eta_2}+L_{\zeta\eta_2}L_{\zeta\eta_2} & L_{n\zeta}L_{n\eta_2}+L_{\zeta\zeta}L_{\zeta\eta_2} \\
                          L_{n\eta_1}L_{n\zeta}+L_{\zeta\eta_1}L_{\zeta\zeta} & L_{n\eta_2}L_{n\zeta}+L_{\zeta\eta_2}L_{\zeta\zeta} & L_{\zeta k}L_{\zeta k}+L_{n\zeta}L_{n\zeta}+L_{\zeta\zeta}L_{\zeta\zeta}
                        \end{matrix}\ \ \right]_{(9\times9)}
\end{align} }
In the definition of ${\bf B}_3$, the repeated indices $n$ and $k$ imply  summation. Further, the derivative of the invariants $I_1$ and $I_2$ with respect to the column vector form of deformation or displacement gradient can be given as:
\begin{eqnarray}\label{eq:B7}
{\partial I_1\over \partial\tilde{\bf F}}={\partial I_1\over \partial\tilde{\bf L}} &=& 2(\tilde{\bf I}+\tilde{\bf L})\cr\cr
{\partial I_2\over \partial\tilde{\bf F}}={\partial I_2\over \partial\tilde{\bf L}} &=&4\left(\tilde{\bf I}+{\bf B}_1\tilde{\bf L}+{1\over 2}{\bf B}_2\tilde{\bf L}+{1\over 3}{\bf B}_3\tilde{\bf L}\right)
\end{eqnarray}

\bibliographystyle{unsrt}
\bibliography{references}

\begin{thebibliography}{10}

\bibitem{love1888xvi}
Augustus Edward~Hough Love.
\newblock Xvi. the small free vibrations and deformation of a thin elastic
  shell.
\newblock {\em Philosophical Transactions of the Royal Society of London.(A.)},
  179:491--546, 1888.

\bibitem{naghdi1973theory}
Paul~Mansour Naghdi.
\newblock The theory of shells and plates.
\newblock In {\em Linear Theories of Elasticity and Thermoelasticity}, pages
  425--640. Springer, 1973.

\bibitem{ventsel2002thin}
Eduard Ventsel, Theodor Krauthammer, and E~Carrera.
\newblock Thin plates and shells: theory, analysis, and applications.
\newblock {\em Appl. Mech. Rev.}, 55(4):B72--B73, 2002.

\bibitem{gauss1902general}
Carl~Friedrich Gauss.
\newblock {\em General investigations of curved surfaces of 1827 and 1825}.
\newblock Princeton university library, 1902.

\bibitem{green1971theoretical}
AE~Green and W~Zerna.
\newblock Theoretical elasticity,(1968), 1971.

\bibitem{basar1996finite}
Y~Basar and Y~Ding.
\newblock Finite-element analysis of hyperelastic thin shells with large
  strains.
\newblock {\em Computational Mechanics}, 18(3):200--214, 1996.

\bibitem{bacsar1998finite}
Yavuz Ba{\c{s}}ar and Mikhail Itskov.
\newblock Finite element formulation of the ogden material model with
  application to rubber-like shells.
\newblock {\em International Journal for Numerical Methods in Engineering},
  42(7):1279--1305, 1998.

\bibitem{campello2011exact}
EMB Campello, PM~Pimenta, and P~Wriggers.
\newblock An exact conserving algorithm for nonlinear dynamics with rotational
  dofs and general hyperelasticity. part 2: shells.
\newblock {\em Computational Mechanics}, 48(2):195--211, 2011.

\bibitem{kiendl2015isogeometric}
Josef Kiendl, Ming-Chen Hsu, Michael~CH Wu, and Alessandro Reali.
\newblock Isogeometric kirchhoff--love shell formulations for general
  hyperelastic materials.
\newblock {\em Computer Methods in Applied Mechanics and Engineering},
  291:280--303, 2015.

\bibitem{luo2016nonlinear}
Kai Luo, Cheng Liu, Qiang Tian, and Haiyan Hu.
\newblock Nonlinear static and dynamic analysis of hyper-elastic thin shells
  via the absolute nodal coordinate formulation.
\newblock {\em Nonlinear Dynamics}, 85(2):949--971, 2016.

\bibitem{song2016consistent}
Zilong Song and Hui-Hui Dai.
\newblock On a consistent finite-strain shell theory based on 3-d nonlinear
  elasticity.
\newblock {\em International Journal of Solids and Structures}, 97:137--149,
  2016.

\bibitem{li2019consistent}
Yuanyou Li, Hui-Hui Dai, and Jiong Wang.
\newblock On a consistent finite-strain shell theory for incompressible
  hyperelastic materials.
\newblock {\em Mathematics and Mechanics of Solids}, 24(5):1320--1339, 2019.

\bibitem{amabili2020}
M~Amabaili and JN~Reddy.
\newblock The nonlinear, third-order thickness and shear deformation theory for
  statics and dynamics of laminated composite shells.
\newblock {\em Composite Structures}, 244, 2020.

\bibitem{arciniega2007tensor}
RA~Arciniega and JN~Reddy.
\newblock Tensor-based finite element formulation for geometrically nonlinear
  analysis of shell structures.
\newblock {\em Computer Methods in Applied Mechanics and Engineering},
  196(4-6):1048--1073, 2007.

\bibitem{arciniega2006tensor}
Roman~Augusto Arciniega~Aleman.
\newblock {\em On a tensor-based finite element model for the analysis of shell
  structures}.
\newblock PhD thesis, Texas A\&M University, 2006.

\bibitem{amabili2010new}
Marco Amabili and JN~Reddy.
\newblock A new non-linear higher-order shear deformation theory for
  large-amplitude vibrations of laminated doubly curved shells.
\newblock {\em International Journal of Non-Linear Mechanics}, 45(4):409--418,
  2010.

\bibitem{amabili2015non}
Marco Amabili.
\newblock Non-linearities in rotation and thickness deformation in a new
  third-order thickness deformation theory for static and dynamic analysis of
  isotropic and laminated doubly curved shells.
\newblock {\em International Journal of Non-linear Mechanics}, 69:109--128,
  2015.

\bibitem{rivera2016stress}
Miguel~Gutierrez Rivera and JN~Reddy.
\newblock Stress analysis of functionally graded shells using a 7-parameter
  shell element.
\newblock {\em Mechanics Research Communications}, 78:60--70, 2016.

\bibitem{rivera2016new}
Miguel~Gutierrez Rivera, JN~Reddy, and Marco Amabili.
\newblock A new twelve-parameter spectral/hp shell finite element for large
  deformation analysis of composite shells.
\newblock {\em Composite Structures}, 151:183--196, 2016.

\bibitem{rivera2017nonlinear}
Miguel~Gutierrez Rivera and JN~Reddy.
\newblock Nonlinear transient and thermal analysis of functionally graded
  shells using a seven-parameter shell finite element.
\newblock {\em Journal of Modeling in Mechanics and Materials}, 1(2), 2017.

\bibitem{rivera2020continuum}
Miguel~Gutierrez Rivera, JN~Reddy, and Marco Amabili.
\newblock A continuum eight-parameter shell finite element for large
  deformation analysis.
\newblock {\em Mechanics of Advanced Materials and Structures}, 27(7):551--560,
  2020.

\bibitem{amabili2019nonlinear}
M~Amabili, ID~Breslavsky, and JN~Reddy.
\newblock Nonlinear higher-order shell theory for incompressible biological
  hyperelastic materials.
\newblock {\em Computer Methods in Applied Mechanics and Engineering},
  346:841--861, 2019.

\bibitem{reddy2015introduction}
JN~Reddy.
\newblock {\em An Introduction to Nonlinear Finite Element Analysis: with
  applications to heat transfer, fluid mechanics, and solid mechanics}.
\newblock OUP Oxford, 2nd edition, 2015.

\bibitem{knowles1956derivation}
James~K Knowles and Eric Reissner.
\newblock A derivation of the equations of shell theory for general orthogonal
  coordinates.
\newblock {\em Journal of Mathematics and Physics}, 35(1-4):351--358, 1956.

\bibitem{spencer2004continuum}
A~J~M Spencer.
\newblock {\em Continuum mechanics}.
\newblock Courier Corporation, 2004.

\bibitem{darboux1896leccons}
Gaston Darboux.
\newblock {\em Le{\c{c}}ons sur la th{\'e}orie g{\'e}n{\'e}rale des surfaces et
  les applications g{\'e}om{\'e}triques du calcul infinit{\'e}simal: ptie.
  D{\'e}formation infiniment petite et r{\'e}pr{\'e}sentation sph{\'e}rique.
  Par l'auteur. 1896}, volume~4.
\newblock Gauthier-Villars, 1896.

\bibitem{arbind2018curved}
Archana Arbind, AR~Srinivasa, and JN~Reddy.
\newblock A higher-order theory for open and closed curved rods and tubes using
  a novel curvilinear cylindrical coordinate system.
\newblock {\em Journal of Applied Mechanics}, 85(9):091006, 2018.

\bibitem{reddy2018energy}
JN~Reddy.
\newblock {\em Energy Principles and Variational Methods in Applied Mechanics}.
\newblock John Wiley \& Sons, 3rd edition, 2018.

\bibitem{arbind2018Neo-Hookean}
A~Arbind, JN~Reddy, and Srinivasa~A. R.
\newblock A nonlinear finite element analysis for higher-order theory for
  rod/tubes considering incompressible neo-hookean material.
\newblock {\em To be submitted}, 2018.

\bibitem{crisfield1981fast}
Michael~A Crisfield.
\newblock A fast incremental/iterative solution procedure that handles
  “snap-through”.
\newblock In {\em Computational Methods in Nonlinear Structural and Solid
  Mechanics}, pages 55--62. Elsevier, 1981.

\bibitem{flanders1963differential}
Harley Flanders.
\newblock {\em Differential Forms with Applications to the Physical Sciences by
  Harley Flanders}, volume~11.
\newblock Elsevier, 1963.

\end{thebibliography}
\end{document}